\input amstex
\documentstyle{amsppt}
\magnification=1200
\overfullrule=0pt
\NoRunningHeads
%
\def\A{\bold A}
\def\C{\bold C}
\def\F{\bold F}

\def\K{\bold K}
\def\N{\bold N}

\def\R{\bold R}

\def\Q{\bold Q}
\def\Z{\bold Z}

%
\topmatter
\title
Eisenstein series and quantum affine algebras
\endtitle
%
\author
M.M. Kapranov
\endauthor
%
%
%
%
%
%
%
%
%
%
\dedicatory
\enddedicatory
%
%
%
%
%
%
%
%
%
%

%

\endtopmatter
\document

The classical Eisenstein-Maass series is the sum
\TagsOnRight

$$E(z,s) = \frac 1 2 \sum \Sb c,d \in \Z \\(c,d) =1 \endSb \frac
 {y^{s/2}} {|cz+d|^s},\quad z = x+yi \in \bold C, y > 0,  s\in \bold C .\tag 1$$
It converges for
$\text {Re}(s) > 2$
and analytically continues to a meromorphic function in $s$, 
which satisfies the functional equation
$$\zeta^*(s)E(z,s) = \zeta^*(2-s)E(z,2-s), \tag 2$$
where
$\zeta^*(s) = \pi^{-s/2}\Gamma(s/2)\zeta(s)$
is the full zeta function of
$\text {Spec} (\bold Z)$.

This paper started from the observation that the simplest function field 
analog of (1), in which
$\bold Z$
is replaced by the polynomial ring
$\bold F_q[x]$,
is related to the quantum affine algebra
$U_q(\widehat{sl_2})$.
More precisely, the Eisenstein series in this case is naturally a function on
$\text {Bun}(P^1)$,
the set of isomorphism classes of vector bundles on
$P^1_{\bold F_q}$,
and the space $R$ of such functions has a natural algebra structure (Hall
 algebra, as modified by Ringel [R3]), given by the parabolic induction. The analog of
$E(z,s)$
is in fact the product
$E(t_1)E(t_2)$
of two simpler elements of $R$, and the analog of (2) can be written as 
the commutation relation in $R$:

$$E(t_2)E(t_1) = q \frac {\zeta_{P^1}(t_2/t_1)} {\zeta_{P^1} 
(q^{-1}t_2/t_1)} E(t_1)E(t_2), \tag 3$$
where
$\zeta_{P^1}(t) = 1/(1-t)(1-qt)$
is the zeta-function of
$P^1_{\bold F_q}$.

Relations like this are familiar in the theory of vertex operators 
[FLM] [FJ]. In particular, (3), being written in the polynomial form,
 turns out to be identical to one of the relations written by Drinfeld 
[Dr1-2] for his ``loop realization'' of
$U_q(\widehat{sl_2})$.
So $R$ is identified with the natural ``pointwise uppertriangular'' subalgebra
$U_q(\widehat {\bold {n}}^+) \subset U_q(\widehat{sl_2})$.

\vskip .2cm

The main results of this paper (Theorems 3.3, 3.8.4  and 6.7) show 
that for an arbitrary smooth projective curve $X$ over
$\bold F_q$
the structure of the algebra $R$ formed by unramified automorphic 
forms and of its natural extensions, is strikingly similar to the 
structure of the quantum affine algebras
$U_q(\widehat {\bold g})$
for Kac-Moody algebras
$\bold g$
and of their natural subalgebras. The role of the simple roots of
$\bold g$
is played here by unramified cusp eigenforms, and the role of the
 Cartan matrix is played by the Rankin-Selberg convolution $L$-functions.
 Thus $R$ is analogous to
$U_q(\widehat {\bold n}^+)$
with functional equations of general Eisenstein series providing 
commutation relations among the generating functions. 

These results can be viewed as ``continual" analogs of results of
 Ringel [R1-3] and Lusztig [Lu 1-3]
on realizing quantum Kac-Moody algebras by means of Hall algebras 
associated to the category
of representations of a quiver, instead of the category of coherent 
sheaves on a curve.
The two types of categories have very similar properties, 
in particular, they have homological
dimension 1.

\vskip .2cm

The Langlands conjectures bring cusp forms (i.e., simple roots of our 
generalized root system) in correspondence with irreducible 
representations of the geometric fundamental group of the curve $X$. 
This becomes therefore analogous to the well known McKay correspondence 
for finite subgroups in
$SL(2,\bold C)$.
These subgroups are in correspondence with affine Dynkin graphs of type
$A,D,E$,
and for a subgroup $G$ the vertices of the corresponding graph correspond
 to irreducible representations of $G$. Including all coherent sheaves
 (and not just vector bundles) into the general framework of Hall
 algebras amounts to considering Hecke operators acting on unramified 
automorphic forms. It turns out that some natural generating 
functions for these Hecke operators are analogous to Drinfeld's
 generators for the pointwise-Cartan subalgebra
$U_q(\bold{h}[t]) \subset U_q(\widehat{\bold g})$.

\vskip .2cm

One nice outcome of this analogy is that it finally provides some 
explanation of Drinfeld's construction, which since its publication 
10 years ago, was reproduced and used many times, but without questioning
 its nature (i.e., asking
why the formulas have exactly this form and not some other).
From our point of view,  each of these formulas has a clear conceptual
 meaning. Some of them express functional equations for Eisenstein 
series, others the fact that Eisenstein series are eigenfunctions of
 Hecke operators, still other express the commutativity of the unramified 
Hecke algebras and so on.

\vskip .2cm

What seems more important, though, is the conclusion this analogy brings 
about the theory of automorphic forms. Namely, the algebra formed by all 
unramified automorphic forms (on all the
$GL_n$)
and by the Hecke operators, corresponds in our analogy, to just one half of
$U_q(\widehat {\bold g})$,
namely the quantization of the subalgebra
$\bold {n}^+[t,t^{-1}] \oplus \bold {h}[t]$.
It means that one should ``double'' the whole theory of automorphic forms
 by finding the automorphic analog of the missing subalgebra
$\bold {n}^-[t,t^{-1}] \oplus \bold {h}[t^{-1}]$.
In this paper we do this (in our unramified context) by applying Drinfeld's
 quantum double construction [CP] to the Hopf algebra structure given,
 essentially, by taking the constant terms of an automorphic form. Certainly,
 there should be a better and more conceptual definition of this double, and 
the author plans to address this in a future paper. However, the purely 
algebraic identification of the double given in Theorem 6.7, exhibits it as
 a self-contained quantum group-like structure involving all the automorphic 
$L$-functions at once, and one can begin to study its representation theory,
 which is of great interest because of the relation with $L$-functions. For 
instance, one can generalize Frenkel-Jing bosonization construction [FJ] to
 our automorphic case. This will be done elsewhere; here let us point out 
just one thing: the validity of the most non-trivial relation (6.7.5) for 
the bosonization operators is equivalent to the fact that the zeta-function
$\zeta_X(t)$
has only  two ``trivial'' poles at
$t=1, q^{-1}$
(and they produce the two summands in the RHS of (6.7.5)) and all the other 
$L$-functions associated with cusp forms have no poles.

\vskip .2cm

There are two further contexts in which the described approach can be pursued.
 One is that of ``geometric Langlands correspondence'', in which we consider a 

curve $X$ over the complex numbers rather than
$\bold F_q$.
Instead of functions on the discrete set
$\text {Bun}(X/F_q)$
one should consider here the cohomology of the moduli stack of all vector
 bundles (such stack-theoretic cohomology involves, for instance, the
 cohomology of classifying spaces of stabilizers of points, i.e., behaves 
like equivariant cohomology). Results of Grojnowski [Gro] on equivariant
 cohomology of the spaces of
$\bold C$-representations of quivers suggest that in our case we should 
get quantized double-affine algebras, as in [GKV].

\vskip .2cm

Another direction concerns automorphic forms over number fields. The results
 of this paper make it clear that there should be a number-theoretic analog
 of the theory of quantum
 affine algebras in which curves over
$\bold F_q$
are replaced with spectra of the rings in number fields. For instance, the
 most immediate analog of
$U_q(\widehat{sl_2})$
corresponds to (compactified)
$\text {Spec} (\Z)$
instead of
$P^1_{\bold F_q}$
and is generated by values of operator fields
$E^{\pm}(s)$
(generating the analogs of
$U_q(\widehat{\bold n^{\pm}}))$
subject to relations like

$$E^+(s) E^+(t) = \frac{\zeta^*(s-t)} {\zeta^*(s-t+1)} E^+(t)E^+(s)$$
and similar other relations involving the Riemann zeta. The author hopes 
to be able to say more about these questions in the future.

\vskip .3cm

Let us give a brief overview of the contents of the paper. In Section 1 we
 recall the
general framework of Hall-Ringel algebras, including the recent result of
 Green [Gr] on the categorical
description of the comultiplication in the case of homological dimension 1. 
In Green's
formulation one gets, on the Ringel algebra, a structure of a bialgebra in a certain
twisted sense (familiar from [Lu 1-2]). To get a bialgebra in the ordinary sense,
we add the Cartan generators in the standard way. 
Our analysis is very similar to that of [X2].

\vskip .2cm

 Section 2  serves to fix the notation related to unramified automorphic forms
on an algebraic curve $X/\F_q$ (with respect to all the groups $GL_n$).
 The most important
concept for us is that of the Rankin-Selberg convolution $L$-function associated to
two cusp forms. 

\vskip .2cm

In Section 3 we formulate the main results about the Hopf algebra formed
 by the unramified automorphic forms on all the $GL_n$ together with all
 the Hecke operators. We introduce the
appropriate generating functions and state the main result (Theorem 3.3) describing
multiplication and comultiplication of these generating functions in  a 
completely explicit
way. We also prove (Theorem 3.8.4) that the coefficients of our generating functions
generate the subalgebra formed by unramified automorphic forms. 

\vskip .2cm

Section 4 contains proofs of the results stated but left unproved in Section 3. 

\vskip .2cm

In Section 5 we compare the results of Section 3 with the structure of quantum affine
algebras in the ``new realization" of Drinfeld. We observe an analogy between two theories;
we show also  that  in the simplest instances (curve $P^1$, 
 affine algebra $\widehat{sl_2}$)
the analogy becomes the identity. 

\vskip .2cm

Finally in Section 6 we describe Drinfeld's quantum double of the Hopf
 algebra constructed in Section 3.
This seems to be a very important ``semisimple" object naturally appear
ing in the theory
of automorphic forms. 
The main technical tool here is the use of Heisenberg doubles [AF] [ST] 
which are easier
to handle. In particular, we get a very transparent formula for the
 multiplication
in the Heisenberg double of the Ringel algebra in terms of long exact sequences. 
Then we find the relations in the Drinfeld double by using the recent
 work of Kashaev [Kas]
who found an embedding of the Drinfeld double into the tensor product 
of two Heisenberg doubles.

\vskip .3cm

I would like to thank A. Goncharov, G. Harder, Y. Soibelman and Y. Tschinkel
 for useful discussions. In particular, I owe to G. Harder a crucial suggestion 
for the proof of Theorem 3.8.4. This research was partially supported by an NSF 
grant and by
 A.P. Sloan Fellowship.
The paper was written during my stay at Max-Planck-Institut f\"ur Mathematik
 in Bonn, whose hospitality and financial support are gratefully acknowledged. 
I am grateful to Mrs. M. Sarlette for typing the manuscript.

\newpage

\centerline {\bf  \S 1. Hall algebras.}

\vskip 1cm

\noindent  {\bf  (1.1) The Euler form.} We will say that an Abelian category
$\Cal A$
is of finite type, if for any objects
$A,B \in \text{Ob}( \Cal A)$
all the groups
$\text {Ext}^i_{\Cal A}(A,B)$
have finite cardinality and are zero for almost all $i$. If
$\Cal A$
is an abelian category of finite type, and
$A,B \in \text{Ob}( \Cal A)$,
we denote
\TagsOnLeft

$$\langle A,B\rangle = \sqrt{\prod_{i\geq 0} |\text {Ext}^i_{\Cal A}
 (A,B)|^{(-1)^i}} \tag 1.1.1$$
For
$A \in {\Cal A}$
let
$\bar A$
be the class of $A$ in the Grothendieck group
$\Cal K_0\Cal A$.
Clearly, the quantity
$\langle A,B\rangle$
depends only on
$\bar A$
and
$\bar B$
and descends to a bilinear form (called the Euler form) still denoted by

$$\alpha,\beta \mapsto \langle\alpha,\beta\rangle, \quad \Cal K_0\Cal A 
\otimes \Cal K_0\Cal A \rightarrow \bold Q^* \tag 1.1.2$$

\vskip .3cm

\noindent  {\bf  (1.2) Hall and Ringel algebras.}
Let
$\Cal A$
be an Abelian category of finite type. Its Hall algebra
$H(\Cal A)$
is the 
$\C$
-vector space with basis [A] for all isomorphism classes of objects
$A \in \text{Ob}( \Cal A)$.
The multiplication is given by

$$[A] \circ [B] = \sum_{[C]} g^C_{AB} [C] \tag 1.2.1$$
for a fixed object $C$, where
$g^C_{AB}$
is the number of subjects
$A^{\prime} \subset C$
such that
$A^\prime \simeq A$
and
$C/A \simeq B$,
or equivalently, the number of exact sequences

$$0 \rightarrow A \overset \alpha \to \longrightarrow C \overset \beta \to
 \longrightarrow B \rightarrow 0$$
taken modulo the (free) action of
$\text {Aut}(A) \times  \text {Aut}(B)$.
This multiplication is well known to be associative, the coefficient at $C$ in
$A_1 \circ \ldots \circ A_n$
being the number of filtrations of $C$ with quotients
$A_1, \ldots, A_n$.

The modified multiplication

$$[A] * [B] = \langle B,A\rangle \cdot [A] \circ [B] \tag 1.2.2$$
is still associative. We will call the Ringel algebra of
$\Cal A$
and denote
$R(\Cal A)$
the same vector space as
$H(\Cal A)$
but with $*$ as multiplication.

\vskip .3cm

\noindent  {\bf  (1.2.3) Remark.}
It was C.M. Ringel [R3] who first drew attention to the twist (1.2.2). A little
 earlier and independently, G. Lusztig [Lu2-3] considered several twistings by
 bilinear forms, withour specially distinguishing the Euler  form (1.1.1).
 With a certain hindsight, precursors of (1.2.2) can be traced as far back as
 the relabelling of the principal series
representations so as to make the intertwiners to act between representations 
whose weights differ exactly by permutation,
 see, e.g. [GN].

\vskip .3cm

\noindent  {\bf  (1.3) Moduli space point of view.}
Let
$\Cal M(\Cal A) = \text{Ob}(\Cal A)/\text{Iso}$
be the set of isomorphism classes of objects of
$\Cal A$.
The algebras
$H(\Cal A), R(\Cal A)$
can be identified with the space of functions
$f : M(\Cal A) \rightarrow \C$
with finite support, the operations being

$$(f\circ g)(A) = \sum_{A^{\prime} \subset A} 
f(A^{\prime})g(A/A^{\prime}) \qquad (f*g)(A) = \sum_{A^{\prime} \subset A} \langle A/A^{\prime},A^{\prime} \rangle f(A^{\prime}) g(A/A^{\prime}). \tag 1.3.1$$
This point of view makes very natural the ``orbifold'' Hermitian scalar product on
$H(\Cal A)$
and
$R(\Cal A)$:

$$(f,g) =\sum_{A\in \Cal M (A)} \frac{f(A)\overline{g(A)}} {|\text {Aut}
 (A)|} \tag 1.3.2$$
or, equivalently,

$$([A], [B]) = \delta_{[A][B]}/|\text {Aut}(A)|. \tag 1.3.3$$

\vskip .3cm

\noindent  {\bf  (1.4) Green's comultiplication.}
Let
$\Cal A$
be an Abelian category of finite type, satisfying the following additional
 condition: every object of
$\Cal A$
has only finitely many subobjects. Let
$r : R(\Cal A) \rightarrow R(\Cal A) \otimes R(\Cal A)$
be the map, adjoint to the multiplication map
$m : R(\Cal A) \otimes R(\Cal A) \rightarrow R(\Cal A)$
with respect to the scalar product (1.3.2). It has the form

$$r ([A]) = \sum_{A^{\prime} \subset A} \langle A/A^{\prime},A^{\prime}
 \rangle \frac{|\text {Aut}(A^{\prime})|\cdot |\text {Aut}(A/A^{\prime})|}
 {|\text {Aut}(A)|} [A^{\prime}] \otimes [A/A^{\prime}], \tag 1.4.1$$
or, in the functional language (1.3), for a function
$f : M(\Cal A) \rightarrow \C$,
the element 
$r(f)$
is a function on
$\Cal  M (\Cal A) \times \Cal  M (\Cal A)$
given by

$$r(f)(A^{\prime}, A^{\prime\prime}) = \langle A^{\prime\prime},
 A^{\prime}\rangle \sum_{\xi \in \text {Ext}^1(A^{\prime\prime},
 A^{\prime})} f(\text {Cone} (\xi) [-1]) \tag 1.4.2 $$
where
$\text{Cone}(\xi)[-1]$
is the middle term of the extension corresponding to
$\xi$.

For two objects
$A,B \in \Cal A$
set

$$(A|B) = \langle A,B\rangle \cdot \langle B,A\rangle. \tag 1.4.3$$
One easily verifies that the twisted multiplication on
$R(\Cal A) \otimes R(\Cal A)$
given by

$$([A] \otimes [B])([C] \otimes [D]) = (A|B)(([A] * [C]) \otimes
 ([B] * [D])) \tag 1.4.4$$
is associative.

\proclaim {(1.5) Green's theorem}
Suppose that
$\Cal A$
satisfies the  conditions of (1.4), and, in addition, has homological dimension
$\le 1$,
i.e.,
$\text {Ext}^i_{\Cal A} (A,B) =0$
for
$i\ge 2$
and all
$A,B$.
Then
$r: R(\Cal A) \rightarrow R(\Cal A) \otimes R(\Cal A)$
is an algebra homomorphism, if the multiplication on
$R(\Cal A) \otimes R(\Cal A)$
is given by (1.4.4). \endproclaim

Note that because of the twist (1.4.4), the theorem does \underbar{not} mean that
$R(\Cal A)$
is a bialgebra in the ordinary sense; it can be interpreted, however, by saying that
$R(\Cal A)$
is a bialgebra in an appropriate braided monoidal category of
$\Cal K_0\Cal A$
-graded vector spaces.

In [Gr], Green considered only the case when
$\Cal A$ 
consists of finite modules over an
$\bold F_q$
-algebra. The modification to the case of finite modules over any ring
 (of homological dimension 1) is trivial. The case of general
$\Cal A$,
as in (1.5), can be reduced to this by embedding finite pieces of
$\Cal A$
into the categories of modules over appropriate rings, as in Freyd's
 embedding theorem [Fr].

Sources of Green's result can be found in the works of Lusztig [Lu4]
 and Zelevinsky [Ze] and in more classical studies of the functors of
 parabolic induction and restriction in representation theory [BerZ],
 evaluation of constant terms of Eisenstein series [La3]  and so on.

\vskip .3cm

\noindent  {\bf  (1.6) Reformulation.}
As with
$\langle A,B\rangle$,
the quantity
$(A|B)$
depends only on
$\bar A,\bar B \in \Cal K_0\Cal A$,
giving rise to the form

$$(\alpha |\beta) = \langle \alpha,\beta\rangle \cdot \langle \beta,
 \alpha\rangle : \Cal K_0\Cal A \otimes \Cal K_0\Cal A \rightarrow \Q^* \tag 1.6.1$$
Let
$\C[\Cal K_0\Cal A]$
be the group algebra of
$\Cal K_0\Cal A$,
with basis
$K_{\alpha}$,
$\alpha \in \Cal K_0 \Cal A$
and multiplication
$K_{\alpha} K_{\beta} = K_{\alpha +\beta}$.
Let us extend the algebra
$R(\Cal A)$
by adding to it these symbols
$K_{\alpha}$
which we make commute with
$[A] \in R(\Cal A)$
by the rule

$$[A]K_{\beta} = (\bar A|\beta) K_{\beta} [A]. \tag 1.6.2$$
Denote the resulting algebra
$B(\Cal A)$.
So as a vector space
$B(\Cal A) \simeq \C[\Cal K_0\Cal A] \otimes_{\C} R(\Cal A)$,
with
$K_{\alpha} \otimes [A] \mapsto K_{\alpha}A$
establishing the isomorphism.

\vskip .3cm

The next statement was obtained by Xiao [X2].

\proclaim{(1.6.3) Green's theorem (strengthened form)}
In the assumptions of (1.5), the map
$\Delta : B(\Cal A) \rightarrow B(\Cal A) \otimes B(\Cal A)$
given by:

$$\Delta(K_{\alpha}) =K_{\alpha} \otimes K_{\alpha},$$
$$\Delta([A]) = \sum_{A^{\prime} \subset A} \langle A/A^{\prime}, 
A^{\prime}\rangle \frac {|\text {Aut} (A^{\prime})|\cdot |\text {Aut}
 (A/A^{\prime})|} {|\text {Aut} (A)|} [A^{\prime}] \otimes K_{A^{\prime}}
 [A/A^{\prime}] \tag 1.6.4$$
makes
$B(\Cal A)$
into a bialgebra in the ordinary sense, i.e., 
$\Delta$
is a homomorphism of algebras with respect to the standard (untwisted) 
multiplication in
$B(\Cal A) \otimes B(\Cal A)$.
Moreover, $B(\Cal A)$
is a Hopf algebra with respect to the counit
$\epsilon : B(\Cal A) \rightarrow \C$
given by

$$\epsilon (K_{\alpha} [A]) = \left\{ \aligned &1, \quad \text{if} \quad A=0 \\
&0, \quad \text{if} \quad A \ne 0 \endaligned \right. \tag 1.6.5$$
and antipode
$S : B(\Cal A) \rightarrow B(\Cal A)$
given by

$$\gather S(K_{\alpha} [A] = \sum^{\infty}_{n=1} (-1)^n \sum_{A_0 \subset 
\ldots \subset A_n=A} \prod^n_{i=1} \langle A_i/A_{i-1},A_{i-1}\rangle 
\frac{\prod^n_{j=0} |\text {Aut} (A_j/A_{j-1})|} {|\text {Aut} (A)|} \cdot \\
\cdot [A_0] * [A_1/A_0] * \ldots * [A_n/A_{n-1}] \cdot K^{-1}_{\alpha}
 K^{-1}_A \tag 1.6.6 \endgather$$
where
$A_0 \subset \ldots \subset A_n =A$
runs over arbitrary chains of strict
$(A_i \ne A_{i+1})$
inclusions of length $n$.
\endproclaim

\demo{Proof}
The fact that
$\Delta$
is a homomorphism of algebras, follows at once from Theorem 1.5 and from (1.6.2).
 To prove that
$\epsilon$
is a counit, we must show that it is an algebra homomorphism and
 that the compositions

$$(Id \otimes\epsilon)\Delta, \quad (\epsilon \otimes Id)\Delta :
 B(\Cal A) \rightarrow B(\Cal A) \otimes \C = B(\Cal A)$$
are the identity maps. Both these statements are obvious from the 
nature of multiplication in 
$H(\Cal A), R(\Cal A)$
and
$B(\Cal A)$. To prove that
$S$
is an antipode, we must show that the compositions

$$m(S \otimes Id)\Delta, \quad m(Id \otimes S)\Delta : B(\Cal A)
 \rightarrow B(\Cal A)$$
coincide with 
$i \circ \epsilon$
where
$i : \C \rightarrow B(\Cal A)$
is the embedding of the unit, and $m$ is the multiplication in
$B(\Cal A)$.
Let us show this for the first composition, the second one being similar.

From (1.6.4) and (1.6.6), we find that

$$\gather ((S \otimes Id)\Delta)(K_{\alpha}[A]) =\sum_{A^{\prime}
 \subset A} \sum^{\infty}_{n=1} \sum_{A^{\prime}_0 \subset \ldots
 \subset A^{\prime}_n =A^{\prime}} (-1)^n \cdot \langle A/A^{\prime},
 A^{\prime}\rangle \cdot \prod^n_{i=1} \langle A^{\prime}_i/A^{\prime}_{i-1},
 A^{\prime}_{i-1}\rangle \cdot \\
\cdot \frac {|\text {Aut}(A/A^{\prime})|\cdot \prod^n_{j=0}|\text
 {Aut}(A^{\prime}_j/A^{\prime}_{j-1})|} {|\text {Aut} (A)|}
 \cdot [A^{\prime}_0] \ldots [A^{\prime}_n/A^{\prime}_{n-1}]
 K^{-1}_A K^{-1}_{\alpha} \otimes K_{\alpha} K_A[A/A^{\prime}] \endgather$$
where the first sum is over all subobjects
$A^{\prime} \subset A$. 
We can combine the first and third summations together and write the above quantity as

$$\gather \sum^{\infty}_{m=1}(-1)^{m-1} \sum_{A_0 \subset \ldots 
\subset A_m=A} \prod^m_{i=1} \langle A_i/A_{i-1}, A_{i-1}\rangle 
\frac{\prod^m_{j=0} |\text {Aut} (A_j/A_{j-1})|} {|\text {Aut} (A)|} \cdot \\
\cdot [A_0] \ldots [A_{m-1}/A_{m-2}]K^{-1}_{A_{m-1}} K^{-1}_{\alpha}
 \otimes K_{\alpha} K_{A_{m-1}}[A_m/A_{m-1}], \endgather$$
where
$A_0 \subset \ldots \subset A_m =A$
runs over all chains of subobjects of length $m$, in which all the
 inclusions, except, maybe,
$A_{m-1} \subset A_m$,
are strict. Therefore

$$\gather (m(S \otimes Id)\Delta)(K_{\alpha} [A]) = \\
=\sum^{\infty}_{m-1}(-1)^{m-1} \sum_{A_0 \subset \ldots A_m =A}
 \prod^m_{i=1} \langle A_i/A_{i-1}, A_{i-1} \rangle \frac {\prod^m_{j=0}
 |\text {Aut}(A_j/A_{j-1})|} {|\text {Aut} (A)|} \cdot [A_0] 
\ldots [A_m/A_{m-1}]\endgather$$
Now notice that for
$A \ne 0$
each summand in this sum will appear twice: once for a strictly 
increasing filtration
$A_0 \subset \ldots A_m$
and once for the filtration
$A_0 \subset \ldots \subset A_m = A_{m-1}$.
These summands will enter with opposite signs and so will cancel
 each other, and the whole sum will be equal to
$0 = i(\epsilon(K_{\alpha}[A]))$.
If
$A=0$, we get the sum of the empty set of summands, which is equal to
$1 = i(\epsilon (K_{\alpha}))$.
Theorem is proved.
\enddemo

\vskip .3cm

\noindent  {\bf  (1.7) The bilinear form on $B(\Cal A)$.}
Let us extend the Hermitian form (1.3.2-3) on the Ringel algebra
$R(\Cal A)$
to
$B(\Cal A)$
by putting

$$(K_{\alpha}[A], K_{\beta} [B]) = (\alpha|\beta) ([A], [B])  = 
\frac {(\alpha|\beta) \delta_{[A],[B]}} {|\text {Aut} (A)|} \tag 1.7.1$$
In other words, we introduce the form on
$B(\Cal A) = \C[\Cal K_0\Cal A] \otimes_{\C} R(\Cal A)$
to be the tensor product of the old form on
$R(\Cal A)$
and the form on
$\C[\Cal K_0 \Cal A]$
given by
$(K_{\alpha},K_{\beta}) = (\alpha|\beta)$.

\proclaim{(1.7.2) Proposition}
With respect to the form
$( \quad , \quad )$
the multiplication and comultiplication in the Hopf algebra
$B(\Cal A)$
are adjoint to each other.
\endproclaim

\demo{Proof}
In other words, we need to prove the equality

$$(K_{\alpha}[A]K_{\beta}[B], K_{\gamma}[C]) = (K_{\alpha}[A]
 \otimes K_{\beta} [B], \Delta(K_{\gamma} [C])) \tag 1.7.3$$
where on the right stands the Hermitian bilinear form on
$B(\Cal A) \otimes B(\Cal A)$
given by tensoring
$(\quad , \quad )$
with itself. To prove (1.7.3), notice that the left hand side is

$$ (\bar A|\beta)(K_{\alpha+\beta}[A][B], K_{\gamma}[C]) =
 (\bar A|\beta)(\alpha +\beta|\gamma) ([A] [B], [C]) = $$
$$= (\bar A|\beta)(\alpha+\beta| \gamma)([A] \otimes [B], r([C])) = 
 \leqno (1.7.4)$$
$$= (\bar A|\beta)(\alpha+\beta|\gamma) \sum_{C^{\prime} \subset C} 
\langle C/C^{\prime}, C^{\prime}\rangle \frac {|\text {Aut}(C^{\prime})| 
\cdot |\text {Aut}(C/C^{\prime})|} {|\text {Aut} (C)|} \cdot $$
$$\cdot (A,C^{\prime}) \cdot (B,C/C^{\prime}),  $$
while the right hand side of (1.7.3) is

$$ \biggl(K_{\alpha} [A] \otimes K_{\beta} [B], \quad \sum_{C^{\prime}
 \subset C} \langle C/C^{\prime}, C^{\prime}\rangle \frac{|\text {Aut}(C^{\prime})| 
\cdot |\text {Aut}(C/C^{\prime})|} {|\text {Aut} (C)|} \cdot $$
$$\cdot K_{\gamma}[C^{\prime}] \otimes K_{\bar C+\gamma} 
[C^ {\prime\prime}]\biggl) = \leqno (1.7.5)$$
$$= \sum_{C^{\prime} \subset C} \langle C/C^{\prime}, C^{\prime}\rangle
 \frac {|\text {Aut}(C^{\prime})| \cdot |\text {Aut}(C/C^{\prime})|} 
{|\text {Aut}(C)|} (\alpha|\gamma)(\beta |\gamma)(\beta|\bar C^{\prime}) \cdot $$
$$\cdot (A,C^{\prime}) \cdot (B, C/C^{\prime}). $$
Notice now that in order that
$(A,C^{\prime}) \ne 0$,
we should have 
$A \simeq C^{\prime}$,
and under this assumption the corresponding summands in (1.7.4) and 
(1.7.5) coincide. Proposition is proved.
\enddemo

\newpage

\centerline{\bf  \S 2 Background material related to automorphic forms.} 

\vskip 1cm

\noindent  {\bf  (2.1) Notations and conventions.}
Let $X$ be a smooth projective algebraic curve over a finite field
$\bold F_q$.
In this paper we will be interested in the Hall algebra of
$\Cal A = \text {Coh}_X$,
the category of all coherent sheaves on $X$. Let us start by
 introducing some notations and conventions, to be used
 throughout the rest of the paper.

By $g_X$ we denote the genus of $X$. By a point of $X$ we always mean a
 0-dimensional point $x$ (notation:
$x \in X)$.
For such a point $x$ we denote by
$q_x$
the cardinality of
$\bold F_q(x)$,
the residue field of $x$, and by
$\deg (x)$
the degree
$[\bold F_q(x) : \bold F_q]$.
Thus
$q_x = q^{\deg x}$.
By
$\text{Pic} (X)$
we denote the Picard group of  line bundles on $X$ (defined over
$\bold F_q$).
For
$L \in \text{Pic}(X)$
we denote by
$\deg (L) \in \Z$
the degree (first Chern class) of $L$. Thus for 
$x\in X$
we have
$\deg(x) = \deg (\Cal O_X(x))$.
For a vector bundle $V$ on $X$ of rank $n$ we set
$\deg (V) = \deg(\Lambda^nV)$.
The kernel of the degree homomorphism
$\text{Pic}(X) \rightarrow \Z$
is denoted
$\text{Pic}^0(X)$.
It is a finite Abelian group.

As usual, we identify
$\text{Pic}(X)$
with the group of divisors modulo linear equivalence, by associating
 to a divisor
$D =\sum  n_x \cdot x$
the line bundle
$\Cal O(D)$.
Thus
$\deg D=\deg (\Cal O(D))$.

\vskip .3cm

\noindent  {\bf  (2.2) The adelic language. {Aut}omorphic forms.}
Let
$\text {Bun}_n(X)$
be the set of isomorphism classes of rank $n$ vector bundles on $X$ and
$\text {Bun}_{n,d} (X) \subset \text {Bun}_n(X)$
the set of isomorphism classes of bundles of degree $d$.

Let
$k = \bold F_q(X)$
be the field of rational functions on $X$,
$\A$
its ring of adeles and
$\widehat {\Cal O} \subset \A$
the ring of integer adeles. Then

$$\text {Bun}_n(X) \simeq GL_nk\setminus GL_n\A / GL_n \widehat {\Cal O}.
 \tag 2.2.1$$
Let
$\text{AF}_n$
(resp.
$\text{AF}_{n,d})$
be the space of all complex valued functions on
$\text {Bun}_n(X)$
(resp.
$\text {Bun}_{n,d}(X))$.
By (2.2.1) we can regard such functions as
$GL_n \widehat {\Cal O}$
-invariant automorphic forms on
$GL_n\A$.

A function
$f \in \text{AF}_n$
is called a cusp form if for any vector bundles
$V^{\prime} \in \text {Bun}_{n^{\prime}}(X), V^{\prime\prime}
 \in \text {Bun}_{n^{\prime\prime}}(X), n^{\prime} + n^{\prime\prime} =n$,
$n^{\prime}, n^{\prime\prime} > 0$,
we have

$$\sum_{\xi \in \text {Ext}^1(V^{\prime\prime}, V^{\prime})} f(\text{Cone}
 (\xi)[-1]) =0, \tag 2.2.2$$
compare with (1.4.2). This is equivalent to the standard condition

$$\int_{U(k)\setminus U(\bold A)} f(ug) du=0 \tag 2.2.3$$
where we view $f$, by (2.2.1), as a function on
$GL_n\bold A$
and
$U \subset GL_n\A$
is the unipotent radical of a minimal parabolic subgroup.

\vskip .2cm

Let
$\text{AF}^{cusp}_n \subset \text{AF}_n, \text{AF}^{cusp}_{n,d} \subset 
\text{AF}_{n,d}$
be the subspaces formed by cusp forms. The following fact is a well known
 consequence of the reduction theory of Harder [Ha2] [MW].

\proclaim{(2.2.4) Proposition}
Every function from
$\text{AF}^{cusp}_{n,d}$
has finite support. The space
$\text{AF}^{cusp}_{n,d}$
is finite-dimensional.
\endproclaim

\vskip .2cm

\noindent  {\bf  (2.3) Hecke operators.}
Let
$\text {Coh}_{0,X}$
be the category of coherent sheaves on $X$ with 0-dimensional support. By
$\text{Coh}^{\le n}_0(X)$
we denote the set of isomorphism classes of such sheaves
$\Cal F$
which satisfy the additional property
$\dim (\Cal F \otimes \Cal O_x)\le n$
for all
$x \in X$.
We have an identification

$$\text {Coh}^{\le n}_0(X) \quad \simeq \quad GL_n \widehat {\Cal O}
 \setminus GL_n\A \cap \text{Mat}_n(\widehat {\Cal O})/GL_n\widehat {\Cal O},
 \tag 2.3.1$$
which takes
$g \in GL_n\A \cap \text{Mat}_n \widehat {\Cal O}$
into the sheaf
$\text{Coker} \{g : \widehat {\Cal O}^n\rightarrow \widehat {\Cal O}^n\}$.

For
$\Cal F \in \text {Coh}_{0,X}$
and
$n > 0$
we define the operator
$T_{\Cal F} : \text{AF}_n \rightarrow \text{AF}_n$
by

$$(T_{\Cal F} f)(V) = \sum \Sb V^{\prime} \subset V \\ V/V^{\prime} 
\simeq \Cal F \endSb f(V^{\prime}) \tag 2.3.2$$
where the sum is over coherent subsheaves
$V^{\prime}$
in 
$V \in \text {Bun}_n(X)$
with quotient isomorphic to
$\Cal F$.
(Since
$V^{\prime}$
is locally free,
$f(V^{\prime})$
is defined.) Clearly,
$T_{\Cal F} =0$
on
$\text{AF}_n$
unless
$\Cal F \in \text {Coh}^{\le n}_0(X)$.
If
$\Cal F \in \text {Coh}^{\le n}_0(X)$,
one can describe 
$T_{\Cal F}$
in the adelic language as the operator taking a function
$f : GL_nk \setminus GL_n\A /GL_n \widehat {\Cal O} \rightarrow \C$
into
$T_{\Cal F} f$
given by

$$(T_{\Cal F} f)(g) = \int_{h \in GL_n\A} f(gh^{-1}) \bold 1_{\Cal F} (h)dh \tag 2.3.3$$
where
$\bold 1_{\Cal F}$
is the characteristic function of the double coset corresponding to
$\Cal F$
by (2.3.1). For this reason,
$T_{\Cal F}$
is called the Hecke operator.

\proclaim{(2.3.4) Proposition}
The correspondence
$[\Cal F] \mapsto T_{\Cal F}$
makes 
$\text{AF}_n$
into a left module over the Hall algebra
$H(\text {Coh}_{0,X})$,
and
$\text{AF}^{cusp}_n \subset \text{AF}_n$
is a submodule.
\endproclaim

\demo{Proof}
Let $\Cal M$ be the set of isomorphism classes of all coherent sheaves on $X$, and
$\C [\Cal M]$
be the space of all functions
$\Cal M \rightarrow \C$.
This is just the vector space dual to the Hall algebra
$H(\text {Coh}_X)$
and is therefore an
$H(\text {Coh}_X)$
-bimodule. The left action of
$[\Cal F] \in H(\text {Coh}_{0,X})$
on
$\C[\Cal M]$
(dual to its right action on
$H(\text {Coh}_X))$
is given by the formula identical to (2.3.2), but in which
$V, V^{\prime}$
are arbitrary sheaves. If we view 
$\text{AF}_n$
as a subspace in
$\C[\Cal M]$
(consisting of functions vanishing outside
$\text {Bun}_n(X))$,
then it is preserved by this action, so is an
$H(\text {Coh}_{0,X})$
-module as claimed.

To see that
$\text{AF}^{cusp}_n \subset \text{AF}_n$
is a submodule, note that in the adelic language the condition for
$f \in \text{AF}_n$
to the cuspidal involves left shifts of $f$, while the
 Hecke operators involve right shifts.
\hfill $\square$\enddemo

Let
$\text {Coh}_{x,X}$
be the category of coherent sheaves on $X$ supported at $x$ 
(so each such sheaf has the form
$\bigoplus \Cal O_X/I^{\lambda_i}_x$,
where
$I_x \subset \Cal O_X$
is the ideal of $x$). The following facts are well known.

\proclaim{(2.3.5) Proposition}
(a) For
$\Cal F, \Cal G \in \text {Coh}_{0,X}$
we have
$\langle \Cal F,\Cal G \rangle =1$,
so the multiplications
$\circ , *$
in
$H(\text {Coh}_{0,X})$
and
$R(\text {Coh}_{0,X})$
coincide. \newline
(b) $H(\text {Coh}_{0,X}) = \bigotimes_{x \in X} H(\text {Coh}_{x,X})$
(the restricted tensor product in which almost all factors
 in any decomposable tensor are required to be 1). \newline
(c) Each
$H(\text {Coh}_{0,X})$
is a commutative polynomial algebra in either of
 the following two sets of generators:

$$ [\Cal O^{\oplus n} _x], n \ge 1, \leqno (c1)$$
$$  [\Cal O/I^n_x], n \ge 1.\leqno (c2)$$
(d) Let
$\Lambda = \varprojlim \C[z_1, \ldots, z_n]^{S_n}$
be the ring of symmetric functions, and define an isomorphism

$$\text{Ch} : H(\text {Coh}_{x,X}) \rightarrow \Lambda,
 [\Cal O^{\oplus n}_x] \mapsto q^{-n(n-1)/2}_x e_n(z_1,\ldots, z_n)$$
where
$e_n$
is the elementary symmetric function. Then for any integer sequence
$\mu =(\mu_1\ge\mu_2\ge\ldots\ge\mu_r\ge 0)$
the element
$[\bigoplus \Cal O_X/I^{\mu_i}_x]$
will go into
$q^{-\Sigma(i-1)\mu_i} P_{\mu}(z_1, \ldots, z_n; q^{-1}_x)$
where
$P_{\mu}(z_1,\ldots,z_N,t)$
is the Hall-Littlewood polynomial.
\endproclaim

\demo{Proof}
Part (a) is easily obtained by reduction by devissage to the case
$\Cal F = \Cal O_x, \Cal G = \Cal O_y$.
Part (b) is obvious, while (c) and (d) can be found in Macdonald [Mac].
\enddemo

\vskip .2cm

\noindent  {\bf  (2.4) Cusp eigenforms. The sets ${\text{Cusp}}_n$.}
Let
$\widehat{\text{Pic}(X)}$
be the group of all homomorphisms (characters)
$\mu : \text{Pic}(X) \rightarrow \C^*$.
There is an embedding

$$\C^*\hookrightarrow \widehat{\text{Pic}(X)}, t \mapsto t^{\deg} :
 L \mapsto t^{\deg L} \tag 2.4.1$$
whose cokernel is a finite group of characters of
$\text{Pic}^0(X)$. Let us choose representatives
$\mu_1, \ldots, \mu_h$,
one in each coset by the image of (2.4.1), which are unitary, i.e.,
$|\mu_i(L)| =1$
for any
$L \in \text{Pic}(X)$.

For any character
$\mu : \text{Pic}(X) \rightarrow \C^*$
we denote by
$\text{AF}_n(\mu)$
the space of functions (automorphic forms)
$f : \text {Bun}_n(X)\rightarrow \C$
satisfying the property

$$f(V \otimes L) = \mu(L)f(V), \quad \forall L \in 
\text{Pic}(X). \tag 2.4.2$$
Let
$\text{AF}^{cusp}_n(\mu) \subset \text{AF}_n(\mu)$
be the subspace formed by cusp forms. By (2.2.4)
$\dim \text{AF}^{cusp}_n(\mu) < \infty$.
The Hecke operators
$T_{\Cal F}, \Cal F \in \text {Coh}_0(X)$,
preserve
$\text{AF}_n(\mu)$
and
$\text{AF}^{cusp}_n(\mu)$.
For any algebra homomorphism
$\chi : H(\text {Coh}_{0,X}) \rightarrow \C$
denote by
$\text{AF}^{cusp}_n(\mu)_{\chi}$
the corresponding eigenspace, i.e., the space of
$f \in \text{AF}^{cusp}_n(\mu)$
such that

$$T_{\Cal F} f = \chi([\Cal F]) \cdot f, \qquad \forall
 \Cal F \in \text {Coh}_{0,X} \tag 2.4.3$$
By the multiplicity one theorem for
$GL_n$,
see [Sh], the space
$\text{AF}^{cusp}_n(\mu)_{\chi}$
has dimension at most 1. Let
$\frak X_n(\mu)$
be the set of
$\chi$
such that
$\dim \text{AF}^{cusp}_n(\mu)_{\chi} =1$.
If
$\chi \in \frak X_n(\mu)$,
then for any
$x \in X$
we have

$$\chi([\Cal O^{\oplus n}_x]) = \mu(\Cal O_X(-x)). \tag 2.4.4$$
We have a direct sum decomposition

$$\text{AF}^{cusp}_n(\mu) = \bigoplus_{\chi \in \frak X_n(\mu)}
 \text{AF}^{cusp}_n(\mu)_{\chi}. \tag 2.4.5$$
Choose a non-zero vector
$f_{\chi}$
in each summand in (2.4.5) (in virtue of (2.4.4),
$\mu$
is determined by
$\chi$
so it can be dropped from the notation). Let
${\text{Cusp}}_n$
be the set
$\bigcup^h_{i=1} \{f_{\chi}, \chi \in \frak X_n(\mu_i)\}$
where
$\{\mu_1, \ldots , \mu_h\}$
are our unitary representatives (see above). Let also

$$\text{Cusp} = \coprod_{n \ge 1} {\text{Cusp}}_n.$$

\vskip .3cm

\noindent  {\bf  (2.5) Rankin $L$-functions.}
For
$f = f_{\chi} \in {\text{Cusp}}_n$
and
$x \in X$
introduce the numbers
$\lambda_{i,x} (f), i=1, \ldots, n$
defined up to permutation by the condition

$$e_l\bigl(\lambda_{1,x}(f)^{-1}, \ldots, \lambda_{n,x} 
(f)^{-1}\bigr) = q^{l(l-n)/2}_x \chi([\Cal O^l_x]) \tag 2.5.1$$
where 
$e_l$
is the elementary symmetric function. By (2.3.5) (d) we have, for any
$\mu =(\mu_1\ge\ldots\ge\mu_r\ge 0)$:

$$\chi\left(\left[\bigoplus \Cal O/I^{\mu_i}_x\right]\right) = 
q^{-\Sigma(i-1)\mu_i} P_{\mu} \biggl(q^{\frac {n-1} 2}_x \lambda_{1,x} 
(f)^{-1}, \ldots, q^{\frac {n-1} 2}_x \lambda_{n,x} (f)^{-1}; 
q^{-1}_x\biggr) \tag 2.5.2$$
The $L$-function of $f$ is defined by the product

$$L(f,t) = \prod_{x\in X} \prod_{i=1}^n \frac 1 
{1-\lambda_{i,x}(f)t^{\deg (x)}} \tag 2.5.3$$
Our normalization of the
$\lambda_{i,x}(f)$
is chosen so as to make
$L(f,t)$
satisfy the functional equation exchanging $t$ and
$1/qt$
rather than $t$ and
$1/q^nt$,
as is sometimes done in the theory of automorphic forms.

We will be interested, however, in a more general class of
 $L$-functions, which we call Rankin $L$-functions. Given two cusp eigenforms
$f \in {\text{Cusp}}_n, g \in {\text{Cusp}}_m$,
their Rankin $L$-function
$\text{LHom}(f,g,t)$
is defined by the product

$$\text{LHom}(f,g,t) = \prod_{x\in X} \prod^n_{i=1}
 \prod^m_{j=1} \frac 1 {1-\frac{\lambda_{j,x}(g)} {\lambda_{i,x}(f)}
 t^{\deg (x)}} \tag 2.5.4$$
The notation
$\text{LHom}$
is explained as follows. The Langlands correspondence predicts that to any
$f \in {\text{Cusp}}_n$
one can associate a local system
$\Cal L_f$ on $X$. The function
$\text{LHom}(f,g,t)$
is the automorphic counterpart of the $L$-function of the local system
$\underline{\text {Hom}} (\Cal L_f, \Cal L_g)$.

\vskip .2cm

The following result can be found in [JPS] (see [Bu] for a general
 survey of the Rankin-Selberg method).

\proclaim{(2.5.5) Theorem}
The product (2.5.4) converges for
$|t| < q^{-1}$
and defines a rational function in $t$, still denoted
$\text{LHom}(f,g,t)$.
If
$f \ne g$,
then
$\text{LHom}(f,g,t)$
is a polynomial of degree
$mn(2g_X-2)$.
If
$f=g$,
then

$$\text{LHom}(f,g,t) = \frac {P_{f,g} (t)} {(1-t)(1-qt)} \tag 2.5.6$$
where
$P_{f,g}$
is a polynomial of degree
$(2g_X-2)mn+2$.
In any event,
$\text{LHom}(f,g,t)$
satisfies the functional equation

$$\text{LHom}(f,g,  1/{qt}) = \epsilon_{f,g} t^{(2-2g_X)mn}
 \text{LHom}(f,g,t) \tag 2.5.7$$
where
$\epsilon_{f,g} = \prod_{x\in X} \left( \frac{\prod_j \lambda_{i,x}
 (g)} {\prod_i \lambda_{i,x} (f)} \right)^{\text{ord}_x\omega}$
for any rational differential form
$\omega$ 
on $X$ (this product is independent of
$\omega$, see [De]).
\endproclaim

\vskip .2cm

\noindent  {\bf  (2.6) Scalar products and dual Hecke operators.}
For two automorphic forms
$f,g \in \text{AF}_{n,d}$,
of which at least one has finite support, put

$$(f,g)_d = \sum_{V \in \text {Bun}_{n,d}(X)} \frac {f(V)
 \overline{g(v)}} {|\text {Aut} (V)|}. \tag 2.6.1$$
Thus
$(f,g)_d$
is the degree $d$ part of the orbifold scalar product (1.3.2). If
$f,g$
are automorphic forms defined on all
$\text {Bun}_n(X)$,
not just
$\text {Bun}_{n,d}(X)$,
we denote by
$(f,g)_d$
the scalar product of their restrictions to
$\text {Bun}_{n,d}(X)$.
We also write

$$\parallel f \parallel^2_d = (f,f)_d \tag 2.6.2$$
Note that cusp forms have finite support on each
$\text {Bun}_{n,d}(X)$
by (2.2.4), so their scalar products are defined.

We will be interested in the adjoints of the Hecke operators with respect 
to these scalar products. To describe them, we introduce the concept of an 
overbundle.

Let $V$ be a vector bundle on $X$. By an overbundle of $V$ we mean a vector 
bundle $U$ of the same rank as $V$ containing $V$ as a coherent subsheaf. So the sheaf
$U/V$
lies in
$\text {Coh}_{0,X}$.
Its isomorphism class is called the cotype of the overbundle, and the
 isomorphism class of $U$ is called its type. Two overbundles
$U,U^{\prime} \supset V$
are called equivalent, if there is an isomorphism
$U \rightarrow U^{\prime}$
identical on $V$.

\proclaim{(2.6.3) Proposition}
The number of equivalence classes of overbundles of $V$ of type $U$ 
and cotype
$\Cal F$
is equal to

$$g^U_{V\Cal F} \cdot \frac{|\text {Aut}(V)|} {|\text {Aut}(U)|},$$
where
$g^U_{V\Cal F}$
are the same as in (1.2).
\endproclaim

\demo{Proof}
Let
$e^U_{V\Cal F}$
be the number of all exact sequences

$$0 \rightarrow V \rightarrow U \rightarrow \Cal F \rightarrow 0.$$
Then
$g^U_{V\Cal F} = e^U_{V\Cal F}/|\text {Aut}(V)| \cdot |\text {Aut}(\Cal F)|$,
while the number of equivalence classes of overbundles is equal to
$e^U_{V\Cal F}/|\text {Aut}(U)| \cdot |\text {Aut}(\Cal F)|$,
whence the statement. \enddemo

 For
$\Cal F\in \text {Coh}_{0,X}$
we define the dual Hecke operator
$T^\vee_{\Cal F} : \text{AF}_n \rightarrow \text{AF}_n$ 
by

$$(T^\vee_{\Cal F}f)(V)\quad = \quad \sum \Sb U\supset V \\ U/V
 \simeq \Cal F \endSb f(U) 
\quad =\quad \sum_{U \in \text {Bun}_n(X)} g^U_{V\Cal F}\frac 
{|\text {Aut}(V)|} {|\text {Aut}(U)|} f(U) \tag 2.6.4$$
where the first sum is over equivalence classes of overbundles of $V$ of cotype
$\Cal F$,
and the second sum is over all isomorphism classes of bundles. It is clear that
$T_{\Cal F}$
takes
$\text{AF}_{n,d}$
into
$\text{AF}_{n,d+h^0(\Cal F)}$
and
$T^\vee_{\Cal F}$
takes
$\text{AF}_{n,d}$
into
$\text{AF}_{n,d-h^0(\Cal F)}$.
Here
$h^0(\Cal F) = \dim_{F_q} H^0(X,\Cal F)$.
\vskip .2cm

\proclaim{(2.6.5) Proposition}
For each
$d \in \Z$
the operators

$$T_{\Cal F} : \text{AF}_{n,d} \rightarrow \text{AF}_{n,d+h^0(\Cal F)},
\quad  T^\vee_{\Cal F} : \text{AF}_{n,d+h^0(\Cal F)} \rightarrow \text{AF}_{n,d}$$
are adjoint to each other with respect to the scalar product (2.6.1), i.e., for
$f \in \text{AF}_{n,d}$,
$g \in \text{AF}_{n,d+h^0(\Cal F)}$
with finite support, we have

$$(T_{\Cal F}f,g)_{d+h^0(\Cal F)} = (f,T^\vee_{\Cal F}g)_d.$$
\endproclaim

\demo{Proof}
This follows at once from Proposition 2.6.3 and the definition of
 the orbifold scalar product.
\enddemo

Let now $f$ be a cusp eigenform,
$f \in {\text{Cusp}}_n$,
and let
$\chi_f : H(\text {Coh}_{0,X}) \rightarrow \C$
be the algebra homomorphism describing the action of Hecke operators on $f$:

$$\chi_f([\Cal F]) \cdot f = T_{\Cal F}f .\tag 2.6.6$$
Define a new homomorphism
$\chi^\vee_f : H(\text {Coh}_{0,X}) \rightarrow \C$
by

$$\chi^\vee_f([\Cal O^i_x]) = \left\{ \aligned &\chi_f([\Cal O^{n-i}_x]) 
\chi_f([\Cal O^n_x])^{-1}, i \le n \\
&0, i> n \endaligned \right. \tag 2.6.7$$

\proclaim{(2.6.8) Proposition}
In the above assumptions the action of the dual Hecke operators on $f$ 
is given by

$$T^\vee_{\Cal F} f = \chi^\vee_f([\Cal F]) \cdot f.$$
\endproclaim

\demo{Proof}
Let
$V(x), x\in X$,
be the sheaf whose sections are sections of $V$ which are allowed a 
first order pole at $x$. Then each equivalence class of overbundles of $V$ of cotype
$\Cal O^i_x$
can be realized by a unique subsheaf
$U \subset V(x)$
such that
$V \subset U \subset V(x)$.
In other words, such equivalence classes are in bijection with subsheaves in
$V(x)$
of cotype
$\Cal O^{n-i}_x$.
This means that

$$T^\vee_{\Cal O^i_x} = T_{\Cal O^{n-i}_x} T^{-1}_{\Cal O^n_x} \tag 2.6.9$$
and our statement follows from the definitions.\enddemo

Note that it follows from (2.6.7) and (2.5.2) that

$$\chi^\vee_f\left(\left[\oplus \Cal O/I^{\mu_i}_x\right]\right) = 
q^{-\Sigma(i-1)\mu_i}_x P_{\mu} \biggl(q^{\frac {n-1} 2} \lambda_{1,x} (f),
 \ldots, q^{\frac {n-1} 2}_x \lambda_{n,x} (f); q^{-1}_x\biggr). \tag 2.6.10$$

\proclaim{(2.6.11) Proposition}
Let
$f,g \in {\text{Cusp}}_n$.
Then:
\roster
\item"(a)" For each 
$\Cal F \in \text {Coh}_{0,X}$
we have
$\chi^\vee_f([\Cal F]) = \overline{\chi_f([\Cal F])}$.
\item"(b)" The number
$\parallel f \parallel^2_d$
is independent on
$d \in \Z$.
\item"(c)" If 
$f \ne g$,
then
$(f,g)_d =0$
for any $d$.
\endroster
\endproclaim

\demo{Proof}
(a) Fix 
$h \in \Z_+$.
Then for every
$\Cal F$
with
$h^0(\Cal F) =h$,
we have

$$\chi_f([\Cal F])(f,f)_{d+h} = (T_{\Cal F}f,f)_{d+h} = 
(f,T^\vee_{\Cal F}f)_d = 
\overline{\chi^\vee_f([\Cal F])} (f,f)_d.\tag 2.6.12$$
Note that
$(f,f)_d > 0$
and depends only on $d$  modulo $n$,
since
$f(V \otimes L)=\mu(L) f(V)$
for a line bundle $L$, and
$\mu$
is a unitary character. Thus we conclude that for any
$\Cal F$
with
$h^0(\Cal F)$
divisible by $n$ we indeed have the desired equality
$\chi^\vee_f([\Cal F]) = \overline{\chi_f([\Cal F])}$.
However, any two characters of the Hall algebra coinciding an $[\Cal F]$
with
$h^0(\Cal F) \equiv 0(\text{mod} \, n)$
should be equal. \newline

\vskip .2cm
(b) Apply (a) and (2.6.12) with $h$ now being arbitrary. \newline

\vskip .2cm

(c) If
$f \ne g$
then by the multiplicity one theorem for
$GL_n$,
see [Sh], there is
$x \in X$
such that the set of the
$\lambda_{i,x}(f)$
(with multiplicities) is not equal to the set of the
$\lambda_{i,x}(g)$.
It follows that we can find
$\Cal F \in \text {Coh}_{0,x}(X)$
such that
$\chi_f([\Cal F]) \ne \chi_g([\Cal F])$
and, in addition,
$h^0(\Cal F) \equiv 0 (\text{mod} \, n)$.
Thus

$$\align \chi_f([\Cal F])(f,g)_d &=\chi_f([\Cal F])
 (f,g)_{d+h^0(\Cal F)} = (T_{\Cal F}f,g)_{d+h^0(\Cal F)} = \\
&= (f,T^\vee_{\Cal F}g)_d = \chi_g([\Cal F]) (f,g)_d \endalign$$
whence 
$(f,g)_d = 0$.
\enddemo

\vskip .1cm

\noindent  {\bf  (2.6.13) Remark.}
Part (a) of the above proposition means that for each
$x \in X$
the value of any symmetric function on
$\lambda_{1,x}(f)^{-1}, \ldots, \lambda_{n,x}(f)^{-1}$
is equal to its value on
$\overline{\lambda_{1,x}(f)}, \ldots, \overline{\lambda_{n,x}(f)}$,
in other words, that the set of the
$\lambda_{i,x}(f)^{-1}$
is equal to the set of the
$\overline{\lambda_{i,x}(f)}$.
This is not to be confused with the generalized 
Ramanujan-Petersson conjecture [FK] which, in our normalization,
 asserts that
$\lambda_{i,x} (f)^{-1} = \overline{\lambda_{i,x}(f)}$
for each $i$, i.e., that
$|\lambda_{i,x}(f)|=1$.
\vskip .2cm

\noindent  {\bf  (2.6.14) Normalization convention.}
Recall (2.4.5) that the set
${\text{Cusp}}_n$
was obtained by choosing a nonzero vector
$f_{\chi}$
in each 1-dimensional vector space
$\text{AF}^{cusp}_n(\mu)_{\chi}$.
By Proposition 2.6.11 (b) we can choose
$f_{\chi}$
so that
$\parallel f_{\chi} \parallel^2_d =1$
for any $d$. So in the sequel we will always assume that
$\parallel f \parallel^2_d=1$
for any
$f \in {\text{Cusp}}_n, d\in \Z$.

\vskip .3cm

\noindent  {\bf  (2.7) Dualization of bundles and conjugation of forms.}

\proclaim{(2.7.1) Proposition}
For any
$\Cal F\in \text {Coh}_{0,X}$
equivalence classes of overbundles of $V$ of cotype 
$\Cal F$
are in bijection with (locally free) subsheaves in the dual bundle
$V^*$,
of the same cotype
$\Cal F$.
\endproclaim

\demo{Proof}
Let an overbundle $U$ be given. From the short exact sequence

$$0 \rightarrow V \rightarrow U \rightarrow \Cal F \rightarrow 0$$
we get a long exact sequence for
$\underline{\text {Ext}}^{\bullet}(-,\Cal O_X)$,
a part of which has the form

$$0 \rightarrow U^* \rightarrow V^* \rightarrow \underline{\text 
{Ext}}^1(\Cal F,\Cal O_X) \rightarrow 0$$
Since $X$ is a curve and
$\Cal F\in \text {Coh}_{0,X}$,
the sheaf
$\underline{\text {Ext}}^1(\Cal F,\Cal O_X)$
is (non-canonically) isomorphic to
$\Cal F$.
So
$U^*$
is a subsheaf in
$V^*$
of the same cotype
$\Cal F$.
By applying th dualization twice, we find that our correspondence 
is a bijection.
\enddemo

\proclaim {(2.7.2) Corollary}
If
$V,W \in \text {Bun}_n(X)$
and
$\Cal F \in \text {Coh}_{0,X}$,
then

$$g^{V^*}_{W^*\Cal F} = g^W_{V\Cal F} \frac{|\text {Aut}(V)|} 
{|\text {Aut}(W)|}.$$
For an automorphic form
$f \in \text{AF}_n$
define
$f^D \in \text{AF}_n$
by

$$f^D(V) = f(V^*).$$
\endproclaim

\proclaim{(2.7.3) Proposition}
For any
$\Cal F \in \text {Coh}_{0,X}$
and
$f \in \text{AF}_n$
we have

$$T_{\Cal F}(f^D) = (T^\vee_{\Cal F} f)^D.$$
\endproclaim

\demo{Proof}
This is an immediate consequence of (2.7.1)
\enddemo

\proclaim{(2.7.4) Corollary}
If
$f \in {\text{Cusp}}_n$,
then there is
$\epsilon \in \{\pm 1\}$
such that
$f(V^*) = \epsilon \overline{f(V)}$
for any
$V \in \text {Bun}_n(X)$.
\endproclaim

\demo{Proof}
From (2.7.3), both
$f^D$
and
$\bar f$
are eigenforms of the Hecke algebra with the same character
$\chi^\vee_f = \bar\chi_f$.
So by the multiplicity one theorem
$f^D =\epsilon \bar f$
for some constant
$\epsilon \ne 0$.
Since the dualization and conjugation are involutive,
$\epsilon^2=1$.
\enddemo

\newpage

\centerline {\bf  \S 3. The Hopf algebra of automorphic forms.} 

\vskip 1cm

\noindent  {\bf  (3.1) The setup.}
We keep the notation of (2.1), and are going to apply the general 
formalism of Section 1 to the Abelian category
$\Cal A = \text {Coh}_X$.
We denote
$H=H(\Cal A), R=R(\Cal A)$
and
$B=B(\Cal A)$
its Hall, Ringel and extended Ringel algebra (see Section 1).

The Grothendieck group
$\Cal K_0\Cal A = \Cal K_0X$
is identified with
$\Z \oplus \text{Pic}(X)$
via the map

$$\biggl(n,\,\,\, \sum_{x\in X}m_x \cdot x\biggr)\quad \mapsto 
\quad n \cdot \bar \Cal O_X + \sum_{x \in X} m_x \cdot \bar \Cal O_x \tag 3.1.1$$
where we view
$\text{Pic}(X)$
as the quotient of the group of divisors modulo principal divisors.

We denote the generator
$K_{\bar \Cal O_X} \in \C[\Cal K_0\Cal A]$
simply by $K$ and denote
$K_{\bar \Cal O_x}, x\in X$
simply by 
$c_x$.
For a divisor
$D = \sum m_x\cdot x$
we set
$c_D = \prod c^{m_x}_x$.
By the above,
$c_D =1$
for principal $D$. So for a line bundle $L$ on $X$ there is a
 well-defined element
$c_L = c_D$
where $D$ is any divisor such that
$L \simeq \Cal O_X(D)$.

Note that each
$\bar \Cal O_x$
lies in the kernel of the bilinear form
$(\alpha |\beta)$
on
$\Cal K_0(X)$.
Thus each
$c_x, c_L$
is a central element in $B$. We have there fore a character

$$c: \text{Pic}(X) \rightarrow B^*, \quad L \mapsto c_L \tag 3.1.2$$
of
$\text{Pic}(X)$
with values in the multiplicative group of $B$. As for the
 generator $K$, for any vector bundle $V$ on $X$ we have, by
 Riemann-Roch theorem

$$[V]K=q^{rk(V)(1-g_X)}K[V]. \tag 3.1.3$$
Any coherent sheaf
$\Cal F$
on $X$ can be written as a direct sum
$\Cal F = \Cal F_{tors} \oplus \Cal F_{lf}$,
where
$\Cal F_{lf}$
is locally free and
$\Cal F_{tors}$
is a torsion sheaf (i.e., has 0-dimensional support). The isomorphism
 classes of
$\Cal F_{tors}, \Cal F_{lf}$
depend on
$\Cal F$
only. If
$\Cal F_{lf} \ne 0$,
there are infinitely many subsheaves in
$\Cal F$,
so the formula (1.6.3) for the comultiplicaton in $B$ produces an 
infinite sum, i.e., an element of a certain completion of
$B \otimes B$. More precisely, let
$B \widehat \otimes B$
be the space of possibly infinite sums
$\sum b^{\prime}_i \otimes b^{\prime\prime}_i$
where
$$b^{\prime}_i = [\Cal F^{\prime}_i] \kappa^{\prime}_i, \quad 
b^{\prime\prime}_i = [\Cal F^{\prime\prime}_i]\kappa^{\prime\prime}_i, 
\quad \Cal F^{\prime}_i, \Cal F^{\prime\prime}_i \in 
\text{Coh}(X), \kappa^{\prime}_i, \kappa^{\prime\prime}_i \in
 \C[\Cal K_0 X],$$
satisfying the following condition:

\vskip .2cm

(3.1.4)
For each 
$d\in \Z$
the number of $i$ such that
$\deg \Cal F^{\prime}_{i,lf} =d$,
is finite, and for
$d \gg 0$
this number is 0.

\vskip .2cm

The following fact is easily proved by applying the main lemma of Green
 [Gr] plus the fact that the number of 
coherent subsheaves of given degree in a vector bundle is finite.

\proclaim{(3.1.5) Proposition}
$B \widehat \otimes B$
is an algebra, and
$\Delta : B \rightarrow B \widehat \otimes B$
given by (1.6.3), is a homomorphism of algebras. \endproclaim

So we shall say that $B$ is a topological Hopf algebra.

\vskip .3cm

\noindent  {\bf  (3.2) Generating functions associated to cusp forms.}
Let
$f \in {\text{Cusp}}_n$
be a cusp eigenform on
$\text {Bun}_n(X)$,
and
$\chi =\chi_f : H(\text {Coh}_{0,X}) \rightarrow \C$
be the algebra homomorphism giving the action of Hecke operators on $f$:

$$T_{\Cal F}f = \chi ([\Cal F]) \cdot f \tag 3.2.1$$
Consider the following formal power series with coefficients in the Ringel
 algebra $R$:

$$E_f(t) = \sum_{V\in \text {Bun}_n(X)} f(V) t^{\deg (V)} [V]\quad 
 \in\quad  R[[t,t^{-1}]], \tag 3.2.2$$
the sum over all isomorphism classes of rank $n$ vector bundles on 
$X$. Note that the coefficients at each power of $t$ in
$E_f(t)$
is a finite sum by (2.2.4).

More generally, for any quasi-character
$\mu : \text{Pic}(X) \rightarrow B^*$
taking values in the multiplicative group of the center of $B$ we can 
form the series

$$E_f(\mu t) =\sum_{V \in \text {Bun}_nX} f(V) \mu(\det V)t^{\deg V}[V]
 \quad \in \quad B[[t,t^{-1}]]. \tag 3.2.2'$$
The notation becomes unambigious once we agree to identify $t$ itself
 with the quasicharacter
$L \mapsto t^{\deg L}$
of
$\text{Pic}(X)$.
If $\mu$
takes values in
$\C^* \subset B^*$,
then
$\mu(L) = \mu_i(L) \lambda^{\deg L}$
for some
$\lambda \in \C^*$
and
$\mu_i \in \{\mu_1, \ldots, \mu_h\}$,
our set of unitary representatives (2.4), so
$E_f(\mu t) = E_{f^{\prime}} (\lambda t)$
for some
$f^{\prime} \in {\text{Cusp}}_n$
and we don't get anything new. However, taking
$\mu = c$,
the character defined by (3.1.2), we get new elements, to be used
 later in the formulas for comultiplication.

Let also

$$\psi_f(t) = \sum_{\Cal F\in \text {Coh}_0(X)} \bar\chi_f([\Cal F])
 t^{h^0(\Cal F)} |\text {Aut} (\Cal F)| \cdot [\Cal F] 
\quad \in \quad R[[t]], \tag 3.2.3$$
where the sum is over isomorphism classes of all sheaves with 
0-dimensional support and
$h^0(\Cal F) =\dim H^0(X,\Cal F)$. 
Thus
$\psi_f$
is a generating function for Hecke operators. If $W$ is any rank $n$
 vector bundle such that
$f(W) \ne 0$,
we have, by (2.7.4):

$$\psi_f(t) = \frac{t^{\deg W}} {q^{\deg W} f(W^*)} \sum \Sb V \subset W
 \\ rk(V) =n \endSb f(V^*) t^{-\deg V} |\text {Aut} (W/V)|
 \cdot [W/V] \tag 3.2.4$$
where the sum is over all subsheaves in $W$ of full rank $n$, i.e., over 
``effective matrix divisors''.
As with
$E_f(t)$, we will use the series
$\psi_f(\mu t) \in B[[t]]$
for any central quasicharacter
$\mu : \text{Pic}(X) \rightarrow B^*$.
It is given by

$$\Psi_f(\mu t) =\sum_{\Cal F \in \text {Coh}_0(X)} \bar \chi_f([\Cal F]) 
t^{h^0(\Cal F)} \mu (\bar \Cal F)|\text {Aut} (\Cal F)|
 \cdot [\Cal F] \tag 3.2.5$$
where
$\bar \Cal F$
is the class of
$\Cal F$
in
$\Cal K_0(X)$
which lies in the subgroup
$\text{Pic}(X) \subset \Cal K_0(X)$,
see (3.1.1). The space
$R[[t,t^{-1}]]$
of series
$a(t) =\sum^{+\infty}_{i=-\infty} r_it^i, r_i \in R$,
infinite in both directions, is not a ring, but we do have a well-defined
 multiplication

$$R[[t_1,t^{-1}_1]] \otimes R[[t_2,t^{-1}_2]] \rightarrow R[[t^{\pm 1}_1,
 t^{\pm 1}_2]],
\quad  a(t_1) \otimes b(t_2) \mapsto a(t_1)b(t_2).$$
Also, the use of generating functions is well-suited to the study of the
 topological Hopf algebra
$B \supset R$.
More precisely, we have the following theorem which is the main result of 
this section.

\proclaim{(3.3) Theorem}
Let
$f \in {\text{Cusp}}_n, g \in {\text{Cusp}}_m$
be the two cusp eigenforms. Then: \newline
(a) For each
$\Cal F \in Coh(X)$
the coefficient at the basis vector
$[\Cal F] \in R$
in each of the products
$$E_f(t_1)*E_g(t_2), \quad E_f(t_1)*\psi_g(t_2), \quad
\psi_f(t_1)*E_g(t_2) \in R[[t^{\pm 1}_1, t^{\pm1}_2]]$$
is a power series in
$t_1, t_2$
which converges for
$|t_1| \gg |t_2|$
to a rational function. \newline
(b) These rational functons satisfy the following relations:

$$E_f(t_1)*E_g(t_2) = q^{mn(1-g_X)} \frac {\text{LHom}(f,g,t_2/t_1)} 
{\text{LHom}(f,g,t_2/qt_1)} E_g(t_2)*E_f(t_1) \tag 3.3.1$$

$$E_f(t_1) *\psi_g(t_2) = \frac{\text{LHom}(f,g,q^{\frac m 2} t_2/t_1)} 
{\text{LHom}(f,g,q^{\frac m 2 -1}t_2/t_1)} \psi_g(t_2)*E_f(t_1) \tag 3.3.2$$
(c) In the topological Hopf algebra
$B \supset R$
we have the identities:

$$\Delta\psi_f(t) = \psi_f(t \otimes c)(1 \otimes \psi_f(t)) \tag 3.3.3$$

$$\Delta E_f(t) = 1 \otimes E_f(t) + E_f(t \otimes c)(1 \otimes
 K^n\psi_f(q^{-\frac n 2} t)) \tag 3.3.4$$
where
$E_f(t \otimes c) = \sum_V f(V)t^{\deg (V)} [V]\otimes c_{\det (V)}$.

$$\epsilon (\psi_f(t))=1, \quad \epsilon (E_f(t)) =0, \quad \epsilon(c_L) =1, 
\quad \epsilon (K) =1. \tag 3.3.5$$

$$S(\psi_f(t)) = \psi_f(c^{-1}t)^{-1}, \quad S(E_f(t)) =-E_f(c^{-1}t)\psi_f
 (q^{-\frac n 2} t)^{-1} K^{-n}. \tag 3.3.6$$
\endproclaim

In this section we will prove only the equality (3.3.1), relegating the rest 
to Section 4.

\vskip .3cm

\noindent  {\bf  (3.4) Eisenstein series. Proof of (3.3.1).}
Let
$f \in {\text{Cusp}}_n, g \in {\text{Cusp}}_m$.
Let us write, for the product in the Hall algebra

$$E_f(t_1) \circ E_g(t_2) =\sum_{V\in \text {Bun}_{n+m} (X)} \Cal
 E_V(f,g,t_1,t_2)[V]. \tag 3.4.1$$
Then

$$\Cal E_V(f,g,t_1,t_2 ) =\sum \Sb V^{\prime} \subset V \\ \text {rk}(V) =n
 \endSb f(V^{\prime})g(V/V^{\prime}) t^{\deg(V^{\prime})}_1 
t_2^{\deg(V/V^{\prime})} \tag 3.4.2$$
where
$V^{\prime}$
runs over all subbundles (i.e., subsheaves which are locally
 direct summands) in $V$ of rank $n$. This is nothing but the 
unramified Eisenstein series associated to cusp forms
$f,g$,
see [Ha1] [Mor] [MW]. We recall the following result.

\proclaim{(3.4.3) Proposition}
(a) For any
$V \in \text {Bun}_{n+m} (X)$
the series (3.4.2) converges for
$|t_1| \gg |t_2|$
to a rational function. \newline
(b) These rational function satisfy the functional equations

$$\Cal E_V(f,g,t_1,t_2) = q^{mn(1-g_X)} \frac {\text{LHom}(f,g,qt_2/t_1)}
 {\text{LHom}(f,g,t_2/t_1)}
 \Cal E_V(g,f,q^n t_2, q^{-m}t_1)$$

(c) The poles of the rational function
$\Cal E_V(f,g,t_1,t_2)$
are precisely the poles of the function

$$\frac{\text{LHom}(f,g,qt_2/t_1)} {\text{LHom}(f,g,t_2/t_1)}$$
and the orders of poles are the same.
\endproclaim

\demo{Proof}
(a) See, for instance, [Ha1], [Mor] or [MW], Prop. IV.1.12 . \newline

\vskip .2cm

(b) Let us briefly recall  the general framework of functional equations of
 Eisenstein series [MW] and show what it yields in our particular case. For
 a space $S$ let
$\text{Fun}(S)$
denote the space of locally constant functions on $S$ (we ignore the growth 
conditions in this formal reminder). Let
$\Xi_{n,m} \subset  GL_{n+m} \A$
be the subgroup

$$\Xi_{n,m} = \pmatrix GL_n(k) & \text{Mat}_{n,m}(\A) \\ 0 & GL_m(k) \endpmatrix$$
The Eisenstein series construction defines a map

$$\gather \text{Eis}_{n,m}: \text{Fun}(\Xi_{n,m}\setminus GL_{n+m} (\A)) \rightarrow \text{Fun}(GL_{n+m}(k) \setminus GL_{n+m} (\A)), \\
(\text{Eis} \, \varphi)(g) =\sum_{\gamma \in P_{n,m}(k)\setminus GL_{n+m}(k)} 
\varphi(\gamma g),\quad 
 P_{n,m} = \pmatrix GL_n & * \\ 0 & GL_m \endpmatrix \endgather$$
(when converges). In general,
$\text{Eis}_{n,m}$
is defined by analytic continuation over auxiliary parameters. The functional
 equation, formally, has the form

$$\text{Eis}_{n,m}(\varphi) = \text{Eis}_{m,n} (M\varphi)$$
where

$$M : \text{Fun}(\Xi_{n,m} \setminus GL_{n+m} \A) \rightarrow \text{Fun}(\Xi_{m,n} 
\setminus GL_{n+m} \A)$$
is defined by

$$(M\varphi)(g) = \int_{Z\in \text{Mat}_{n,m}(\A)} \varphi \left( \pmatrix Z & 1_m 
\\ 1_n & 0 \endpmatrix g\right) dZ,$$
where
$dZ = \prod dz_{ij}$
and
$\int_{\A/k} dz_{ij} =1$.
Our case is obtained from here as follows. For
$$f \in \text{AF}_n = \text{Fun}(GL_n k \setminus GL_n \A/GL_n \widehat {\Cal O}),\quad 
 g \in \text{AF}_m$$
let
$f \odot g \in Fun(\Xi_{n,m} \setminus GL_{n+m} \A/GL_{n+m} \widehat {\Cal O})$
be the function defined uniquely (in virtue of the Iwasawa decomposition) by

$$(f \odot g) \pmatrix A & 0 \\ 0 & B \endpmatrix = f(A) g(B),\quad  A \in GL_n \A,
 B\in GL_m \A.$$
For
$A \in GL_n(\A)$
let
$\deg(A) = \sum_{x\in X} \deg(x) \cdot \text{ord}_x(a)$.
Then for
$f \in {\text{Cusp}}_n, g \in {\text{Cusp}}_m$ we have

$$\Cal E_V(f,g,t_1, t_2) = \text{Eis}_{n,m} ((ft^{\deg}_1) \odot (gt^{\deg}_2))$$
so (3.4.3) (b) follows from the identification

$$M((ft^{\deg}_1) \odot (gt^{\deg}_2)) = q^{mn(1-g_X)} \frac
 {\text{LHom}(f,g,qt_2/t_1)} {\text{LHom}(f,gt_2/t_1)} (g \cdot (q^mt_2)^{\deg})
 \odot (f \cdot (q^{-m} t_1)^{\deg}) \tag 3.4.4$$
This can be established by using the fact that the representation of
$GL_n\A$
corresponding to $f$ is, at every
$x \in X$,
a principal series representation, and same for $g$. Thus the operator 
$M$ which has the form
$\bigotimes_{x\in X} M_x$
can be calculated by splitting each
$M_x$
into $mn$ 1-dimensional intertwiners, each evaluated by integration over
$\widehat k_x$.
Each such intertwiner contributes a factor

$$\frac {1-\lambda_{j,x}(g)\lambda_{i,x}(f)^{-1}
 q_x t^{\deg (x)}_2/t^{\deg (x)} _1} 
{1-\lambda_{j,x} (g) \lambda_{i,x} (f)^{-1} t^{\deg (x)}_2/t^{\deg (x)}_1} 
\int_{\widehat {\Cal O}_x} dz$$
which, being all multiplied, give precisely (3.4.4), 
once we recall, that
$\int_{\widehat {\Cal O}} dz = q^{1-g_X}$
if the Haar measure $dz$
on
$\A$
is normalized by
$\int_{\A/k} dz =1$.
 \newline

\vskip .2cm

(c) This follows from (3.4.4) and the general fact about Eisenstein 
series of relative rank 1([MW], $\S$ \, IV.3.10, Remark) which says
 that the singularities of Eisenstein series are in this case precisely 
the singularities of the intertwiner $M$.
\enddemo

Proposition 3.4.3 is proved.

$   $\newline
Now, part (a) of (3.4.3) can be written as an equality in the 
Hall algebra of the category
$\text{Coh}(X):$

$$E_f(t_1)\circ E_g(t_2) = q^{mn(1-g_X)} \frac {\text{LHom}(f,g,qt_2/t_1)}
 {\text{LHom}(f,g,t_2/t_1)} E_g(q^nt_2) \circ E_f(q^{-m}t_1) \tag 3.4.5$$
By using the definition of the Ringel product $*$ and the Riemann-Roch
 theorem, we find

$$E_f(t_1) * E_g(t_2) = q^{mn(1-g_X)} E_f(q^{\frac m 2}t_1) 
\circ E_g(q^{-\frac n 2} t_2)$$
whence the validity of (3.3.1).

\vskip .3cm

\noindent  {\bf  (3.5) Algebraic relations in $B$.}
For
$f \in {\text{Cusp}}_n$
we define elements
$E_{f,d}, d \in \Z$
and
$a_{f,d}, d \in \N$,
by

$$E_f(t) = \sum_{d\in \Z} E_{f,d}t^d;\quad  a_f(t) := \log \psi_f(t)
= \sum^{\infty}_{d=1} a_{f,d} t^d \tag 3.5.1$$
Thus
$E_{f,d}, a_{f,d}$
are some (finite) elements of $B$. We now proceed to find some algebraic 
relations among them by using the relations (3.3.1-2) for generating
 functions. For two different cusp forms
$f \in {\text{Cusp}}_n, g \in {\text{Cusp}}_m$
let

$$Q_{f,g} (t_1,t_2) = t_1^{(2g_X-2)mn} \text{LHom}(f,g,t_2/t_1)
 \tag 3.5.2$$
be the homogeneization of their Rankin $L$-function, and

$$Q_f(t_1,t_2) = t^{(2g_X-2)n^2+2}_1 P_f(t_2/t_1), \quad \text{where}
\quad  \text{LHom} (f,f,t) = \frac {P_f(t)} {(1-t)(1-qt)},\tag 3.5.3$$
be the homogeneization of the numerator of
$\text{LHom}(f,f,t)$.
The relations (3.3.1) can be written in the polynomial form:

$$(t_1-qt_2)Q_f(qt_1,t_2)E_f(t_1) * E_f(t_2) =q^{1-mn(1-g_X)}
 (qt_1-t_2)Q_f(t_1,t_2) \cdot E_f(t_2) * E_f(t_1), \tag 3.5.4$$

$$Q_{f,g} (qt_1, t_2) E_f(t_1) * E_g(t_2) = q^{-mn(1-g_X)} Q_{f,g}
 (t_1,t_2) E_g(t_2) * E_f(t_1), f\ne g \tag 3.5.5$$
where the equality is understood in the same sense as in Theorem 3.3:
 as the equality of rational functions constituting the coefficients of
 the LHS and the RHS at any given $[V]$,
$V \in \text {Bun}_{n+m} (X)$.

\proclaim{(3.5.6) Theorem}
(a) For each
$i,j \in \Z$
comparing the coefficients at
$t^i_1t^j_2$
in both sides of (3.5.4) or (3.5.5) gives a valid relation among 
the elements
$E_{f,d}, E_{g,d} \in R \subset B$. \newline
(b) If we write

$$\text{LHom}(f,g,t) = \prod_i (1-\alpha_i(f,g)t)^{\nu_i}, \tag 3.5.7$$
then we have an equality in $B$:

$$[a_{g,d}, E_{f,l}] = \frac 1 d \left (\sum_i \nu_i
 \alpha_i(f,g)^d(q^d-1)q^{(n+m-\frac 3 2)}\right) E_{f,l+d}. \tag 3.5.8$$
\endproclaim

\demo{Proof}
(a) Consider, say, the equality (3.5.4), and let some
$V \in Bun_{n+m}(X)$
be fixed. The coefficients at $V$ of the two sides of (3.5.4)
 are power series in
$t_1, t_2$
converging to the same rational function, denote it
$\varphi_V(t_1,t_2)$,
but in different regions:
$|t_1| \gg |t_2|$
for the left hand side and
$|t_1| \ll |t_2|$
for the right hand side. So the coefficient at
$t^i_1 t^j_2$
in the left hand side is
$\int_{|t_1| =R,|t_2|=r} \varphi_V(t_1,t_2)t^{-i-1}_1 t^{-j-1}_2 
dt, dt_2$,
$R \gg r$,
while the coefficient at
$t^i_1 t^j_2$
in the right hand side is similar integral, but taken over the torus
$|t_1|=r, |t_2| =R$.
However, Proposition 3.4.3 (c) shows that
$\varphi_V(t_1,t_2)$
has no singularities and thus is a Laurent polynomial in
$t_1,t_2$. So the two integrals coincide, and comparing coefficients
 indeed gives a valid relation. \newline

\vskip .2cm

(b) From (3.3.2) we find

$$[a_g(t_2), E_f(t_1)] = \log \frac {\text{LHom}(f,g,q^{n+m- \frac 3 2
} t_2/t_1)} {\text{LHom} (f,g,q^{n+m-\frac 1 2}(t_2/t_1)} E_f(t_1) \tag 3.5.9$$
Moreover, the series
$a_g(t_2)$
going only in one direction, the coefficients at each $[V]$ in
 both sides of (3.5.9) are Laurent polynomials in
$t_1,t_2$,
so we can proceed to comparing coefficients at each
$t^i_1 t^j_2$.
This is done by applying the formula
$\log \, (1-z) = -\sum_{d\ge 1} z^d/d$
to (3.5.9) and (3.5.7) and yields the claimed answer (3.5.8).
\enddemo

\vskip .2cm

\noindent  {\bf  (3.6) The Hermitian form on $B(\text {Coh}_X)$.}
We now proceed to describe the values of the Hermitian form (1.7.1)
 on our generating functions
$E_f(t), \psi_f(t)$ 
(or
$a_f(t) = \log \, \psi_f(t))$.
This will automatically give us scalar products of any products of
 the generating functions because of the identities

$$(xy,z) =(x\otimes y,\Delta (z)), \quad
 (x,yz) = (\Delta (x), y \otimes z)\tag 3.6.1$$
expressing the fact that the multiplication and the comultiplication in
$B(\text {Coh}_X)$
are conjugate to each other. 

In the following we will make use of the formal power series

$$\delta (z) = \sum^{+\infty}_{n=-\infty} z^n \tag 3.6.2$$
representing (in the sense of distribution theory) the Dirac
$\delta$-function at
$z=1$.

\proclaim{(3.6.3) Proposition}
The scalar products of the generating functions
$E_f(t), a_f(t) \in B(\text {Coh}_X)$
are given by

$$(E_f(t_1), a_g(t_2)) =0, \tag 3.6.4$$

$$(E_f(t_1), E_g(t_2)) = \delta_{f,g} \cdot \delta (t_1\bar t_2), 
f \in {\text{Cusp}}_n, g \in {\text{Cusp}}_m.\tag 3.6.5$$

$$(a_f(t_1), a_g(t_2)) = \log \frac {\text{LHom} (g,f,q^{\frac {n+m} 2 -1} 
t_1\bar t_2)} {\text{LHom}
 (g,f,q^{\frac{n+m} 2} t_1\bar t_2)} \tag 3.6.6$$

$$(c_L,x) =0, \quad  \forall x \in B(\text {Coh}_X) \tag 3.6.7$$

$$(K^i, K^j) = q^{ij(1-g_X)}, \quad (K^i, E_f(t)) = (K^i, a_f(t)) =0 
\tag 3.6.8$$
\endproclaim

The proof will be given in $\S$ \, 4.

\vskip .3cm

\noindent  {\bf  (3.7) More general scalar products.}
A general element of the subalgebra in
$R(\text {Coh}_X)$
generated by the coefficients
$E_{f,d}$,
has the form

$$E(\bold f, \varphi) = \int_{|t_i| =a_i} E_{f_1} (t_1) \ldots E_{f_r} 
\varphi(t_1, \ldots, t_r) \prod \frac {dt_i} {t_i}, \quad a_1 \gg \ldots
 \gg a_r \tag 3.7.1$$
where
$\bold f = (f_1, \ldots, f_r), f_i \in \text{Cusp}$,
is a sequence of cusp forms, and 
$\varphi$
is a Laurent polynomial. In the standard terminology of the theory of
 automorphic forms, such elements are called theta-series [Go] or
 pseudo-Eisenstein series [MW]. The well-known formula [La3] [MW]
 for the scalar product of two pseudo-Eisensetin seris can be easily 
deduced from Theorem 3.3, Proposition 3.6.5 and general properties of
 Hopf algebras. Let us recall this formula in our notation.

Let
$\bar B(\text {Bun}_X)$
be the quotient algebra of
$B(\text {Coh}_X)$,
obtained by putting each [$\Cal F$],
where
$\Cal F$
is not a vector bundle, to be equal to 0, and each
$c_L$
to be equal to 1, and let

$$p_{\text{Bun}} : B(\text {Coh}_X) \rightarrow \bar B(\text {Bun}_X)
 \tag 3.7.2$$
be the natural projection. Notice that if
$b \in B(\text {Coh}_X)$
is a linear combination of $[V]$ with $V$ being vector bundles, then
 for any
$a \in B(\text {Coh}_X)$
the scalar product
$(b,a)$
depends only on
$p_{\text {Bun}} (a)$.

For any Hopf algebra $A$ with comultiplication
$\Delta$
we denote
$\Delta^{(r)} : A \rightarrow A^{\otimes r}$
the 
$(r-1)$
fold iteration of
$\Delta$,
i.e.,
$\Delta^{(2)} = \Delta$
and
$\Delta^{(r+1)}= Id_{A^{\otimes (r-1)}} \otimes \Delta$.
If a scalar product
$( \quad , \quad)$
on $A$ satisfies the identity (3.6.1), then by iterating these
 identities we find, in particular, that

$$(x,y_1 \ldots y_r) = (\Delta^{(r)}(x), y_1 \otimes \ldots 
\otimes y_r) \tag 3.7.3$$
For a sequence
$\bold f = (f_1, \ldots, f_r), f_i \in {\text{Cusp}}_{n_i}$
of cusp forms and a permutation
$\sigma \in S_r$
let

$$M^{\bold f}_{\sigma} (t_1, \ldots, t_r) = \prod \Sb i<j \\ \sigma(i) > 
\sigma(j) \endSb q^{n_in_j(1-g_X)} \frac{\text{LHom}(f_i,f_j, t_j/t_i)}
 {\text{LHom}(f_i, f_j, t_j/q t_i)}, \tag 3.7.4$$
so that the functional equation for Eisenstein series yields

$$E_{f_1} (t_1) \ldots E_{f_r} (t_r) = M^{\bold f}_{\sigma} (t_1, \ldots,
 t_r) E_{f_{\sigma(1)}} (t_{\sigma(1)}) \ldots E_{f_{\sigma(r)}}
 (t_{\sigma(r)}). \tag 3.7.5$$
The classical formula of Langlands for the constant term of a (pseudo) 
Eisenstein series has, in our notation, the form

$$ p^{\otimes r^{\prime}}_{\text {Bun}} \Delta^{r^{\prime}} (E_{f_1} (t_1) 
\ldots E_{f_r} (t_r)) = \leqno (3.7.6)$$
$$= \delta_{rr^{\prime}} \sum_{\sigma \in S_r} M^{\bold f}_{\sigma}
 (t_1, \ldots, t_r) K^nE_{f_{\sigma^{-1}(1)}} (t_{\sigma^{-1}(1)})
 \otimes \ldots \otimes K^nE_{f_{\sigma^{-1}(r)}} (t_{\sigma^{-1}(r)}),
 \, n = \Sigma n_i. $$
This formula follows at once from (3.3.2) and (3.3.4). By using
 (3.6.7) and (3.7.3), we deduce that for
$\bold g = (g_1, \ldots, g_{r^{\prime}}), g_j \in {\text{Cusp}}_{m_j}$,
we have

$$ \biggl( E_{f_1}(t_1) \ldots E_{f_r} (t_r), \,\,\,  E_{g_1}
 (t^{\prime}_1) \ldots E_{g_{r^{\prime}}} (t^{\prime}_{r^{\prime}})
\biggl) = \leqno (3.7.7)$$
$$= \delta_{rr^{\prime}} \sum\Sb \sigma \in S_r \\ g_j =
 f_{\sigma^{-1}(j)}, \forall j \endSb M^{\bold  f}_{\sigma}
 (t_1, \ldots, t_r) \prod^r_{j=1} \delta (t_{\sigma^{-1}(j)}
 \cdot \bar t^{\prime}_j). $$
Then by integrating (3.7.7) against a Laurent polynomial
$\psi(t^{\prime}_1 \ldots t^{\prime}_r)$
 we find

$$ \biggl(E_{f_1}(t_1) \ldots E_{f_r} (t_r), \,\,\, 
E(\bold g, \psi)\biggl ) = \leqno (3.7.8)$$
$$= \sum\Sb \sigma \in S_r \\ g_j = f_{\sigma^{-1}(j)} 
\endSb M^{\bold f}_{\sigma}(t_1, \ldots, t_r)
\cdot \bar \psi (\bar t^{-1}_{\sigma^{-1}(1)}, \ldots, 
\bar t^{-1}_{\sigma^{-1}(r)}) $$

\vskip .2cm

\noindent  {\bf  (3.8) The algebras $\tilde \Cal B$ and 
$\Cal B$.}
Let
$\Cal B \subset B(\text {Coh}_X)$
be the subalgebra generated by
$K, c_L$
and the coefficients 
$E_{f,d}, a_{f,d}, f\in \text{Cusp}$.
One would like to have a complete description of
$\Cal B$
by generators and relations. To address this problem, let us introduce the algebra
$\tilde \Cal B$
generated by formal symbols
$\tilde E_{f,d} , \tilde a_{f,d}, \tilde K, \tilde c_L$
which are subject only to the relations that the
$\tilde c_L$
are central, that

$$ \tilde c_L \tilde c_M = \tilde c_{L \otimes M},\quad \tilde K 
\tilde c_L = \tilde c_L \tilde K, \quad \tilde K \tilde E_{f,d} =
 q^{-n(1-g)} \tilde E_{f,d} K, \leqno (3.8.1)$$
$$\tilde K \tilde a_{f,d} = \tilde a_{f,d} \tilde K, \quad f \in 
{\text{Cusp}}_n, $$
plus the relations obtained from (3.4.4), (3.4.5), (3.4.8) by
 replacing $E$ with
$\tilde E$
and $a$ with
$\tilde a$.
So we have a natural surjection

$$\pi_{\Cal B} : \tilde \Cal B \rightarrow \Cal B , \tag 3.8.2$$
and generating functions
$$\tilde E_f(t) = \sum_{d\in \Z} \tilde E_{f,d}t^d,\quad \tilde a_{f}
 (t) = \sum_{d \ge 1} \tilde a_{f,d} t^d,
\quad \tilde \psi_f(t) = \exp(\tilde a_f(t))$$
satifying (3.3.1) and (3.4.2). The problem of explicit 
determination of the ideal
$\text{Ker}(\pi_{\Cal B})$
seems very interesting. Basically,  elements of $\text{Ker}(\pi_{\Cal B})$
 are some subtle relations between residues
of Eisenstein series.
 As we will see in $\S$ \, 5,  these elements 
should be thought of analogs of Serre relations in quantum affine
 algebras
which makes it plausible that one can give a completely explicit
 description of the generators of the ideal. Here
we will use an approach similar to one used by G. Lusztig [Lu1]
 for ordinary quantum groups.

\vskip .1cm

Denote by
$\tilde\Cal R$
(resp.,
$\tilde \Cal R_{\text {Bun}}, \tilde \Cal R_0)$
the subalgebra in
$\tilde \Cal B$
generated only by the
$\tilde E_{f,d}$
and
$\tilde a_{f,d}$
(resp. only by
$\tilde E_{f,d}$,
only by
$\tilde a_{f,d}$)
and by
$\Cal R$
(resp.
$\Cal R_{\text {Bun}}, R_0)$
the image of
$\tilde \Cal R$
(resp.
$\tilde \Cal R_{\text {Bun}}, \tilde \Cal R_0)$
under
$\pi_{\Cal B}$.

Note that formulas (3.3.4-6) (modified by putting tildes over
 generating functions) make
$\tilde \Cal B$
into a topological Hopf algebra, and
$\pi_{\Cal B}$
into a Hopf homomorphism. Note also that formulas (3.6.6-10) define
 a Hermitian scalar product on
$\tilde \Cal B$,
with respect to which the multiplication and comultiplication are
 conjugate, i.e., (3.6.1) holds. Let
$N \subset \tilde \Cal B$
be the kernel of this scalar product. Since
$(\pi_{\Cal B} (x), \pi_{\Cal B}(y)) = (x,y)$
for any
$x,y \in \tilde \Cal B$,
we have
$\text{Ker} (\pi_{\Cal B}) \subset N$.

\proclaim{(3.8.3) Proposition}
$\text{Ker}  (\pi_{\Cal B}) = \C[\Cal K_0X] \cdot (N \cap \tilde
 \Cal R)$.
\endproclaim

\demo{Proof}
This follows from the fact that the form
$( \quad , \quad )$
on
$\Cal R \subset R(\text {Coh}_X)$,
being the restriction of a positive definite Hermitian form, is 
itself positive definite and hence non-degenerate. Thus
$(\text{Ker} \, \pi_{\Cal B}) \cap \tilde \Cal R = N \cap \tilde
 \Cal R$,
and we should only use the fact that both
$\tilde \Cal B$
and
$\Cal B$
have the form
$\tilde \Cal B = \C[\Cal K_0X] \otimes \tilde \Cal R, \Cal B =
 \C[\Cal K_0X] \otimes \Cal R$.
\enddemo

This proposition shows that the algebra
$\Cal B$
can be defined entirely in terms of the $L$-function data made explicit 
in Theorem 3.3 and Proposition 3.6.5. 
We want to finish this section by discussing the relationship between
$\Cal B$, the ``free'' algebra
$\tilde \Cal B$
and the bigger algebra
$B(\text {Coh}_X)$
a little bit more closely.

Let
$R(\text {Bun}_X) \subset R(\text {Coh}_X)$
be the subalgebra generated by elements $[V]$ where $V$ is a vector
 bundle. By restricting
$\pi_{\Cal B}$
to the two subalgebras below and composing it with natural embeddings, 
we get homomorphisms

$$\pi_0 : \tilde \Cal R_0 \rightarrow R(\text {Coh}_{0,X}),\quad 
 \pi_{\text {Bun}} : \tilde \Cal R_{\text {Bun}} \rightarrow R(\text {Bun}_X).$$

\proclaim{(3.8.4) Theorem}
The homomorphism 
$\pi_{\text {Bun}}$
is surjective, i.e.
$\Cal R_{\text {Bun}} = R(\text {Bun}_X)$. \endproclaim

In other words, any element $[V]$,
$V \in \text {Bun}_n(X)$
can be expressed as a polynomial in the
$E_{f,d}$,
i.e., in the form (3.7.1). This statement is similar to the
 spectral decomposition
theorem [MW] but is different since here we consider the coefficients of the
Laurent expansion of the rational functions $E_{f_1}(t_1)...E_{f_r}(t_r)$
in the domain of convergence of the series defining these functions
(and these coefficients are automorphic forms with finite support on each 
$\text{Bun}_{n,d}(X)$), while spectral decomposition theorem deals, in our
notation, with the coefficients of Laurent expansions of the same functions
but on the unit torus $|t_i| = 1$ (and these coefficients are, in general,
only square-integrable).

We will prove (3.8.4) a little later. As for the homomorphism
$\pi_0$,
we would like to state the following conjecture.

\proclaim {(3.8.5) Conjecture}
The map
$\pi_0$
is injective, i.e., the (commuting) elements
$a_{f,d}, f \in \text{Cusp}, d\ge 1$,
are algebraically independent over
$\C$. \endproclaim

Note that since 
$\langle \Cal F, \Cal G \rangle = 1$
 for any
$\Cal F, \Cal G \in \text {Coh}_{0,X}$,
the formula (1.4.1) makes
$R(\text {Coh}_{0,X}) = H(\text {Coh}_{0,X})$
into a Hopf algebra (with comultiplication denoted by $r$).
 This is just the tensor product over all
$x \in X$
of the Hopf algebras studied by Zelevinsky [Ze]. With respect 
to $r$ every
$\psi_f(t)$
is group-like:
$r(\psi_f(t)) = \psi_f(t) \otimes \psi_f(t)$,
so
$a_f(t)$
is primitive:
$r (a_f(t)) = a_f(t) \otimes 1 +1 \otimes a_f(t)$,
and the same is true for each
$a_{f,d}$.
Now, primitive elements in a commutative and cocommutative Hopf
$\C$
-algebra are algebraically dependent if and only if they 
are linearly dependent. 
By using (3.6.6), we get the following reformulation of
 the conjecture. Let
${\text{Cusp}}_{\le N} = \coprod_{i\le N} {\text{Cusp}}_i$.

\proclaim {(3.8.6) Reformulation}
Let $N$ be fixed, and
$\epsilon$
be small enough. Consider the self-adjoint matrix integral
 operator on the circle
$|t| =\epsilon$
whose matrix indices run over the set
${\text{Cusp}}_{\le N}$
and whose kernel is given by

$$\K_{f,g} (t_1,t_2) = \log \frac{\text{LHom} (g,f,q^{\frac {m+n}
 2 -1} t_1\bar t_2)} {\text{LHom} (g,f,q^{\frac{m+n} 2} t_1 \bar t_2)},
 \quad f\in {\text{Cusp}}_n, g\in {\text{Cusp}}_m, m,n \le N$$
Then this operator is strictly positive definite.
\endproclaim

This conjectural property can be regarded as a certain strengthening
 of the multiplicity one theorem for
$GL_n$.

\vskip .3cm

\noindent  {\bf  (3.9) Proof of Theorem 3.8.4.} We will proceed 
in three steps.

\vskip .2cm

\demo{ \underbar {Step 1: Use of the reduction theory}} We recall
 some basic facts
about stable vector bundles, see [HN], [Stu] for more details.
For any vector bundle
$V$ on $X$ let $\mu(V) = \deg(V)/\text{rk}(V)$ denote its ``slope". 
A bundle $V$ is
called semistable, if for any subbundle $V'\i V$ we have $\mu(V')
 \leq \mu(V/V')$. 
An arbitrary bundle $V$ possesses a canonical Harder-Narasimhan 
filtration
$$V_\bullet = (V_1 \i ... \i V_r = V)$$
with the properties that  for each $i$ the quotient
$\text{gr}_i(V) = V_i/V_{i-1}$ is semistable and $\mu(\text{gr}_i(V)) >
\mu(\text{gr}_{i+1}(V))$. Thus $V$ is semistable if and only if 
its Harder-Narasimhan
filtration consists of only one layer. 

Let $\lambda\in {\bold Z}$ be a positive integer. Let us say that
 a vector bundle $V$
is $\lambda$-unstable, if it is not semistable and for at least one 
$i$
we have $\mu(\text{gr}_i(V)) + \lambda < \mu(\text{gr}_{i-1}(V))$.
 Denote
by $\text{Bun}_{n,d}^{>\lambda}$ the set of isomorphism classes of
 $\lambda$-unstable
bundles on $X$ of rank $n$ and degree $d$, and by
 $\text{Bun}_{n,d}^{\leq\lambda} = 
\text{Bun}_{n,d}(X) - \text{Bun}_{n,d}^{>\lambda}$ its complement.

\proclaim{(3.9.1) Lemma} (a) For each $\lambda > 0$ the set
 $\text{Bun}_{n,d}^{\leq\lambda}$ is finite. \newline
(b) For fixed $n,d$ we can find $\lambda > 0$ such that whenever
 $V\in \text{Bun}_{n,d}^{>\lambda}$ and $i$ is such that
 $\mu(\text{gr}_i(V)) + \lambda < \mu(\text{gr}_{i-1}(V))$, 
then in the Hall algebra we have the equality
$[V] = [V_i] \circ [V/V_i]$. \endproclaim

\demo{Proof of (3.9.1)} (a) This is an easy consequence of the
 reduction theory,
 see, e.g., [Stu]. 

\vskip .1cm

(b) If $E$ is any vector bundle on $X$, then
there is $\delta\in {\bold Z}$ such that for any line bundle on
 $X$ of degree
$\geq \delta$ we have $H^1(X, E\otimes L) = 0$, $H^0(X, (E\otimes L)^*) =
 0$
(Serre's theorem). Moreover, if we have an algebraic family of 
bundles parametrized
by a scheme $S$ of finite type, we can find $\delta$ good for
 all the bundles in
the family.

We apply this to the ``universal" situation. 
Let $\Cal{B}un_{n,d}(X)$ be the moduli space of semistable vector
bundles over $X$ of rank $n$ and degree $d$. It is a projective
algebraic variety defined over $\bold F_q$. Tensoring with a line bundle
of degree 1 defines an isomorphism $\Cal Bun_{n,d}(X) \rightarrow 
\Cal Bun_{n, n+d}(X)$. We conclude therefore the following:

\proclaim{(3.9.2) Lemma} For each $m_1, m_2$ there exists a number
$\delta_{m_1, m_2} > 0$ with the following property: If $W_i$, $i=1,2$,
is a semistable
bundle on $X$ of rank $m_i$ and $\mu(W_1) +\delta_{m_1, m_2} < \mu(W_2)$,
then
$$\text{Ext}^1_X(W_2, W_1) = \text{Hom}_X(W_1, W_2) = 0.$$
\endproclaim

To return to the proof of (3.9.1), take $\lambda$ greater than all 
the $\delta_{m_1, m_2}$, $m_i \leq n$. If $V$ satisfies the condition of (3.9.1),
then by (3.9.2) we have
$$\text{Ext}^1(\text{gr}_j(V), \text{gr}_k(V)) = \text{Hom}(\text{gr}_k(V),
\text{gr}_j(V)) = 0, \quad j> i, k\leq i$$
from which we deduce that
$$\text{Ext}^1(V/V_i, V_i) = \text{Hom}(V_i, V/V_i) = 0.$$
The vanishing of $\text{Ext}^1$ means that $V\simeq V_i \oplus (V/V_i)$
and the vanishing of Hom means that there exists a unique subbundle
in $V_i \oplus (V/V_i)$ of type $V_i$ and cotype $V/V_i$. Using the vanishing
of $\text{Ext}^1$ one more time, we find that $[V_i]\circ [V/V_i]
= [V_i \oplus (V/V_i)] = [V]$. Lemma (3.9.1) is proved.

\enddemo\enddemo

\demo{\underbar{Step 2: Representation of a function as an infinite
sum of Eisenstein coefficients}}
Let
$R(\text {Bun}_{n,X})$
be the subspace in
$R(\text {Coh}_X)$
spanned by basis vectors $[V]$ where $V$ is a rank $n$ vector bundle.
 We are going to use the formula for the
scalar product of two Eisenstein series in order to represent an element
$h\in R(\text {Bun}_{n,X})$ as a possibly infinite linear combination
of coefficients of Eisenstein series.

Fix $n > 0$,
fix some real numbers
$a_1 \gg \ldots \gg a_n > 0$
(enough to take
$a_i > qa_{i+1})$
and a total order
$\le$
on the set
${\text{Cusp}}_{\le n} = \coprod_{i \le n} {\text{Cusp}}_i$.
If
$f \in {\text{Cusp}}_i$,
write
$n(f) =i$.
Denote also
$\bold{Cusp}_n$
the set of sequences
$\bold f = (f_1 \le \ldots \le f_r), f_i \in {\text{Cusp}}_{\le n},
 \Sigma n(f_i) =n$.
We will also write
$E_{\bold f} (t)$
for
$E_{f_1}(t_1) \ldots E_{f_r} (t_r), t=(t_1, \ldots, t_r)$.

\proclaim{(3.9.3) Lemma}
If
$h \in R(\text {Bun}_{n,X})$
is orthogonal to each element
$E_{\bold f}(t), f=(f_1, \ldots, f_r) \in \bold{Cusp}_n, |t_i| =a_i$,
then
$h=0$.
\endproclaim

\demo{Proof}
The assumptions imply that
$(E_{\bold f}(t),h) =0$
identically as a rational function in
$t_1 \ldots t_r$.
Thus from the functional equations (3.7.5) we find that for any permutation
$\sigma \in S_r$
and
$|t_i| = a_i$
the product
$E_{f_{\sigma (1)}}(t_1) \ldots E_{f_{\sigma (r)}}(t_r)$
is orthogonal to $h$ as well. This means that $h$ is orthogonal to
 all pseudo-Eisenstein series and thus
$h=0$
by [MW], Th. II.11.2.
\enddemo

Now, given 
$h \in R(\text {Bun}_{n,X})$,
we can try to find an element
$h^{\prime} = \sum_{\bold f \in \bold{Cusp}_n} E(\bold f, \varphi_{\bold f})$
such that
$(E_{\bold f} (t), h-h^{\prime}) =0$
for any
$f \in \bold{Cusp}_n$.
For a given
$\bold f$
we denote by
$S(\bold f) \subset S_r$
the subgroup of permutations
$\sigma$
such that
$f_{\sigma (i)} = f_i, \forall i$,
and will look for 
$\varphi_{\bold f}$
invariant under
$S(\bold f)$.
By (3.7.9) the condition for
$(E_{\bold f}(t), h) = (E_{\bold f}(t), h^{\prime})$
is that

$$\biggl(\sum_{\sigma\in S(\bold f)} M^{\bold f}_{\sigma}
 (t)\biggr) \bar \varphi_{\bold f} (\bar t^{-1})\quad  =
 \quad (E_{\bold f}(t), h) \tag 3.9.4$$
Thus
$\varphi_{\bold f}(t) =\bar\Phi_{\bold f}(\bar t^{-1})$,
where

$$\Phi_{\bold f}(t) = \frac{(E_{\bold f}(t), h)} {\sum_{\sigma 
\in S(\bold f)} M^{\bold f}_{\sigma} (t)}. \tag 3.9.5$$

The function $\varphi_{\bold f}(t)$ is indeed invariant under
 $S(\bold f)$, but it is only
a rational function, not a Laurent polynomial, because of the 
zeros of
the denominator in $\Phi_{\bold f}(t)$. So in order to get a 
representation of $h'$ as a sum of
coefficients of Eisenstein series, we should first expand
 $\Phi_{\bold f}$.

Fix $\bold f = (f_1 \leq ... \leq f_r)$ and consider the 
coordinate vector space $\bold R ^r$
with the standard basis $e_1, ..., e_r$.
Let $\Gamma_{\bold f} \i \R ^r$ be the convex cone
with apex 0 generated by the vectors $e_j-e_i$ for $i< j$ 
such that $f_i = f_j$. 
For a Laurent series $\sum_{\omega\in \Z ^r} a_\omega t^\omega$ 
in $r$ variables
we call the set of $\omega$ such that $a_\omega \neq 0$ the support
 of the series.

\proclaim{(3.9.6) Lemma} The function $\Phi_{\bold f}(t)$ can
 be expanded into a Laurent series
$$\sum_{\omega\in \Z ^r \cap (\tau_{\bold f} + \Gamma_{\bold f})}
 a^{\bold f}_\omega t^\omega$$
whose support is contained in some translation of the cone 
$\Gamma_{\bold f}$. \endproclaim

\demo{Proof} First, the numerator $(E_{\bold f}(t), h)$ is a 
finite linear combination of
$(E_{\bold f}(t), [V]), V\in \text{Bun}_n(X)$, i.e., of Eisenstein 
series
$\Cal E_V(f_1, ..., f_r, t_1, ..., t_r)$. The poles of any such 
series come from 
the poles of 
$$M_\sigma^{\bold f}(t) = \prod_{i<j: \, \sigma(i) > \sigma(j)}
 q^{n_in_j(1-g_X)}\frac{ \text{LHom}(f_i,f_j, t_j/t_i)}
{\text{LHom}(f_i, f_j, t_j/qt_i)}, \quad \sigma\in S(\bold f)\leqno
 (3.9.7)$$
Each factor here is a  series in  non-negative powers of
 $t_j/t_i$, starting with 1.
So each $M_\sigma^{\bold f}(t)$ has support in a translation
 of $\Gamma_{\bold f}$. Now, for the
denominator $\sum M_\sigma^{\bold f}(t)$ we have, by the
 same reason, an expansion
$\sum_{\omega\in \Gamma_{\bold f}\cap \Z ^r} c_\omega 
t^\omega$ with $c_0\neq 0$.
Thus its inverse admits the geometric
 series expansion
$$\frac{1}{\sum M_\sigma^{\bold f}(t)} =
 \frac{1}{c_0} \sum_{m=0}^\infty \biggl(
-\sum_{\omega\in (\Gamma_{\bold f}\cap \Z ^r) - 
\{0\}} (c_\omega/c_0) t^\omega\biggr)^m$$
which is a series supported in
 $\Gamma_{\bold f}$. Lemma 3.9.6 is proved.\enddemo

Thus we can write $h'$ as an infinite series
$$h' = \sum_{\bold f\in \bold{Cusp}_n} E(\bold f, 
\varphi_{\bold f}) = 
\sum_{\bold f\in \bold{Cusp}_n}\sum_{\omega\in \Z ^r
 \cap (\tau_{\bold f} + \Gamma_{\bold f})} \bar
 a^{\bold f}_\omega E_{f_1, \omega_1} ... E_{f_r, \omega_r},
 \leqno (3.9.8)$$
where $E_{f,d}$ is the $d$-th coefficient of $E_f(t)$.

\proclaim{(3.9.9) Lemma} For any $\lambda > 0$ all but
 finitely many terms of the series 
(3.9.8), regarded as functions on $\text{Bun}_n(X)$,
 vanish outside $\text{Bun}_{n,d}^{\leq\lambda}$.
\endproclaim

This implies, in particular, that the series converges, 
as a series of functions on $\text{Bun}_n(X)$.

\demo{Proof} It is enough to conisder one $\bold f$ at a 
time. In order that 
$E_{f_1, \omega_1} ... E_{f_r, \omega_r}$ be nonzero
 on $V$, there should be a flag 
$V_1\i ... \i V_r = V$ of subbundles with
 $\deg(V_i/V_{i-1}) = \omega_i$. Thus, if
for at least one $i$ we have
$$\frac{\omega_i}{n_i} + \lambda < \frac
 {\omega_{i+1}}{n_{i+1}} \quad (\text{where}\quad f_i
\in \text{Cusp}_{n_i}), \leqno (3.9.10)$$
then $V$ must have a destabilizing flag with at least
 one gap in the slopes greater than $\lambda$,
which implies that in the Harder-Narasimhan filtration
 at least 
one gap in the slopes  will be greater than $\lambda$,
 i.e., that $V\in \text{Bun}_{n,d}^{>\lambda}$. 
It remains to notice that if $f_i=f_j$ then, of course 
$n_i=n_j$ so for all but finitely many
integer points lying in a translation of $\Gamma_{\bold f}$, 
the condition (3.9.10) will
be satisfied for some $i$. Lemma 3.9.10 is proved. \enddemo

So $h'$ is a well defined function on $\text{Bun}_n(X)$ and 
the scalar product $(E_{\bold f}(t), h')$
is indeed defined for any $\bold f\in \bold{Cusp}_n$ by (3.7.9)
 and is equal to 
$(E_{\bold f}(t), h)$ by construction. So we get the following 
fact.

\proclaim {(3.9.11) Lemma} The series in the right hand
 side of (3.9.8) converges to $h$.\endproclaim

\enddemo

\demo{\underbar{Step 3: Completion of the proof}} Fix $n$
 and assume, by induction, that any
$h\in R(\text{Bun}_{m,X})$ with $m<n$ lies in the image
 of $\pi_{\text{Bun}}$. 
Fix $d\in {\bold Z}$ and a big enough $\lambda\in {\bold Z}$.
 By Step 1, all $[V]$ with $V \in \text{Bun}_{n,d}^{>\lambda}$ 
lie in $\text{Im}(\pi_{\text{Bun}})$. So we $h$ to be equal to
$[V]$, $V\in \text{Bun}_{n,d}^{\leq\lambda}$. By Lemmas 3.9.9 
and 3.9.11, we can, by truncating the
series (3.9.8), find a finite linear combination $h''$ of
 elements from $\text{Im}(\pi_{\text{Bun}})$
which coincides with $h$ on $\text{Bun}_{n,d}^{\leq\lambda}$.
 But then $h-h''$ lies in
$\text{Im}(\pi_{\text{Bun}})$ by Step 1. This proves the theorem. 
\enddemo

\newpage

\centerline{\bf  \S 4.  Proof of Theorem 3.3 and Proposition 3.6.3.}

\vskip 1cm

Recall that the equality (3.3.1) has already been proved.
 So we first prove the rest 
of the assertions of Theorem 3.3.

\vskip .3cm

\noindent{\bf  (4.1) Proof of (3.3.2).}
Let
$\hat H(\text{Coh}_X)$
be the completion of the Hall algebra
$H(\text{Coh}_X)$
consisting by all, possibly infinite, sums
$\sum_{A \in \text{Coh}_X} f(A)[A]$.
Clearly
$\hat H(\text{Coh}_X)$
is a bimodule over the algebra
$H(\text{Coh}_X)$.
We denote the bimodule structure by $\circ$.

\proclaim{(4.1.1) Proposition}
For any
$\Cal F \in \text{Coh}_{0,X}$
and any
$f : \text{Bun}_n(X) \rightarrow \C$
we have equalities in
$\hat H(\text{Coh}_X)$:

$$[\Cal F] \circ \sum_{V \in \text{Bun}_n(X)} f(V)[V] =
 \sum_{V \in \text{Bun}_n(X)} f(V) [V \oplus \Cal F], \tag 4.1.1a$$

$$ \left(\sum_{V \in \text{Bun}_n(X)} f(V)[V]) \right)
 \circ [\Cal F]
\quad  = \quad \sum_{\Cal F^{\prime} \subset \Cal F}
 \sum_{W \in \text{Bun}_n(X)} (T_{\Cal F/\Cal
 F^{\prime}} f)(W) \cdot \tag 4.1.1b$$
$$\cdot \frac{|\text{Aut} (\Cal F^{\prime})|\cdot 
 |\text{Aut} (\Cal F/\Cal F^{\prime})|} {|\text{Aut}
 (\Cal F)|} \cdot q^{nh^0(\Cal F^{\prime})} [W \oplus \Cal F^{\prime}]. $$
\endproclaim

\demo{Proof}
 The statement (a) follows from the equality
$[\Cal F] \circ [V] = [V \oplus \Cal F]$
holding for any
$\Cal F \in \text{Coh}_{0,X}, V \in \text{Bun}_n(X)$,
because
$\text{Ext}^1(V,\Cal F) = \text{Hom}(\Cal F,V) =0$.
To see part (b), recall (1.2.1) the notation
$g^C_{AB}$
for the structure constants in the Hall algebra.
 Our statement follows from the next lemma about these constants.
\enddemo

\proclaim{(4.1.2) Lemma}
Let
$V,W \in \text{Bun}_nX$
and
$\Cal F, \Cal F^{\prime} \in \text{Coh}_{0,X}$.
Then

$$g^{W \oplus \Cal F^{\prime}}_{V,\Cal F} = \sum_{\Cal
 F^{\prime\prime}\in \text{Coh}_{0,X}} g^W_{V\Cal 
F^{\prime\prime}} \frac{|\text{Aut}(\Cal F^{\prime})
|\cdot |\text{Aut}(\Cal F^{\prime\prime})|} {|\text{Aut}
(\Cal F)|} \cdot q^{nh^0(\Cal F^{\prime})}$$\endproclaim 

Indeed, supposing the lemma true, we have

$$\gather \left(\sum_v f(V) [V])[ \Cal F]\right) = 
\sum \Sb V,W \in \text{Bun}_nX \\ \Cal F^{\prime} 
\in \text{Coh}_{0,X} \endSb g^{W \oplus\Cal F^{\prime}}_{V\Cal F} f(V) 
[W \oplus \Cal F^{\prime}] = \\
= \sum_{V,W} \sum_{\Cal F^{\prime}, \Cal F^{\prime\prime}}
 \frac {|\text{Aut} (\Cal F^{\prime} )|
\cdot |\text{Aut} (\Cal F^{\prime\prime})|} {|\text{Aut} 
(\Cal F)|} f(V) 
[W \oplus \Cal F^{\prime}] \cdot q^{nh^0(\Cal F^{\prime})} \\
= \sum_{\Cal F^{\prime} \subset \Cal F} \sum_W \sum \Sb V 
\subset W \\ W/V \simeq \Cal F^{\prime}/\Cal F^{\prime} 
\endSb f(V) \frac {|\text{Aut} (\Cal F^{\prime})|\cdot
 |\text{Aut} (\Cal F/\Cal F^{\prime})|} 
{|\text{Aut} (\Cal F)|} \cdot [W \oplus \Cal F^{\prime}] 
\cdot q^{nh^0(\Cal F)} \\
= \sum_{\Cal F^{\prime} \subset \Cal F} \sum_W 
(T_{\Cal F/\Cal F^{\prime}} f)(W) \frac {|\text{Aut} (\Cal F^{\prime}) 
|\cdot |\text{Aut} (\Cal F/\Cal F^{\prime})|} {|\text{Aut} (\Cal F)|} 
\cdot q^{nh^0(\Cal F^{\prime})} [W \oplus \Cal F^{\prime}], \endgather$$
as claimed. 

\demo{Proof of Lemma 4.1.2}
If $G$ is a finite group acting on a finite set $S$, we call the 
orbifold number of elements in $S$ modulo $G$ the quantity

$$\frac {|S|} {|G|} = \sum_{\{s\}} \frac 1 {|\text{Stab} (s)|} $$
where
$\{s\}$
are all the $G$-orbits in $S$, with
$s \in \{s\}$
being a chosen representative and $\text{Stab}(s)\subset G$ being
 its stabilizer.

Now, let
$V,W,\Cal F^{\prime},\Cal F$
be given. The subsheaf
$\Cal F^{\prime} \subset W \oplus \Cal F^{\prime}$
is defined intrinsically, as the maximal torsion subsheaf, so
$g^{W\oplus \Cal F^{\prime}}_{\Cal F^{\prime},W} =1$.
It follows that
$$g^{W \oplus \Cal F^{\prime}}_{V,\Cal F} = |\text{Aut}(W \oplus 
\Cal F^{\prime})| \cdot C,$$
where $C$ is the orbifold number of diagrams (``crosses'')

$$\CD @. @. 0 \\
@. @. @AAA \\
@. @. W \\
@. @. @AAA \\
0 @>>> V @>\alpha >> W \oplus \Cal F^{\prime} @>\beta >> \Cal F @>>> 0 \\
@. @. @AA\gamma A \\
@. @. \Cal F^{\prime} \\
@. @. @AAA\\
@. @. 0
\endCD \leqno (4.1.3)$$
modulo the product of the groups of automorphisms of all five objects.
 By elementary homological algebra, every cross can be completed to a 
$3 \times 3$
diagram with exact rows and columns, from which we want to retain the outer frame

$$\CD @. 0 @. @. 0 \\
@. @AAA @. @AAA \\
0 @>>> V^{\prime\prime} @>>> W @>>> \tilde\Cal F^{\prime\prime} @>>> 0 \\
@. @AAA @. @AAA \\
@. V @. @. \Cal F \\
@. @AAA @. @AAA \\
0 @>>> V^{\prime} @>>> \Cal F^{\prime} @>>> \tilde \Cal F^{\prime} @>>> 0 \\
@. @AAA @. @AAA \\
@. 0 @. @. 0 \endCD \tag4.1.4$$
Here, for instance,
$\tilde \Cal F^\prime = \text{Im} (\beta \gamma)$
etc. Notice that we have
$V^{\prime} =0$
since it should be a torsion sheaf embedded in $V$. Thus
$\tilde \Cal F^{\prime} \cong \Cal F$
and
$V^{\prime\prime} \simeq V$.
For each
$\Cal F^{\prime\prime}$
let
$C(\Cal F^{\prime\prime})$
be the orbifold number of crosses giving frames with the upper right corner
$\tilde \Cal F^{\prime\prime}$
isomorphic to
$\Cal F^{\prime\prime}$,
so that
$C = \sum_{\Cal F^{\prime\prime}} C(\Cal F^{\prime\prime})$.
Let also
$F(\Cal F^{\prime\prime})$
be the orbifold number of frames, i.e., arbitrary diagrams
 of the form (4.1.4) in which
$\tilde \Cal F^{\prime\prime} = \Cal F^{\prime\prime}, \tilde 
\Cal F^{\prime} = \Cal F^{\prime}, V^{\prime\prime} = V$,
(modulo the product of the groups of automorphisms of all 8
 objects constituting the frame).

The main result of Green ([Gr], Theorem 2) says that
$C(\Cal F^{\prime\prime}) = \langle \Cal F^{\prime\prime},
 V^{\prime}\rangle^2 \cdot F(\Cal F^{\prime\prime})$
which is equal to just
$F(\Cal F^{\prime\prime})$
since
$V^{\prime} =0$
and thus
$\langle \Cal F^{\prime\prime}, V^{\prime}\rangle =1$.
On the other hand,

$$F(\Cal F^{\prime\prime}) = g^{\Cal F}_{\tilde \Cal F^{\prime}
 \tilde \Cal F^{\prime\prime}} \cdot g^W_{V^{\prime\prime}\tilde
 \Cal F^{\prime\prime}} \cdot \frac {|\text{Aut}(V^{\prime\prime})
|\cdot|\text{Aut}(\tilde \Cal F^{\prime\prime})|\cdot
 |\text{Aut}(\tilde \Cal F^{\prime})|\cdot |\text{Aut}(V^{\prime})|}
 {|\text{Aut} (V)| \cdot |\text{Aut} (W)| \cdot |\text{Aut} (\Cal F)|
 \cdot |\text{Aut} \Cal F^{\prime}|} $$
as it follows from the general fact that
$g^C_{AB}$
is the orbifold number of exact sequences

$$0 \rightarrow A \rightarrow C \rightarrow B \rightarrow 0$$
modulo
$\text{Aut} (A) \times \text{Aut}(B)$.

Recalling our identifications, we find that

$$F(\Cal F^{\prime\prime}) = g^{\Cal F}_{\Cal F^{\prime} 
\Cal F ^{\prime\prime}} g^W_{V\Cal F^{\prime\prime}} \cdot \frac
 {|\text{Aut} (\Cal F^{\prime\prime})|} {|\text{Aut} (W)|
 \cdot |\text{Aut} (\Cal F)|}.$$
Notice also that

$$\gather |\text{Aut} (W \oplus \Cal F^{\prime})| = 
|\text{Aut} (W)| \cdot |\text{Aut} (\Cal F^{\prime})
| \cdot |\text{Hom}(W,\Cal F^{\prime})| = \\
= q^{nh^0(\Cal F^{\prime})} |\text{Aut} (W)| \cdot
 |\text{Aut}(\Cal F^{\prime})|. \endgather$$
Thus

$$\gather g^{W\oplus \Cal F^{\prime}}_{V,\Cal F} =
 |\text{Aut}(W \oplus \Cal F^{\prime})| \cdot C = 
|\text{Aut} (W \oplus \Cal F^{\prime})| \cdot 
\sum_{\Cal F^{\prime\prime}} C(F^{\prime\prime}) = \\
=|\text{Aut}(W \oplus \Cal F)| \cdot
 \sum_{\Cal F^{\prime\prime}} F(\Cal F^{\prime\prime}) = \\
= q^{nh^0(\Cal F^{\prime})} |\text{Aut}(W)| \cdot
 |\text{Aut}(\Cal F^{\prime})| \cdot \sum_{\Cal F^{\prime\prime}} 
g^{\Cal F} _{\Cal F^{\prime} \Cal F^{\prime\prime}} g^W_{V\Cal F^{\prime \prime}}
 \frac{|\text{Aut}(\Cal F^{\prime\prime})|} {|\text{Aut}(W)| 
\cdot |\text{Aut}(\Cal F)|} =\\
= q^{nh^0(\Cal F^{\prime})} \sum_{\Cal F^{\prime\prime}} 
g^{\Cal F}_{\Cal F^{\prime}\Cal F^{\prime\prime}} 
g^W_{V\Cal F^{\prime\prime}} \frac{|\text{Aut}(\Cal F^{\prime})
|\cdot |\text{Aut}(\Cal F^{\prime\prime})|} {|\text{Aut} (\Cal F)|} 
\endgather$$
as claimed.  Lemma 4.1.2 and thus Proposition 4.1.1
 are proved.\enddemo

For the product in the Ringel algebra $R$ we find from (4.1.1):

$$[\Cal F]  *  \sum_{V\in \text{Bun}_nX} f(V)[V] =
 q^{nh^0(\Cal F)/2} \sum_V f(V)[V \oplus \Cal F] \tag 4.1.5a$$

$$\left( \sum_{V \in \text{Bun}_n X} f(V)[V] \right)  * 
 [\Cal F] = q^{-nh^0(\Cal F)/2} \cdot \left( \sum_V f(V)[V]
 \right) \circ [\Cal F] =  \tag 4.1.5b $$
$$= \sum_{\Cal F^{\prime}\subset \Cal F}
 q^{-nh^0(\Cal F\Cal F^{\prime})/2} (T_{\Cal F/\Cal F^{\prime}} f)(W) 
\frac {|\text{Aut} (\Cal F^{\prime})|\cdot |\text{Aut}(\Cal F/\Cal F^{\prime})|}
 {|\text{Aut} (\Cal F)|} \cdot q^{\frac{nh^0(\Cal F^{\prime})} 2}
 [W \oplus \Cal F^{\prime}]. $$
Passing now from  just one basis vector 
$[\Cal F]$
to the generating function

$$\psi_g(t) = \sum_{\Cal F} \bar \chi_g([\Cal F]) t^{h^0(\Cal F)}
 |\text{Aut}(\Cal F)| \cdot [\Cal F],
\quad g \in \text{Cusp}_m$$
we find, for
$f \in \text{Cusp}_n$:

$$ \psi_g(t_2) *  E_f(t_1) = \leqno (4.1.6)$$
$$\sum \Sb \Cal F\in \text{Coh}_{0,X} \\ V \in \text{Bun}_n(X)
 \endSb \bar \chi_g([\Cal F]) f(V) q^{+nh^0(\Cal F)/2} |\text{Aut}(\Cal F)| \cdot 
\cdot t^{h^0(\Cal F)}_2 t^{\deg (V)}_1 [V \oplus \Cal F]. $$
To find the product in the opposite order, we use (4.1.5b) 
together with the fact that $f$ is a Hecke eigenform, and find

$$ E_f(t_1)  *  \psi_g(t_2) = \sum_{\Cal F \in \text{Coh}_{0,X}}
 \sum_{W\in \text{Bun}_n X} \sum_{\Cal F^{\prime},
 \Cal F^{\prime\prime} \in \text{Coh}_{0,X}} g^{\Cal F}_{\Cal F^{\prime}
 \Cal F^{\prime\prime}} \chi_f([\Cal F^{\prime\prime}])f (W) \cdot \leqno (4.1.7)$$
$$\cdot \bar \chi_g([\Cal F]) q^{-n(h^0(\Cal F^{\prime\prime})-
h^0(\Cal F^{\prime}))/2} |\text{Aut}(\Cal F^{\prime})| \cdot
 |\text{Aut}(\Cal F^{\prime\prime})| \cdot $$
$$\cdot t^{h^0(\Cal F)}_2 t^{\deg W-h^0(\Cal F^{\prime\prime})} _1 
[W \oplus \Cal F^{\prime}],  $$
where we have also replaced the summation over subsheaves
$\Cal F^{\prime} \subset \Cal F$
by the summation over arbitrary pairs of isomorphism classes
$\Cal F^{\prime}, \Cal F^{\prime\prime}$,
with the factor
$g^{\Cal F}_{\Cal F^{\prime}\Cal F^{\prime\prime}}$
counting the number of subsheaves of type
$\Cal F^{\prime}$
and cotype
$\Cal F^{\prime\prime}$.
Since
$\bar \chi_g$
is a character of the Hall algebra, and
$\sum_{\Cal F} g^{\Cal F}_{\Cal F^{\prime}\Cal F^{\prime\prime}}
 [\Cal F] = [\Cal F^{\prime}] \circ [\Cal F^{\prime\prime}]$,
we can replace the RHS of (4.1.7) by
$$\sum_{W\in \text{Bun}_nX} \,\, \sum_{\Cal F^{\prime}, 
\Cal F^{\prime\prime} \in \text{Coh}_{0,X}}
\chi_f([\Cal F^{\prime\prime}]) f(W) \bar
 \chi_g([\Cal F^{\prime}]) \bar \chi_g([\Cal F^{\prime\prime}]) 
\cdot \leqno (4.1.8)$$
$$\cdot q^{nh^0(\Cal F^{\prime})/2-nh^0(\Cal F^{\prime\prime})/2}
 |\text{Aut}(\Cal F^{\prime})| \cdot |\text{Aut}(\Cal F^{\prime\prime})| 
\cdot t_2^{h^0(\Cal F^{\prime})+h^0(\Cal F^{\prime\prime})}
 \cdot [W \oplus \Cal F^{\prime}]. $$
This can be written as the product

$$ \left\{ \sum_{W \in \text{Bun}_nX} \sum_{\Cal F^{\prime}
 \in \text{Coh}_{0,X}} \bar \chi_g([\Cal F^{\prime}]) f(W)
 q^{nh^0(\Cal F^{\prime})/2} \cdot |\text{Aut}(\Cal F^{\prime})| 
\cdot t^{h^0(\Cal F^{\prime})}_2 t^{\deg (W)}_1 \cdot
 [W \oplus \Cal F^{\prime}] \right\} \times \leqno (4.1.9)$$
$$\times \left\{ \sum_{\Cal F^{\prime\prime}\in \text{Coh}_{0,X}}
 \chi_f([\Cal F^{\prime\prime}]) \bar \chi_g([\Cal F^{\prime \prime}])
 q^{nh^0(\Cal F^{\prime\prime})/2} \cdot
 |\text{Aut}(\Cal F^{\prime\prime})|(t_2/t_1)^{h^0(\Cal
 F^{\prime\prime})} \right\},$$
of which the first factor is identical, up to renaming
 the summation arguments, with the RHS of (4.1.6), i.e.,
 it is equal to
$\psi_g(t_2)  *  E_f(t_1)$.
The second factor has the form
$\Lambda (t_2/q^{\frac n 2} t_1)$,
where

$$\Lambda (t) = \sum_{\Cal F \in \text{Coh}_{0,X}} \chi_f([\Cal F])
 \bar \chi_g ([\Cal F]) \cdot |\text{Aut} (\Cal F)|t^{h^0(\Cal F)} \tag 4.1.10$$
For
$x \in X$
and an integer sequence
$\mu = (\mu_1 \ge \ldots \ge \mu_r \ge 0)$
denote
$\Cal F_{x,\mu} = \bigoplus \Cal O/I^{\mu_i}_x$.
Since every
$\Cal F \in \text{Coh}_{0,X}$
has a unique decomposition into the direct sum of sheaves of the form
$\Cal F_{x,\mu}$ (one for each $x\in X$),
we can decompose
$\Lambda (t)$
into the Euler product:

$$\Lambda (t) = \prod_{x \in X} \Lambda_x(t^{\deg (x)}), \tag 4.1.11$$
$$\Lambda_x(t) = \sum_{\mu} \chi_f([\Cal F_{x,\mu}]) \bar
 \chi_g([\Cal F_{x,\mu}])
 \cdot |\text{Aut}(\Cal F_{x,\mu} )| t^{|\mu|},$$
where
$|\mu| = \sum \mu_i$.
By (2.5.2) and (2.6.13),

$$\chi_f([\Cal F_{x,\mu}]) = q^{-\sum(i-1)\mu_i}_x P_{\mu} 
\left(q^{\frac{n-1}2}_x \lambda_{1,x}(f)^{-1}, \ldots,
 q^{\frac{n-1}2}_x \lambda_{n,x} (f)^{-1}; q^{-1}_x \right)
 \tag 4.1.12$$
$$\bar \chi_g([\Cal F_{x,\mu}]) = q_x^{-\sum (i-1)\mu_i} 
P_{\mu} \left( q^{\frac{m-1}2}_x \lambda_{1,x}(g), \ldots,
 q^{\frac{m-1}2}_x \lambda_{m,x} (g); q^{-1}_x \right) $$
where
$P_{\mu}$
is the Hall-Littlewood polynomial. It is known ([Mac], Ch. 
III, formula (4.4)) that

$$\sum_{\mu} b_{\mu} P_{\mu} (z_1 \ldots z_n; q^{-1}_x) 
P_{\mu} (w_1,\ldots, w_m; q^{-1}_x) = \prod_{i,j} \frac
 {1-q^{-1}_xz_iw_j} {1-z_iw_j}, \tag 4.1.13$$
where

$$b_{\mu} = q_x^{-|\mu|-2\Sigma(i-1)\mu_i} \cdot |\text{Aut}
 (\Cal F_{x,\mu})|.\tag 4.1.14$$
Therefore

$$\Lambda_x(t) = \prod_{i,j} \frac{1-q^{\frac{n+m}2 -1}_x 
\frac{\lambda_{j,x}(g)} {\lambda_{i,x}(f)} t} {1-q_x^{\frac{n+m} 2} 
\frac{\lambda_{j,x}(g)} {\lambda_{i,x} (f)} t},$$
and so

$$\Lambda (t) = \frac{\text{LHom}(f,g,q^{\frac{n+m} 2} t)} 
{\text{LHom} (f,g,q^{\frac{n+m} 2 -1} t)}
 \tag 4.1.15$$
Since, as we already saw,

$$E_f(t_1)  *  \psi_g(t_2) = \Lambda(t_2/q^{\frac n 2} t_1)
 \Psi_g(t_2)  *  E_f(t_1),$$
we have proved the formula (3.3.2).

\vskip .3cm

\noindent{\bf  (4.2) Proof of (3.3.3).}
Recalling that the bilinear form
$\langle \alpha, \beta \rangle$
vanishes on
$\Cal K_0(\text{Coh}_{0,X})$
and that
$\bar \chi_f$
is a homomorphism
$H(\text{Coh}_{0,X}) \rightarrow \C$,
we find, from (1.6.3):

$$\gather \Delta \psi_f(t) = \sum_{\Cal F \in \text{Coh}_{0,X}} 
\bar \chi_f([\Cal F])t^{h^0(\Cal F)} |\text{Aut}(\Cal F)| \cdot
 \Delta ([\Cal F]) = \\
= \sum_{\Cal F} \sum_{\Cal F^{\prime}, \Cal F^{\prime\prime}}
 \bar \chi_f ([\Cal F])t^{h^0(\Cal F)} |\text{Aut}(\Cal F)| 
g^{\Cal F}_{\Cal F^{\prime} \Cal F^{\prime\prime}} \frac 
{|\text{Aut}(\Cal F^{\prime})|\cdot |\text{Aut}(\Cal F^{\prime\prime})|}
 {|\text{Aut} (\Cal F)|} \cdot [\Cal F^{\prime}] \otimes c_{\Cal F^{\prime}} 
[\Cal F^{\prime\prime}] = \\
= \sum_{\Cal F^{\prime}, \Cal F^{\prime\prime}} \bar 
\chi_f([\Cal F^{\prime}]) \bar \chi_f([\Cal F ^{\prime\prime}])
t^{h^0(\Cal F^{\prime}) + h^0(\Cal F^{\prime\prime})} \cdot 
|\text{Aut} (\Cal F^{\prime}|\cdot |\text{Aut}(\Cal F^{\prime\prime}) \cdot \\
\cdot [\Cal F^{\prime}] \otimes c_{\Cal F^{\prime}}
 [\Cal F^{\prime\prime}] = \psi_f(t \otimes c)(1 \otimes \Psi_f(t))
\endgather$$
as claimed.
\vskip .3cm

\noindent{\bf  (4.3) Proof of (3.3.4).}
Let
$r_0 : AF_n \rightarrow \bigoplus_{i+j=n} AF_i \otimes AF_j$
be the linear map dual to the Hall multiplication
$AF^0_i \otimes AF^0_j \rightarrow AF^0_{i+j}$
(with respect to the orbifold scalar product (1.3.2) on each
$AF^0_i$).
Thus for
$f \in AF_n$
and vector bundles
$W^{\prime}, W^{\prime\prime}$
of ranks
$i,j$
we have

$$(r_0f)(W^{\prime}, W^{\prime\prime}) = \sum_{V \in 
\text{Bun}_nX} g^V_{W^{\prime}W^{\prime\prime}} 
\frac{|\text{Aut}(W^{\prime})| \cdot |\text{Aut}(W^{\prime\prime})|} 
{|\text{Aut}(V)|} f(V), \tag 4.3.1$$
and $f$ is a cusp form if
$r_0(f) = 1 \otimes f + f \otimes 1$.
Let now
$f \in \text{Cusp}_n$
be given. Then, by definition

$$ \Delta E_f(t) =\sum_{V\in \text{Bun}_n(X)} f(V)t^{\deg(V)}
 \Delta ([V]) = \tag 4.3.2 $$
$$= \sum_V \sum_{U,\Cal G} \langle \Cal G, U\rangle \,
 g^V_{U\Cal G} f(V)t^{\deg(V)} \frac {|\text{Aut} (U)|\cdot 
|\text{Aut}(\Cal G)|} {|\text{Aut} (V)|} [U] \otimes K_U[\Cal G], $$
where $U$ and
$\Cal G$
can be, a priori, any coherent sheaves on $X$. However, if
$g^V_{U\Cal G} \ne 0$,
$U$ can be embedded as a subsheaf into $V$ and so is locally free. The sheaf
$\Cal G$
may not be locally free, and we write it as
$\Cal G = W \oplus \Cal F$
where $W$ is a vector bundle,
$\Cal F$
is a torsion sheaf (their isomorphism classes are determined by
$\Cal G$).
Further, if we have an exact sequence

$$0 \rightarrow U \rightarrow V \rightarrow W \oplus \Cal F
 \rightarrow 0 \tag 4.3.3$$
we can consider the unqiue subbundle (i.e., subsheaf which is
 locally a direct summand)
$\bar U \supset U,\,  \text{rk} (\bar U) = \text{rk} \, U$.
The sequence (4.3.3) gives thus two sequences

$$\align &0 \rightarrow \bar U \rightarrow V \rightarrow W 
\rightarrow 0 \\
&0 \rightarrow U \rightarrow \bar U \rightarrow \Cal F \rightarrow 0
 \endalign$$
and we conclude that

$$g^V_{U,W \oplus \Cal F} = \sum_{\bar U \in \text{Bun} (X)}
 g^{\bar U} _{U\Cal F} g^{\Cal F}_{\bar U W} \tag 4.3.4$$
Note also that

$$|\text{Aut}(W \oplus \Cal F)| = |\text{Aut}(W)| \cdot
 |\text{Aut}(\Cal F)| \cdot q^{\text{rk}(W) \cdot h^0(\Cal F)} \tag 4.3.5$$
Keeping these two formulas in mind and  substituting
$\Cal G = W \oplus \Cal F$
into (4.3.2), we find

$$\gather \Delta E_f(t) = \sum_{V,U,W,\bar U, \Cal F}
 \langle W \oplus \Cal F, U\rangle g^{\bar U} _{U\Cal F}
 g^V_{\bar U W} \frac {|\text{Aut}(\bar U)| \cdot |\text{Aut}(W)|
 \cdot |\text{Aut}(U)| \cdot |\text{Aut} (\Cal F)|} {|\text{Aut}(V)|
 \cdot |\text{Aut}(\bar U)|} \cdot \\
\cdot q^{rk(W) \cdot h^0(\Cal F)} \cdot f(V)t^{\deg (V)} [U] \otimes
 K_U[W \oplus \Cal F] \endgather$$
where, in addition, we have multiplied and divided by
$|\text{Aut} (\bar U)|$.
In this sum, 
$V \in \text{Bun}_n(X)$,
while
$U,W,\bar U$
are isomorphism classes of vector bundles of arbitrary rank, and
$\Cal F \in \text{Coh}_{0,X}$.
By (4.3.1) we can write the result of summation over $V$ in terms of
$r_0$,
getting

$$ \Delta E_f(t) = \sum_{U,W,\bar U,\Cal F} g^{\bar U}_{U\Cal F}
 (r_0f)(\bar U,W) \frac {|\text{Aut}(U)| \cdot |\text{Aut}(\Cal F)|}
 {\text{Aut}(\bar U)|} \cdot \tag 4.3.7$$
$$\cdot \langle W \oplus \Cal F,U \rangle q^{rk(W) \cdot h^0(\Cal F)} 
t^{\deg \bar U +\deg W} [U] \otimes K_U [W \oplus \Cal F] $$
Since $f$ is a cusp form,
$(r_0f)(\bar U,W) =0$
unless
$\bar U =0$
or
$W =0$.
Summation with
$\bar U=0$
gives
$1 \otimes E_f(t)$.
If
$W =0$,
the summation gives

$$\sum_{U,\bar U,\Cal F} g^{\bar U} _{U\Cal F} f(\bar U) \frac 
{|\text{Aut} (U)| \cdot |\text{Aut}(\Cal F)|} {|\text{Aut} \bar U|}
\langle \Cal F, U\rangle t^{\deg \bar U} [U] \otimes K_U[\Cal F] =
 \tag 4.3.8 $$
$$= \sum_{U,\Cal F} (T^V_{\Cal F} f)(U) \cdot |\text{Aut} (\Cal F)
|q^{\frac n 2 h^0(\Cal F)} t^{h^0(\Cal F) +\deg U} [U] \otimes K_U 
[\Cal F] = $$
$$= \sum_{U, \Cal F} \bar \chi_f([\Cal F]) f(U) |\text{Aut}(\Cal F)
|q^{-\frac n 2 h^0(\Cal F)} t^{h^0(\Cal F) +\deg U} [U] \otimes K_U[\Cal F] $$
where we used Propositions (2.6.8) and 2.6.11 to identify the action of
$T^V_{\Cal F}$.
Now, the last expression in (4.3.8) factors into the product

$$\left(\sum_{U \in \text{Bun}_n(X)} f(U)t^{\deg(U)} [U] \otimes
 K_U \right) \left( 1 \otimes \sum_{\Cal F\in \text{Coh}_{0,X}} 
\bar\chi_f([\Cal F])|\text{Aut}(\Cal F)|(q^{-\frac n 2}t)^{h^0(\Cal F)}
 [\Cal F]\right)$$
By using the equality
$K_U = K^n c_{\det(U)}$
we can rewrite this as
$$E_f(t \otimes c)(1 \otimes K^n\psi_f(q^{-\frac n 2} t))$$
thus proving (3.3.4).
\vskip .3cm

\noindent{\bf  (4.4) End of the proof of Theorem 3.3.}
It remains to prove the equalities (3.3.5) and (3.3.6) describing
 the counit and the antipode. Now, (3.3.5) is obvious from the
 definition (1.6.5) of
$\epsilon$.
To see the first equality in (3.3.6), notice that

$$ S(\psi_f (t)) = \sum_{\Cal F \in \text{Coh}_{0,X}} \bar
 \chi_f([\Cal F]) \cdot |\text{Aut}(\Cal F)| \cdot t^{h^0(\Cal F)} S([\Cal F]) 
= \tag 4.4.1 $$
$$= \sum_{\Cal F} \sum^{\infty}_{m=1} (-1)^m \sum_{\Cal F_0 \subset
 \ldots \subset \Cal F_m = \Cal F} \bar \chi_f([\Cal F]) \prod_i 
|\text{Aut}(\Cal F_i/\Cal F_{i-1})|\cdot \prod_i t^{h^0(\Cal F_i/\Cal F_{i+1})} \cdot $$
$$\cdot \prod_i[\Cal F_i/\Cal F_{i-1}] \cdot
 \prod_i c^{-1} _{\overline{\Cal F_i/\Cal F_{i-1}}}. $$
By replacing the summation over
$\Cal F$
and then over flags of subsheaves
$\Cal F_0 \subset \ldots \subset \Cal F_m = \Cal F$
in a given 
$\Cal F$
by the summation over (independent) isomorphism classes of
$\Cal G_i = \Cal F_i/\Cal F_{i-1}$
with coefficient
$g^{\Cal F}_{\Cal G_0 \ldots \Cal G_m}$ (the number of filtrations 
on $\Cal F$ with quotients
$\Cal G_0, ..., \Cal G_m$)
and using the fact that
$\bar \chi_f$
is a character of the Hall algebra, we bring (4.6.1) to the form

$$1+\sum^{\infty}_{m=1} (-1)^m (\psi_f(c^{-1}t)-1)^m$$
i.e., to the geometric series for
$\psi_f(c^{-1}t)^{-1}$,
as claimed. To prove the second equality in (3.3.6), we write, for
$f \in \text{Cusp}_n$:

$$\gather S (E_f(t)) = \sum_{V \in \text{Bun}_n(X)} f(V)t^{\deg (V)}
 S([V]) = \\
= \sum_{V \in \text{Bun}_n(X) } \sum^{\infty}_{n=1} (-1)^m \sum_{\Cal
 G_0 \subset \ldots \subset \Cal G_m = V} f(V)t^{\deg(V)} \prod^m_{i=1}
\langle \Cal G_i/\Cal G_{i-1}, \Cal G_{i-1} \rangle \cdot \\
\cdot \frac{\prod^m_{j=0} |\text{Aut} (\Cal G_j/\Cal G_{j-1})|} {|\text{Aut} (V)|} 
[\Cal G_0] \ldots [\Cal G_m/\Cal G_{m-1}] \cdot K^{-n}c^{-1}_{\det(V)} \endgather$$
Since $f$ is a cusp form, the reasoning similar to that in (4.3), shows
 that only flags
$\Cal G_0 \subset \ldots \subset \Cal G_m = V$
with
$\text{rk} (\Cal G_i) = \text{rk} (V), \forall i$,
contribute to the total sum. But when we restrict the summation
 to such flags (in which we thus have
$\Cal G_i/\Cal G_{i-1} \in \text{Coh}_{0,X}$
for
$i > 0)$
we immediately get
$-E_f(c^{-1}t)\psi_f(q^{-\frac n 2} t)^{-1}K^{-n}$
by the same reason as above (summation of geometric series) plus 
the application  of the Riemann-Roch theorem to account for
$\prod \langle \Cal G_i/\Cal G_{i-1}, \Cal G_{i-1} \rangle$.

Theorem 3.3 is completely proved.
\vskip .3cm

\noindent{\bf  (4.5) Proof of Proposition 3.6.3.}
The equality (3.6.4) is obvious: non-isomorphic objects give orthogonal 
elements in the Hall algebra. To see (3.6.5) notice that cusp eigenforms
 with different eigenvalues of Hecke operators are orthogonal, by (2.6.11), 
so
$(E_f(t_1), E_g(t_2)) =0$
for
$f \ne g$.
If
$f=g$,
the equality (3.6.5) follows at once from the definition and the
 assumption (2.6.14) that
$\parallel f \parallel^2_d =1$.
So we concentrate on the proof of (3.6.6). Notice that
$\psi_f(t)$
can be written as the Euler product

$$\psi_f(t) =\prod_{x\in X} \psi_{f,x} (t^{\deg (x)}), 
\quad \psi_{f,x} (t) =\sum_{\mu} \bar \chi_f ([\Cal F_{x,\mu}])
 \cdot |\text{Aut}(\Cal F_{x,\mu})| \cdot t^{|\mu|} \cdot [\Cal F_{x,\mu}] \tag 4.5.1$$
where
$\mu$
runs over all partitions
$\mu_1 \ge \ldots \ge \mu_n \ge 0$,
with
$|\mu| = \sum \mu_i$
and
$\Cal F_{x,\mu} = \bigoplus \Cal O_x/I^{\mu_i}_x$
having the same meaning as in (4.1), more precisely, after (4.1.10). Thus we have:

$$a_f(t) = \sum_{x\in X} a_{f,x} (t^{\deg (x)}), \quad a_{f,x} (t)
 = \log \psi_{f,x} (t), \tag 4.5.2$$
and therefore

$$(a_f(t_1), a_g(t_2)) = \sum_{x\in X} \bigl(a_{f,x} (t^{\deg (x)}_1),
 a_{g,x} (t^{\deg (x)}_2)\bigl). \tag 4.5.3$$
We will evaluate each summand in this sum. Let

$$\text{Ch}: H(\text{Coh}_{x,X}) \rightarrow \Lambda,  \quad
[\Cal F_{x,\mu}] \mapsto q^{-\sum (i-1)\mu_i} P_{\mu} (z_1, \ldots, z_N; q_x^{-1}) \tag 4.5.4$$
be the isomorphism discussed in (2.3.5d).

In [Mac], Ch. III, formula (4.8), Macdonald defined
 a certain scalar product on
$\Lambda [t]$
with values in
$\Q[t,t^{-1}]$.
We denote by
$( \quad , \quad)_{Macd}$
the
$\Q$
-valued scalar product on
$\Lambda$
obtained from this
$\Q[t,t^{-1}]$
-valued scalar product by specializing to
$t = q_x^{-1}$.
Let us, for the time being, abbreviate
$H(\text{Coh}_{x,X})$
to simply $H$ and introduce on the algebra $H$ the grading
$H = \bigoplus_{d \ge 0} H_d$,
where
$H_d$
is linearly spanned by
$\Cal F_{x,\mu}, |\mu| =d$.

\proclaim{(4.5.5) Lemma}
Let
$u_i \in H_{d_i}, i=1,2$.
If
$d_1 \ne d_2$,
then
$(u_1,u_2) =0$,
and if
$d_1 = d_2 = d$,
then

$$(u_1, u_2) = q_x^{-d} (\text{Ch}(u_1), \text{Ch}(u_2))_{Macd}.$$
\endproclaim

\demo{Proof}
The first statement is obvious. The second follows from
 the equality (formula III (4.9) of [Mac])

$$\bigl(P_{\mu} (z;q_x^{-1}), Q_{\mu^{\prime}}(z; q_x^{-1})\bigr)_{Macd}
 = \delta_{\mu\mu^{\prime}},$$
where
$Q_{\mu^{\prime}} (z;q_x^{-1}) = b_{\mu^{\prime}} P_{\mu^{\prime}}
 (z,q_x^{-1})$
and
$b_{\mu}$
is defined in (4.1.14), while on the other hand,

$$(\Cal F_{x,\mu}, \Cal F_{x,\mu^{\prime}}) =
 \frac{ \delta_{\mu\mu^{\prime}}} {|\text{Aut} (\Cal F_{x,\mu})|}.$$
The lemma implies that

$$\bigl(a_{f,x}(t_1), a_{g,x}(t_2)\bigl) = 
\bigl(\text{Ch}(a_{f,x} (q_x^{-1}t_1)), \text{Ch}(a_{g,x}(t_2))\bigl)_{Macd}. 
\tag 4.5.6$$
Let us therefore find
$\text{Ch}(a_{f,x} (t))$. Let us omit the mention of the parameter 
in the Hall polynomials, which we will
always implicitly understand to be equal to $q_x^{-1}$.
We start by writing

$$\gather \text{Ch} \, (\psi_{f,x} (t)) = 
\sum_{\mu} q^{-2\sum (i-1)\mu_i}_x P_{\mu}
 \bigl(q_x^{\frac{n-1}2} \lambda_x(f)\bigr)P_{\mu} (z)|\text{Aut}
(\Cal F_{x,\mu})|t^{|\mu|} = \\
= \sum_{\mu} q^{|\mu|}_x b_{\mu} P_{\mu} \bigl(q_x^{\frac{n-1}2}
 \lambda_x(f)\bigr)
 P_{\mu} (z) |t|^{\mu} \endgather$$
where we used (4.1.12) and set
$\lambda_x(f) = (\lambda_{x,1}(f), \ldots, \lambda_{x,n} (f))$
as well as
$z =(z_1, \ldots, z_N)$.
By the formula quoted in (4.1.13), we find:

$$\text{Ch} \, (\psi_{f,x} (t)) \quad  = \quad 
\prod_{i,j} \frac{1-q^{\frac {n-1}2}_x \lambda_{i,x} (f)z_j
 \, t} {1-q^{\frac {n-1}2}_x \lambda_{i,x} (f) z_j \, qt}$$
and hence

$$\gather\text{ Ch}  (a_{f,x}(t)) \quad = \quad  \log\, \text{ Ch} 
 (\psi_{f,x} (t))
\quad  = \quad 
\sum^{\infty} _{d=1} \frac{1-q^d}d t^d \sum_{i,j} \left( q^{\frac
 {n-1}2}_x \lambda_{i,x} (f) z_j\right)^d = \\
= \sum_d \frac{1-q^d}d p_d(q_x^{\frac{n-1}2} \lambda_x(f)) p_d(z) t^d \endgather$$
where

$p_d$
is the
$d-th$
power sum symmetric function. Now, by formula III (4.11) of [Mac],

$$(p_d, p_{d^{\prime}})_{Macd}= \frac d {1-q_x^{-d}} \,
 \delta_{dd^{\prime}}.$$
By using (4.5.6), we now find

$$ (a_{f,x}(t_1), a_{g,x}(t_2)) =\sum^{\infty}_{d=1}
 \frac{q^d-1} d p_d(q_x^{\frac {n-1}2} \lambda_x(f))
p_d(q_x^{\frac {m-1}2}\overline{\lambda_x(g)} )t^d_1 \bar t^d_2 = $$
$$= \log \prod_{i,j} \frac {1-q_x^{\frac{n+m}2} \lambda_{x,i} (f)
 \overline{\lambda_{x,j}(g)} t_1\bar t_2} {1-q_x^{\frac {n+m}2 -1}
 \lambda_{x,i}(f) \overline{\lambda_{x,j}(g)} t_1\bar t_2},$$
and (3.6.8) is obtained by performing the summation over
$x \in X$,
and noticing (2.6.13) that the set of the
$\overline {\lambda_{x,j}(g)}$
is the same as the set of the
$\lambda_{x,j} (g)^{-1}$.
\enddemo

\newpage

\centerline{\bf  \S 5. Quantum affine algebras.} 

\vskip 1cm

In this section we compare Hopf algebras formed by
 automorphic forms with quantum affine algebras.
 For our purposes it is convenient to treat the
 quantization parameter $q$ as a constant rather
 than as an indeterminate variable.

\vskip .3cm

\noindent{\bf  (5.1) Drinfeld's realization of $U_q(\hat \Cal G)$.}
Fix a  nonzero complex number $q$. Let
$\bold g$
be a finite-dimensional semisimple Lie algebra whose Cartan matrix
$A = \parallel a_{ij} \parallel_{i,j=1,\ldots,r}$
is symmetric and let
$\widehat {\bold g}$
be the corresponding Kac-Moody algbera (central extension of
$\bold g[t,t^{-1}]$).
The quantization
$U_q(\widehat {\bold g})$
can be defined in two ways: the first (root realization)
 uses the system of simple roots for the affine root system of
$\widehat{\bold g}$
and proceeds directly from the affine Cartan matrix
$\hat A$
of
$\widehat {\bold g}$
(of size
$(r+1) \times (r+1))$.
The other, the so-called loop realization of Drinfeld [Dr1-2]
 (see also [CP]) is the 
$\bold C$-algebra generated by the symbols

$$x^+_i(n), x^-_i(n), i =1,\ldots, r,n \in \bold Z, \quad 
 k^{\pm}_{i,n}, i=1, \ldots, r, n \ge 0$$
and the central element $c$. They are subject to some relations
 which are best written in terms of formal generating functions

$$F^{\pm}_i(t) = \sum_{n\in \Z} x^{\pm}_i(n)t^{+n}, \quad
\varphi^{\pm}_i(t) = \sum^{\infty}_{n=0} k^{\pm}_{i,n} t^{\pm n}.$$
The relations have the form:

$$k^+_{i,0}k^-_{i,0} = k^-_{i,0}k^+_{i,0}=1, \quad [\varphi^{\pm}_i(t_1), 
\varphi^{\pm}_j(t_2)] =0, \tag 5.1.1$$

$$ (t_1-\sqrt q^{\pm a_{ij}} t_2) X^{\pm}_i (t_1)X^{\pm}_j(t_2) =
 (\sqrt q^{\pm a_{ij}} t_1-t_2) X^{\pm}_j(t_2)X^{\pm}_i(t_1), \tag 5.1.2 $$

$$F^{\pm}_i(t_1) \varphi^+_j(t_2) = \left( \frac{c^{\mp 1/2}
 t_2/t_1-\sqrt q^{a_{ij}}} {\sqrt q^{a_{ij}} c^{\mp 1/2}
 (t_2/t_1)-1} \right) ^{\pm 1} \varphi^+_j(t_2)F^{\pm}_i(t_1),\tag 5.1.3$$

$$F^{\pm}_i(t_1) \varphi^-_j(t_2) = \left( \frac{c^{\mp 1/2} 
t_1/t_2-\sqrt q^{a_{ij}}} {\sqrt q^{a_{ij}} c^{\mp 1/2}
 (t_1/t_2)-1} \right) ^{\mp 1} \varphi^-_j(t_2)F^+_i(t_1),\tag 5.1.4$$

$$[F^+_i(t_1), F^-_j(t_2)] = \delta_{ij} \left\{ \delta \left(
 \frac{t_1}{ct_2} \right) \varphi^-_i (c^{1/2}t_2)-\delta 
\left(\frac{ct_1}{t_2}\right) \varphi^+_i (c^{1/2}t_1) \right\}
 \cdot (q^{1/2} -q^{-1/2})^{-1}, \tag 5.1.5$$

$$\varphi^+_i(t_1) \varphi^-_j(t_2) = \frac {(\sqrt q^{a_{ij}} 
c^{-1} t_1/t_2-1)(ct_1/t_2-\sqrt q^{a_{ij}})} {(c^{-1}t_1/t_2
 -\sqrt q^{a_{ij}})(\sqrt q^{a_{ij}} ct_1/t_2 -1)} \varphi^-_j(t_2)
\varphi^+_i(t_1), \tag 5.1.6$$

$$\underset {t_1,\ldots,t_m} \to {\text{Sym}} \sum^m_{l=0} (-1)^l 
\bmatrix m \\ l \endbmatrix_q F^{\pm}_i(t_1)
 \ldots F^{\pm}_i(t_l)F^{\pm}_j(s) F^{\pm}_i(t_{l+1}) \ldots F^{\pm}_i 
(t_m)=0, \tag 5.1.7$$
where
$m=1-a_{ij}$
and
Sym
stands for the sum over all the permutations of
$t_1, \ldots, t_m$.

The relations (5.1.7) are the analogs of the Serre relations in 
finite-dimensional semisimple Lie algebras.

Let
$U_q(\widehat{\bold n^+})$,
(resp.
$U_q(\widehat{\bold n^-})$)
denote the subalgebra in
$U_q(\widehat {\bold g})$
generated by the elements
$x^+_i(n)$
(resp.
$x^-_i(n))$
only, which are subject to relations (5.1.2), (5.1.7). Let also
$U_q(\widehat{\bold b^+})$
denote the subalgebra generated by the
$x^+_i(n), n \in \bold Z$
and by the elements
$k^+_i(n), n \ge 0$
Similarly for
$U_q(\widehat{\bold b^-})$.
They are the quantizations of the enveloping algebras of the Lie algebras

$$\widehat{\bold n^{\pm}} = \bold n^{\pm}[t,t^{-1}], \quad
\widehat{\bold b^{\pm}} = \bold n^{\pm} [t,t^{-1}] \oplus \bold h[t^{\pm 1}],$$
where
$\bold g = \bold n^+ \oplus \bold h \oplus \bold n^-$
is the standard decomposition of
$\bold g$
into the nilpotent and Cartan subalgebras.

\vskip .3cm

\noindent{\bf  (5.2) The case
$\bold g = sl_2$
and sheaves on
$P^1$.}
Consider the simplest case
$\bold g = sl_2$.
In this case there is only one root, so we will denote the
 generators and generating functions as
$x^{\pm}(n), F^{\pm}(t)$
etc. The algebra
$U_q(\widehat {\bold b^+})$
is generated by $c$ and by the coefficients of the two power series
$F^+(t), \varphi^+(t)$.

On the other hand, let us suppose that $q$ is a prime power
 and specialize the theory of $\S$ \, 3 to the case of the
 simplest algebraic curve $X$, namely the projective line
$P^1_{\F_q}$.
Let
$B =B(\text{Coh}_{P^{1}})$
be the extended Ringel algebra of the category of coherent sheaves on
$P^1$.
The group
$\text{Pic} (X)$
consists of sheaves
$\Cal O(n), n \in \Z$.
Thus the algebra $B$ is obtained from the algebra
$R(\text{Coh}_{P^{1}})$
by adding two elements
$K =K_{\Cal O}$
and
$c=c_{\Cal O(1)}$.
The set 
$\text{Cusp}$
consists of one element: the trivial character
$\bold 1$
of
$\text{Pic} (P^{1}) = \Z$.
Thus there are only two generating functions

$$E(t) =\sum_{n\in\Z} [\Cal O(n)]t^n \quad \in \quad B[[t,t^{-1}]], $$
$$\psi(t) = \sum_{D\in \text{Div}^+(P^{\prime})}
 |\text{Aut} \, \Cal O_D|\cdot t^{\deg D}[\Cal O_D]
\quad  \in \quad 1 + tB[[t]]$$
where
$\text{Div}^+(P^1)$
is the set of effective divisors on
$P^1$
and
$\Cal O_D$
is the torsion sheaf
$\Cal O_{P^1}/\Cal O_{P^1} (-D)$.

Let
$R(\text{Bun}(P^1)) \subset B$
be the subalgebra generated by elements
$[V]$, 
where $V$ is a vector bundle. Let
$\Cal B \subset B$
be the subalgebra generated by
$[V]$
as above together with
$K,c$
and the coefficients of
$\psi (t)$.

\proclaim{(5.2.1) Theorem}
The algebra
$R(Bun(P^1))$
is isomorphic to
$U_q(\widehat{\bold n^+}) \subset U_q(\widehat{sl}_2)$,
and
$\Cal B$
is isomorphic to
$U_q(\widehat{\bold b^+})$.
\endproclaim

\demo{Proof}
The function
$\text{LHom}(\bold 1, \bold 1, t)$,
i.e., the zeta-function of
$P^1$, has the form

$$\zeta (t) = \frac 1 {(1-t)(1-qt)} .\tag 5.2.2$$
Thus the functional equation for the Eisenstein series
 (3.3.1), when brought to the polynomial form (3.5.4) gives:

$$(t_1-qt_2) E(t_1)E(t_2) = (qt_1-t_2)E(t_2)E(t_1), \tag 5.2.3$$
which is identical to the defining relation (5.1.2) for
$U_q(sl_2)$
(in which
$i=j=1, a_{11} =2$).
Thus the correspondence
$x^+_i(n) \mapsto[\Cal O(n)]$
gives a homomoprhism
$\gamma : U_q(\widehat {\bold n^+}) \rightarrow R(\text{Bun}(P^1))$.
To show that
$\gamma$
is an isomorphism, note that every vector bundle $V$ on
$P^1$
can be represented as
$V = \bigoplus_{i\in\Z} \Cal O(i)^{m_i}$.
Moreover, the filtration
$V^j = \bigoplus_{i \ge j} \Cal O(i)^{m_i}$
is defined intrinsically (this is the Harder-Narasimhan filtration,
 already discussed in (3.9)). It follows that

$$[V] = q^{-\frac 1 2 \sum_{i<j}m_im_j(i-j+1)} [\Cal O(b)^{m_b}] * 
\ldots * [\Cal O(a)^{m_a}]$$
where
$\{a,a+1, \ldots, b\}$
is any interval in
$\Z$
containing all $i$ such that
$m_i \ne  0$.
Note also that

$$[\Cal O(i)^m]*[\Cal O(i)^{m^{\prime}}] = q^{\frac {mm^{\prime}} 2} 
\bmatrix m+m^{\prime} \\ m \endbmatrix_q \Cal O(i)^{m+m^{\prime}},$$
which implies that
$[V]$
is actually a monomial in the
$[\Cal O(i)]$.
Thus
$\gamma$
is surjective. To see that
$\gamma$
is injective, i.e., that there are no further relations among the
$[\Cal O(i)]$
except those following from (5.2.2), remark that we just have shown
 that the monomials

$$[\Cal O(i_1)]^{m_1}[\Cal O(i_2)]^{m_2} \ldots [\Cal O(i_r)]^{m_r}, 
\quad i_1 > \ldots > i_r \tag 5.2.4$$
constitute a 
$\C$-basis in
$R(\text{Bun}(P^1))$.
On the other hand, the relations (5.2.3) (after passing to the
 coefficients) can be used to express any monomial as a linear
 combination of monomials (5.2.4). Since the  latter monomials are
 linearly independent in
$R(\text{Bun}(P^1))$,
there can be no further relations and so 
$\gamma$
is an isomorphism.
\enddemo

Let us now turn to the algebra
$\Cal B$.
The commutation relation (3.3.2) reads in our case (again, use (5.2.2)) as follows:

$$E(t_1)\psi(t_2) = \frac{1-q^{-1/2}t_2/t_1} {1-q^{3/2}t_2/t_1} \psi 
(t_2) E(t_1), \tag 5.2.5$$
while in
$U_q(\widehat{\bold b^+})$
we have, by (5.1.3):

$$F^+(t_1)\varphi^+(t_2) = \frac{q-c^{-1/2} t_2/t_1} {1-q \, c^{-1/2} 
t_2/t_1} \varphi^+(t_2) F^+(t_1) \tag 5.2.6$$
Thus the correspondence

$$\varphi^+(t) \mapsto K \, \psi (c^{-1/2}q^{-1/2}t), \quad F^+(t) \mapsto E(t)$$
defines a homomorphism
$U_q(\widehat{\bold b^+}) \rightarrow \Cal B$.
It is surjective by the definition of
$\Cal B$.
Its injectivity follows from the injectivity of the homomorphism
$\gamma : U_q(\widehat{\bold n^+}) \rightarrow R(\text{Bun}(P^1))$
above and the fact that the coefficients of
$\psi(t)$
are  algebraically independent over
$\C$.
This latter fact is proved by the same argument as given in (3.8.5): 
the coefficients of
$\log \psi(t)$,
being linearly independent and primitve with respect to a Hopf algebra
 structure, are algebraically independent. Theorem (5.2.1) is proved.

\vskip .3cm

\noindent{\bf  (5.2.7) Remark.}
It is instructive to compare Theorem 5.2.1 with results of Ringel
 [R1-3] and Lusztig [Lu1-3] on representations of quivers. Namely, let
$\Gamma$
be a quiver, i.e., a finite oriented graph without edges-loops. Let
$\Gamma -\text{mod}$
be the category of representations of
$\Gamma$
over
$\bold F_q$,
i.e., rules which associate to any vertex
$i\in \text{Vert}(\Gamma)$
a finite-dimensional
$\bold F_q$
-vector space
$V_i$,
and to every oriented edge
$i \overset e \to \longrightarrow j$
a linear map
$V_e : V_i \rightarrow V_j$.
This is an Abelian category satisfying our finiteness conditions
 (1.1), so its Ringel algebra
$R(\Gamma - \text{mod})$
is defined. For
$i \in \text{Vert}(\Gamma)$
let
$V(i)$
be the representation assigning
$\bold F_q$
to the vertex $i$ and 0 to other vertices. Let
$A = \parallel a_{ij}\parallel_{i,j \in Vert (\Gamma)}$
be the Cartan matrix associated to
$\Gamma$,
i.e.,
$a_{ii} = 2$
and
$a_{ij}$
is minus the number of edges (regardless of the orientation) 
joining $i$ and $j$, if
$i \ne j$.
Let
$\bold g_{\Gamma}$
be the Kac-Moody Lie algebra with the Cartan matrix $A$ and
 $U_q(\bold g_\Gamma)$ be its
$q$-quantization. This is (see, e.g.,  [Lu 1]) the algebra
 generated by symbols $e^\pm_i,  K_i$
(with $K_i$ being invertible)
subject to the relations:
$$ K_i e_j^\pm K_i^{-1} = q^{\pm a_{ij}/2} e_j, \quad [e^+_i, e^-_j] =
 \delta_{ij} {K_i-K_i^{-1}\over
q-q^{-1}},\leqno (5.2.8)$$
$$\sum_{l=0}^{1-a_{ij}} (-1)^l \left[ 1-a_{ij}\atop l\right]_q 
(e_i^\pm)^l e_j^\pm (e_i^\pm)^{1-a_{ij}-l}
=0.\leqno (5.2.9)$$
Let $U_q(N_{\Gamma})$ be the subalgebra generated by the $e_i^+$.
Then the subalgebra in
$R(\Gamma -\text{mod})$
generated by the 
$[V(i)]$,
is isomorphic to
$U_q(N_{\Gamma})$
so that
$[V(i)]$
corresponds to
$e^+_i$.

Now,
$\widehat{sl}_2$
is the Kac-Moody algebra associated to the Cartan matrix
$\pmatrix 2 & -2 \\ -2 & 2 \endpmatrix$,
and its Dynkin graph is
$A^{(1)}_1  = \{ \bullet \rightrightarrows \bullet\}$.
So the general results about representations of quivers are
 applicable to this case and realize
$U_q(N_{A^{(1)}_1}) \subset U_q(\widehat{sl}_2)$
inside
$R(A^{(1)}_1 -\text{mod})$.
So we have realizations of two ``nilpotent'' subalgebras,
$U_q(\widehat{\bold n^+}), U_q(N) \subset U_q(\widehat{sl}_2)$
in terms of Ringel algebras of two Abelian categories:
$\text{Coh}_{P^1}$
(sheaves on
$P^1_{F_q}$)
and
$A^{(1)}_1 - \text{mod}$.
One may wonder what these two categories have in common, 
and the answer is
 that their derived categories are equivalent. This is a
 particular case of a theorem of Beilinson [Be], as reformulated by Bondal 
[Bo] and Geigle-Lenzing [GL]. This is another indication of a deeper relation between
 derived categories and quantum group-like objects, see [X1].

\vskip .3cm

\noindent{\bf  (5.3) The case of a general curve $X$.}
Considering now the case of an arbitrary smooth projective curve
$X/\F_q$
we find, by comparing Theorem 3.5 with formulas of (5.1), that the algebra
$B(\text{Coh}_X)$
(or, rather, its subalgebra
$\Cal B$
defined in (3.8)) is analogous to the subalgebra
$U_q(\widehat{\bold b^+}) \subset U_q(\widehat {\bold g})$
for a huge Kac-Moody algebra
$\bold g$.
Let us sketch this analogy in the following table:

$$\matrix \format \l & \qquad   \l \\
\text{\bf  Theory of quantum affine algebras} & \qquad 
 \text{\bf  Theory of automorphic forms} \\
\text{  } \\
\text{The set of positive roots of
$\bold g$} & \qquad \text{The set $\text{Cusp}$ of cusp
 eigenforms} \\
\text{   }\\
\text{The entries of the Cartan matrix of
$\bold g$} & \qquad \text{The coeffcients of Rankin $L$-functions} \\
\text{   } & \qquad \text{$\text{LHom}(f,g,t)$} \\
\text{   }\\
\text{Symmetry of Cartan matrix} & \qquad \text{Functional
 equations of $L$-functions} \\
\text{   } \\
\text{Pointwise-uppertriangular subalgebra} & \qquad 
\text{The algebra $R(\text{Bun}(X))$, i.e., the}\\
\text{$U_q(\widehat{\bold n}^+)$} & \qquad \text{algebra
 of unramified forms} \\
\text{   }\\
\text{Pointwise-Cartan subalgebra} & \qquad \text{The algebr
a of classical Hecke}\\
\text{$U_q(\bold h [t])$} & \qquad \text{operators} \\
\text{   } \\
\text{Root decomposition of $U_q(\bold n^+)$} & \qquad
 \text{Spectral decomposition of the space}\\
\text{   } & \qquad \text{of automorphic forms} \\
\text{   } \\
\text{Components of $x^+_{i_1}(t_1) \ldots x^+_{i_n}(t_n)$} & 
\qquad \text{Eisenstein series} \\
\text{at a given basis vector of $U_q(\hat \bold g)$} \\
\text{   } \\
\text{Serre relations} & \qquad \text{?} \\
\text{   } \\
\text{?} & \qquad \text{Selberg trace formula} \\
\text{  } \\
\text{Full algebra $U_q(\widehat{ \bold g})$} & \qquad \text{?} \\
\text{   } \\
\text{Mc Kay correspondence} & \qquad \text{Langlands correspondence}
\endmatrix$$ 

\vskip .3cm

Let us comment a little on the last line. On the left, we have the McKay correspondence
  between finite subgroups $\Gamma\subset SL_2({\bold  C})$
and finite Dynkin diagrams of type A-D-E. According to McKay, one can construct,
for a given subgroup $\Gamma$, the {\it affinization} $\hat \Delta$ of the
corresponding Dynkin diagrams as follows. Vertices of $\hat \Delta$ are
isomorphism classes of irreducible representations of $\Gamma$. Two vertices
$[V]$ and $[W]$ are joined by an edge if $W$ enters into the decomposition of
$V\otimes {\bold  C}^2$, where ${\bold  C}^2$ is the restriction of the standard
representation of $SL_2({\bold  C})$.

On the right, we have the (conjectural)
Langlands correspondence between irreducible representations
of the unramified Galois group $\pi_1(X)$ and automorphic forms. Now, automorphic
forms are, as we have seen,  analogous to simple roots of a Lie algebra, i.e., vertices of some huge Dynkin diagram while the Galois group is analogous to a subgroup
$\Gamma\subset SL_2({\bold  C})$. The role of the standard representation ${\bold  C}^2$
is played by the algebraic closure $\overline k$ where $k$ is the global field of
functions on $X$. Since the role of entries of the Cartan matrix (describing edges
of the Dynkin diagram) is played in our situation by the $L$-functions,
this suggests that one can get some information about them by studying
how a Galois representation splits when tensor multiplied with the algebraic closure.
As was pointed out to me by D. Kazhdan, such a situation appears
in the theory of Hodge-Tate decompositions of $p$-adic Galois representations
of $p$-adic fields.

We have also included in this table several concepts whose
 counterparts are not immediately clear. Thus, Serre relations
 should correspond to ``extra'' functional equations for
 Eisenstein series, i.e., elements of the kernel of the map
$\pi_{\Cal B} : \tilde \Cal B \rightarrow \Cal B$
in (3.8.2), which was described in (3.8) only implictly,
 by means of the kernel of a suitable quadratic form on
$\tilde \Cal B$.
Note that the presence of such relations does not contradict
 the spectral decomposition theorem for automorphic forms [MW]
 which, as we already pointed out (3.8.4), is only concerned with 
the values
 of the (analytically continued) products
$E_{f_1} (t_1) \ldots E_{f_r} (t_r)$
for
$|t_i|=1$,
while the extra relations may look like distributions whose 
support does not meet the torus
$|t_i| =1$
at all. This is exactly the case with Serre relations (5.1.7).
 Namely, if we allow ourselves to divide equalities involving
 generating functions by any polynomials in
$t_1, \ldots, t_n$,
then (5.1.7) would ``follow'' from (5.1.2), as one can easily check
by ``bringing"  the LHS of (5.1.7) to a normal form in which  the factor
$F_j^\pm(s)$ goes last in every monomial (of course, such a manipulation is illegal!).
 Moreover, in this particular case it is enough to divide by polynomials of the form
$(t_i-q^{\alpha}t_j), \alpha \in \frac 1 2 \Z-\{0\}$,
so the  need for (5.1.7) cannot be observed after restriction
 to the torus
$|t_i|=1$
where such polynomials do not vanish. However, there is no way
 to deduce (5.1.7) from the quadratic relations obtained by 
expanding (5.1.2). This just means that in the algebra defined by (5.1.2)
 alone the right hand side of (5.1.7) is a distribution with support
 of positive codimension.

The analog of the Selberg formula in the theory of affine Lie algebras
 is completely unclear at the moment. According to our analogy, 
it should be a statement about the trace of a Cartan element
 acting on the enveloping algebra of the uppertriangular subalgebra.

As for the analog of the
$U_q(\widehat{\bold g})$,
we will construct such an analog in the next section by following Drinfeld's
 quantum double construction and using Green's comultiplication.

\newpage

\centerline {\bf  \S 6. The quantum double of the algebra of automorphic forms.} 

\vskip 1cm

\noindent{\bf  (6.0) Motivation.}
Our aim in this section is to construct the automorphic analog
 of the full quantum affine algebra
$U_q(\widehat{\bold g})$ out of the Hopf algebra $B(\Cal A)$
 which should be, morally, just one half
of it. This will be done by the Drinfeld double construction.
 For any Hopf algebra $\Xi$
 its Drinfeld double [Dr 3] [CP] is the tensor product 
$\Xi\otimes_{\bold C} \Xi^*$
with certain twisted multiplication (see below). 
The reason for using this construction is that it is known
 to solve, in a certain sense, the following
model problem: recover the quantum Kac-Moody algebra
 $U_q(\bold g)$ from its Hopf subalgebra
$U_q(\bold b^+)$. Recall (5.2.7) that $U_q(\bold g)$
 is generated by symbols $e_i^\pm, K_i$.
The subalgebra $U_q(\bold b^+)$ is, by definition,
 generated by $e_i^+$ and $K_i$. The comultiplication
in it has the form:
$$\Delta(K_i) = K_i\otimes K_i, \quad \Delta(e_i^+) = 
1\otimes e_i^+ + e_i^+\otimes K_i.
\leqno (6.0.1)$$
However, the recovery of $U_q(\bold g)$ from the $U_q(\bold b^+)$ 
by this method
is not so straightforward. To help motivate the constructions of
 this section,
let us recall the principal features of (and subtleties involved in) 
this recovery.

\vskip .2cm

\noindent {\bf  (6.0.2)} The subalgebras $U_q(\bold n^\pm)$
 generated by the $e_i^\pm$,
are in natural duality, so adding $(U_q(\bold n^+))^*$ as
 a part of $(U_q(\bold b^+))^*$
supplies  the missing subalgebra $U_q(\bold n^-)$.

\vskip .1cm

\noindent {\bf  (6.0.3)}  However, we get also the dual of the ``Cartan"
subalgebra $U_q(\bold h)$ generated by the $K_i$. So by formally applying the double
construction to $U_q(\bold b^+)$, we get an algebra bigger that $U_q(\bold g)$
because of this dual. So in order to get $U_q(\bold g)$, we need to impose certain
identifications in the double. (The construction of the restricted double later in this
section does just that, in the context of Hall algebras.)

\vskip .1cm

\noindent {\bf  (6.0.4)} In addition, the Hopf algebra
 $U_q(\bold h)$ is not self-dual.
It is the algebra of functions on $T$, the maximal torus
 in the Lie group corresponding to $\bold g$,
and the dual algebra will be algebra of functions on
 $\widehat{T}$ (the lattice  of characters of $T$) with
 pointwise multiplication. Formal dualization of $U_q(\bold h)$
 with respect to the basis
of monomials in the $K_i$ gives a basis of delta-functions on
 $\widehat{T}$, in particular,
the unit element in the algebra of functions on $\widehat{T}$ 
will be represented by an infinite
sum of such delta-functions. If the Cartan matrix $A$ is non-degenerate,
 one can realize
elements of $U_q(\bold h)$ by certain infinite sums in the
 completion of the algebra of functons
on $\widehat{T}$. After this one can perform the identification 
mentioned in (6.0.3).

\vskip .1cm

\noindent {\bf  (6.0.5)} If the Cartan matrix is degenerate 
(as is the case with affine algebras),
then one has to add to $U_q(\bold h)$ some essentially new
 elements from $U_q(\bold h)^*$
which pair with the center (the kernel of the Cartan matrix)
 in a nondegenerate way.
In the affine case the new element (which is essentially one)
has important geometric meaning: it corresponds to
the ``rotation of the loop" and is responsible for the appearance 
of elliptic functions
in the theory of representations of affine algebras.  

\vskip .2cm

The general problem of constructing the Drinfeld double for the Hall-Ringel algebra
of a category of homological dimension 1 was studied in the recent paper of J. Xiao [X2].
In the analysis below we complement Xiao's approach by adding an intermediate step,
the so-called Heisenberg double, introduced in [AF] [ST]. The use of  the Heisenberg
double allows one to split the calculations into two more manageable steps.

\vskip .3cm

\noindent {\bf  (6.1) Heisenberg and Drinfeld doubles: generalities.}
We start by recalling some basic 
properties of the doubles. Since Hall algebras come with a natural
 basis, we will use the coordinate-dependent approach, as in [Kas]
(see [ST] for a more standard exposition).
 In this reminder we will ignore subtleties related to dualization of
 infinite-dimensional spaces, since we will pay due attention
 to them in the concrete situations in which we will use the doubles. 

\vskip .2cm 

Let $\Xi$ be a Hopf algebra over $\C$, with comultiplication
 $\Delta$, counit $\epsilon$ and
antipode $S$. Let $\{e_i\}, i\in I,$ be a basis of $\Xi$ and
 $\{e^i\}$ be the dual basis of $\Xi^*$. Let
us write the structure constants of $\Xi$ with respect to our basis:
$$e_ie_j = \sum_k m_{ij}^k e_k, \quad \Delta(e_k) = \sum_{i,j}
 \mu_k^{ij} e_i\otimes e_j, \leqno (6.1.1)$$
$$1 = \sum_i \epsilon^i e_i, \quad \epsilon(e_i) = \epsilon_i, 
\quad S(e_i) =
\sum_j S_i^j e_j, \quad S^{-1}(e_i) = \sum_j \sigma_i^j e_j.$$
As usual, $\Delta$ makes $\Xi^*$ into an algebra by
$$e^ie^j = \sum_k \mu_k^{ij} e^k. \leqno (6.1.2)$$
The Heisenberg double $HD(\Xi)$ is, by definition, the algebra
 generated by the symbols
$Z_i, Z^i, i\in I$, subject to the relations:
$$Z_iZ_j = \sum_k m_{ij}^k Z_k, \quad Z^iZ^j = \sum_k \mu_k^{ij} Z^k, 
\leqno (6.1.3)$$
$$Z_iZ^j = \sum_{a,b,c} m_{ab}^j \mu_i^{bc} Z^a Z_c. \leqno (6.1.4)$$
Thus the map $\Xi^* \otimes_\C \Xi \rightarrow HD(\Xi)$ given by
 the multiplication, i.e., by
$e^i\otimes e_j\mapsto Z^iZ_j$, is an isomorphism of vector spaces:
 any other product can be brought by
(6.1.4) to the normal form in which the $Z^i$ stand on the left. 

\vskip .2cm

There is a certain asymmetry in the definition of the
 Heisenberg double: replacing $\Xi$ with $\Xi^*$ gives a 
different algebra. We will denote $HD(\Xi^*)$ by $\check HD(\Xi)$.
 Explicitly, it is generated by the symbols
$\check Z_i, \check Z^i, i\in I$ subject to the relations:
$$\check Z_i \check Z_j = \sum_k m_{ij}^k \check Z_k, \quad \check 
Z^i \check Z^j = \sum_k \mu_k^{ij} \check Z^k, \leqno (6.1.3')$$
$$\check Z^i \check Z_j = \sum_{a,b,c} \mu_j^{ab} m_{bc}^i \check
 Z_a \check Z^c \leqno (6.1.4')$$
The Drinfeld double $DD(\Xi)$ is the algebra generated by the symbols
 $W_i, W^i, i\in I$,
subject to the relations
$$W_iW_j = \sum_k m_{ij}^k W_k, \quad W^iW^j = \sum_k \mu_k^{ik} W^k,
 \leqno (6.1.5)$$
$$\sum_{a,b,c} \mu_i^{ab} m_{bc}^j W_aW^c = \sum_{a,b,c}
 m^j_{ab}\mu_{i}^{bc} W^aW_c. \leqno (6.1.6)$$

 By using the antipode one can transform (6.1.6) into a
 relation giving
a kind of normal form for elements of $DD(\Xi)$:
$$W_iW^j = \sum_{a,b,c,d,e,f,g} \sigma_b^a m_i^{bc} m_c^{de}
 \mu_j^{af}\mu_f^{ge} W^gW_d.
\leqno (6.1.7)$$
This is the standard description [Dr 3] [CP] of the Drinfeld double. 

The algebra $DD(\Xi)$ is a Hopf algebra with respect to the
 comultiplication given by
$$\Delta(W_k) = \sum_{i,j} \mu_k^{ij} W_i\otimes W_j, \quad
 \Delta(W^k) = \sum_{i,j} m_{ji}^k W^i\otimes W^j, \leqno (6.1.8)$$
(note the transposition in the second formula), the counit
$$\epsilon(W_k) = \epsilon_k, \quad \epsilon(W^k) = \epsilon^k, \leqno (6.1.9)$$
and the antipode given by
$$S(W_i) = \sum_j S_i^jW_j, \quad S(W^j) = \sum_i S_i^jW^i.\leqno (6.1.10).$$
The algebras $HD(\Xi), \check HD(\Xi)$ are not Hopf algebras.
 However, there is the following
result due to R. Kashaev [Kas].

\proclaim {(6.1.11) Proposition}  The correspondence
$$W_k\rightarrow \sum_{i,j} \mu_k^{ij} Z_i\otimes \check Z_j, 
\quad W^k \rightarrow \sum_{i,j} m_{ji}^k Z^i\otimes \check Z^j$$
defines an embedding $\kappa: DD(\Xi) \hookrightarrow HD(\Xi)
\otimes HD(\Xi^*)$.\endproclaim

It is this embedding that makes the Heisenberg doubles important to
 us.

\vskip .2cm

\demo{Proof} Let us first prove the existence of the homomorphism
 $\kappa$ and then its injectivity.
To prove the existence, we need to check that the relations (6.1.5-6)
 defining $DD(\Xi)$ are satisfied for
the proposed images of $W_i, W^i$. For (6.1.5) it is obvious since
 $\Delta: \Xi\rightarrow \Xi\otimes\Xi$ and
$m^*: \Xi^*\rightarrow \Xi^*\otimes\Xi^*$ (the map dual to the
 multiplication) are algebra homomorphisms. To verify
(6.1.6) denote by $L$ and $R$ its left and right hand sides. Then:
$$L = \sum_{a,b,c} \mu_i^{ab} m_{bc}^j\, \kappa( W_a) \kappa(W^c) =
 \sum_{a,b,c,p,q,r,s} \mu_i^{ab} m_{bc}^j \mu_a^{pq} m_{sr}^c \, Z_pZ^r 
\otimes \check Z_q \check Z^s =$$
$$= \sum_{a,b,c,p,q,r,s\atop u,v,w} \mu_i^{ab} m^j_{bc} \mu_a^{pq}
 m_{sr}^c m_{uv}^r \mu_p^{vw}\,
 Z^uZ_w\otimes \check Z_q \check Z^s = \sum_{b,s,u,w,q, v}
 \mu_i^{vwqb} m_{bsuv}^j \, Z^uZ_w\otimes \check Z_q \check Z^s,$$
where $\mu_i^{vwqb}$ and  $m_{bsuv}^j$ are the structure constants
 for the 3-fold comultiplication and multiplication. Similarly,
$$R = \sum_{a,b,c} m^j_{ab} \mu_i^{bc} \kappa(W^a)\kappa(W_c) = 
\sum_{a,b,c,d,p,q,r,s} m^j_{ab} \mu_i^{bc} m_{qp}^a \mu_c^{rs}
Z^pZ_r \otimes \check Z^q \check Z_s =$$
$$=\sum_{a,b,c,d,p,q,r,s\atop u,v,w} m_{ab}^j \mu_i^{bc} m_{qp}^a 
\mu_c^{rs} \mu_s^{uv} m_{vw}^q \, Z^pZ_r\otimes \check Z_u
\check Z^w = \sum_{b,r,u,w,p,v} \mu_i^{bruv} m_{vwpb} \, Z^pZ_r
 \otimes \check Z_u\check Z^w,$$
which is the same as L up to renaming the summation indices. So
 we indeed have the claimed homomorphism $\kappa$.
To see that $\kappa$ is injective, note that $\Delta: \Xi\rightarrow 
\Xi\otimes\Xi$ and $m^*: \Xi^*\rightarrow \Xi^*\otimes\Xi^*$ 
are injective because of the unit and counit in $\Xi$. Since $DD(\Xi) 
= \Xi\otimes \Xi^*$ as a vector space, and $\kappa|_\Xi = \Delta$ while
 $\kappa|_{\Xi^*}$ is the composition of $m^*$ and the permutation,
 $\kappa$ is injective. \enddemo

\vskip .2cm

\noindent{\bf  (6.2) The restricted Heisenberg doubles of the Ringel algebra.}
We now specialize to the case when $\Xi = B(\Cal A) = \C[\Cal K _0 \Cal A]
\otimes R(\Cal A)$ is the extended Ringel algebra of an Abelian category
$\Cal A$ satisfying the assumptions of (1.4) (we assume, in particular, that
$\Cal A$ has homological dimension 1). The bilinear form $(\alpha|\beta)$
on $\Cal K_0\Cal A$ may be degenerate, and we denote by $I$ its kernel.
Since $\Q^*$, the target of our form, has no torsion, $I$ has the following
property: if $n\alpha\in I$ for some $n\in \Z$, $\alpha\in \Cal K_0\Cal A$,
then $\alpha\in I$. This implies that $I$ has a complement, a subgroup
$J\subset \Cal K_0\Cal A$ such that $\Cal K_0\Cal A = I\oplus J$. In the
sequel we fix such a complement and denote $\pi_I, \pi_J$ the projections
to $I$ or $J$ along the other summand. 

\vskip .2cm

A natural basis of $\Xi = B(\Cal A)$ is given by $e_{\alpha A} = K_\alpha [A]$,
$\alpha\in \Cal K_0\Cal A$, $A \in \Cal A$. In this basis, the multiplication
and comlutiplication have the form:
$$e_{\alpha, A}\cdot e_{\beta, B} = (\bar A|\beta) \langle B,A\rangle \sum_C
g_{AB}^C e_{\alpha+\beta, C}, \leqno (6.2.1)$$
where $g_{AB}^C$ is the same as in (1.2),
$$\Delta(e_{\gamma, C}) = \sum_{A,B} \langle B,A\rangle
{ |\text{Aut}(A)|\cdot |\text{Aut}(B)| \over |\text{Aut}(C)|} g_{AB}^C
e_{\gamma, A}\otimes e_{\gamma+\bar A, B}.\leqno (6.2.2)$$
Thus the structure constants are given by
$$m_{\alpha A, \beta B}^{\gamma C} = \left\{ \aligned &(\bar A|\beta)
 \langle B,A\rangle g_{AB}^C, \quad \text{if} \quad \gamma = \alpha+\beta \\
&0, \quad \text{if} \quad  \gamma\neq \alpha+\beta  \endaligned \right. \tag 
6.2.3$$
$$\mu_{\gamma C}^{\alpha A, \beta B} = \left\{
\aligned & \langle B,A\rangle  { |\text{Aut}(A)|\cdot |\text{Aut}(B)|
 \over |\text{Aut}(C)|} g_{AB}^C, \quad \text{if}\quad \alpha=\gamma,
 \beta=\gamma+\bar A,\\
&0, \quad\text{otherwise}\endaligned\right.\tag 6.2.4$$
This gives the multiplication of the dual generators in $HD(\Xi)$ in 
the form
$$ Z^{\alpha A}Z^{\beta B} = \left\{
\aligned & \langle B,A\rangle \sum_C { |\text{Aut}(A)|\cdot 
|\text{Aut}(B)| \over |\text{Aut}(C)|} g_{AB}^C Z^{\alpha C},
 \quad \text{if}\quad
 \beta=\gamma+\bar A,\\
&0, \quad\text{otherwise}\endaligned\right.\tag 6.2.5$$
In particular, the $Z^{\alpha 0}$ are orthogonal idempotents.
 They form the algebra dual to the coalgebra $\C[\Cal K_0\Cal A]$, 
i.e., the algebra of functions on $\Cal K_0\Cal A$ with pointwise
 multiplication. More precisely, $Z^{\alpha 0}$
corresponds to the function equal to 1 at $\alpha$ and to 0 elsewhere.
 In particular, the unit element is represented by the infinite sum $\sum_\alpha
Z^{\alpha 0}$.

 This shows that we need to work with certain infinite sums lying in 
the completion of the
algebra $HD(B(\Cal A))$ generated by the $Z_{\alpha A}, Z^{\alpha A}$. 
For $A \in\text{Ob}(\Cal A)$ and a homomorphism $\chi: \Cal K_0\Cal A
\rightarrow \C^*$ we introduce the elements
$$Z^-_A = \sum_{\alpha\in \Cal K_0\Cal A}{ Z^{\alpha A}\over
 |\text{Aut}(A)|},
\quad K^\chi = \sum_{\Cal K_0\Cal A} \chi(\alpha) Z^{\alpha 0}.
\leqno (6.2.6)$$
We set also $K_\alpha = Z_{\alpha 0}, Z^+_A = Z_{0 A}$.

If $\alpha\in \Cal K_0\Cal A$, we let $\chi_\alpha: \Cal K_0\Cal A\rightarrow
\Q^*$ be the character taking $\beta$ to $(\alpha|\beta)$, and write
$K^\alpha$ for $K^{\chi_\alpha}$. Note that we have the following identities:
$$Z^\pm_A Z^\pm_B = \sum_C g_{AB}^C Z^\pm_C, \tag 6.2.7$$
$$Z^\pm_A K_\beta = (\bar A|\beta)^{\pm 1} K_\beta Z^\pm_A,
\quad Z^\pm_A K^\chi = \chi(\bar A)^{\mp 1} K^\chi Z^\pm_A,
\tag 6.2.8$$
$$K^\chi K^{\chi'} = K^{\chi\chi'}, \quad K^\chi K_\alpha = K_\alpha K^\chi, \tag 6.2.9$$
which show that we have two copies of the Ringel algebra $R(\Cal A)$
inside the completion of $HD(B(\Cal A))$.

\proclaim{ (6.2.10) Proposition} The subspace in the completion 
of $HD(B(\Cal A))$ spanned by elements of the form $Z^-_A K_\alpha K^\chi Z^+_B$,
is an algebra, i.e., the product of any two such elements is a
 linear combination of finitely many elements of this form. \endproclaim

We will call this algebra the algebraic Heisenberg double of
 $B(\Cal A)$ and denote it $HD^{alg}(B(\Cal A))$. 

\vskip .2cm

To prove (6.2.10), note that the identities (6.2.7-9) give an 
almost complete recipe for multiplying any  $Z^-_A K_\alpha K^\chi Z^+_B$
with any $Z^-_{A'}K_{\alpha'} K^{\chi '} Z^+_{\beta'}$, except for the
product $Z^+_B Z^-_{A'}$ which we need to express through $Z^-_M Z^+_N$,
In order to find this expression, denote, for any objects $A,B,M,N$ of
$\Cal A$, by $g_{AB}^{MN}$ the orbifold number of exact sequences
$$0\rightarrow M\buildrel u\over\rightarrow B\buildrel \varphi\over
\rightarrow
A\buildrel v\over\rightarrow N\rightarrow 0 \tag 6.2.11$$
modulo automorphisms of $A$ and $B$, i.e., (see the proof of Lemma 4.2)
the total number of such sequences divided by $|\text{Aut}(A)|\cdot
|\text{Aut}(B)|$. Proposition 6.2.10 follows from the next fact which is
also useful by itself.

\proclaim{(6.2.12) Proposition} In the completion of $HD(B(\Cal A))$ we
 have the identities:
$$Z^+_A Z^-_B = \sum_{M,N} g_{AB}^{MN} \langle \bar N - \bar M, \, 
\bar B - \bar M\rangle K_{\bar B - \bar M} Z^-_M Z^+_N.$$ \endproclaim

\demo{Proof} Let $\Cal F_{AB}^{MN}$ be the set of all exact sequences
(6.2.11). Also, for any three objects $A,B,C$ let $\Cal E_{AB}^C$ be the
set of exact sequences
$$0\rightarrow A\rightarrow C\rightarrow B \rightarrow 0.$$
Thus
$$g_{AB}^{MN} = { |\Cal F_{AB}^{MN}|\over |\text{Aut}(A)|\cdot
|\text{Aut}(B)|}, \quad g_{AB}^C = {|\Cal E_{AB}^{C}|\over |\text{Aut}(A)|\cdot
|\text{Aut}(B)|}.$$
Notice now that
$$\Cal F_{AB}^{MN} \quad = \quad \coprod_{L\in \text{Ob}(\Cal A)/\text{Iso}}
\bigl( \Cal E_{MB}^L \times \Cal E_{LA}^N\bigr) \bigl/ \text{Aut}(L),
\tag 6.2.13$$
with $\text{Aut}(A)$ acting freely. This just means that any
 long exact sequence (6.2.10) can be split into two short sequences with $L =
 \text{Im}(\varphi)$.
From  the cross-symmetry relations (6.1.4) in the Heisenberg double
and the explicit form of the structure constants given above, we find:
$$Z_{0A} Z^{\beta B} = \sum_{M,L,N} \langle L,M\rangle \cdot \langle N,L\rangle
\cdot {|\text{Aut}(L)|\cdot |\text{Aut}(N)|\over |\text{Aut}(A)|}\cdot$$
$${ |\Cal E_{ML}^B|\over |\text{Aut}(L)|\cdot |\text{Aut}(M)|}\cdot
{ |\Cal E_{LN}^A|\over |\text{Aut}(L)|\cdot |\text{Aut}(N)|}\cdot Z^{\beta M}
Z_{\bar L, N}=$$
$$= \sum_{M,N} \langle \bar B - \bar M, \bar M\rangle\cdot \langle \bar N,
 \bar B-\bar M\rangle\cdot {|\Cal F_{AB}^{MN}|\over |\text{Aut}(A)|\cdot
|\text{Aut}(M)|} Z^{\beta M} Z_{\bar B - \bar M, N},$$
where we used (6.2.12) and the identity $\bar L = \bar B -\bar M$
holding whenever $\Cal E_{ML}^B \neq\emptyset$. From this we deduce that
$$Z_{0A} \cdot {Z^{\beta B}\over |\text{Aut}(B)|} = 
\sum_{M,N}  \langle \bar B - \bar M, \bar M\rangle\cdot \langle \bar N,
 \bar B-\bar M\rangle (\bar M|\bar B - \bar M)^{-1} g_{AB}^{MN} K_{\bar B - \bar M}
{ Z^{\beta M}\over  |\text{Aut}(M)|} \cdot Z_{0,N},$$
whence the statement of the proposition.\enddemo

Finally, we define the algebra $\text{Heis}(\Cal A)$, called the
 restricted Heisenberg double
of $B(\Cal A)$, by imposing the following identifications in
 $HD^{alg}(B(\Cal A))$:
$$K^\alpha = K^{-1}_{\pi_J(\alpha)}, \leqno (6.2.14)$$
where $\pi_J$ was defined in (6.1).
Thus, in the case when the form $(\alpha|\beta)$ is non-degenerate, 
this identifies each $K^\chi$ with some product of (complex) powers
 of the $K_\alpha$, while in the  degenerate case we have also elements
$K^\chi$ for characters $\chi$ of $\Cal K_0 \Cal A = I\oplus J$ trivial 
on $J$, i.e., characters of
$I$. These elements are analogs of the ``rotation of the loop" generators 
in the (quantum) 
Kac-Moody algebras.

\vskip .2cm

The other , ``checked", Heisenberg double $\check HD(B(\Cal A))$ can be
treated in a similar way. Its
generators are denoted $\check Z_{\alpha A}, \check Z^{\alpha A}$, and
 we set
$$\check Z^+_A = \check Z_{0A}, \quad \check Z^-_A = \sum_{\alpha\in
 \Cal K_0\Cal A}{ \check Z^{\alpha A}
\over |\text{Aut}(A)|}, \quad \check K_\alpha = \check Z_{\alpha 0},
 \quad \check K^\chi = 
\sum_{\alpha\in \Cal K_0\Cal A} \chi(\alpha)\check Z^{\alpha 0}. \leqno
 (6.2.15)$$
The subspace spanned by the $\check Z^-_A \check K_\alpha \check K^\chi
 \check Z^+_B$ is again 
a subalgebra, denoted $\check HD^{alg}(B(\Cal A))$. Its generators satisfy
 the same identities as in (6.2.7-9), while (6.2.12) is now replaced by
$$\check Z^-_A\check Z^+_B = \sum_{M,N} \langle \bar B - \bar M, \bar M\rangle \cdot
\langle \bar N, \bar B - \bar M\rangle \cdot g_{AB}^{MN} K^{\bar B - \bar M}
\check Z^+_M \check Z^-_N. \leqno (6.2.16)$$
As before, we define the restricted double ${\text{Heis}}^\vee
 (\Cal A)$ by quotienting
$\check HD^{alg}(B(\Cal A))$ by the relations $\check K^\alpha = 
\check K^{-1}_{\pi_J(\alpha)}$.

\vskip .3cm

\noindent {\bf  (6.3) The restricted Drinfeld double of the Ringel
 algebra.} As in the case of Heisenberg
doubles, let $W_{\alpha A}, W^{\alpha A}$ be the generators of
 $DD(B(\Cal A))$ corresponding to
the basis vector $e_{\alpha A} = K_\alpha [A]$. Similarly to
 the above, we set
$$W^+_A = W_{0A}, \quad W^-_A = \sum_{\alpha\in
 \Cal K_0\Cal A}{W^{\alpha A}\over|\text{Aut}(A)|},
\quad K_\alpha = W_{\alpha 0}, \quad K^\chi = \sum_{\alpha\in \Cal
 K_0\Cal A}\chi(\alpha)W^{\alpha 0}.$$
We also write $K^\alpha = K^{\chi_\alpha}$, where $\chi_\alpha(\beta) 
= (\alpha|\beta)$.
We denote by $DD^{alg}(B(\Cal A))$ the subspace in the completion of
 $DD(B(\Cal A))$ spanned by elements
of the form $W^+_A K_\alpha K^\chi W^-_B$.

\proclaim {(6.3.1) Proposition} (a) $DD^{alg}(B(\Cal A))$ is an algebra.

(b) $DD^{alg}(B(\Cal A))$ is a Hopf algebra with respect to the
 comultiplication given on generators by
$$\Delta(W^+_A) = \sum_{A^\prime\i A} \langle A/A', A'\rangle 
{|\text{Aut}(A')|\cdot |\text{Aut}(A/A')|\over |\text{Aut}(A)|}
 W^+_{A'}\otimes K_{\bar A'}W^+_{A/A'}
\leqno (6.3.2)$$

$$\Delta(W^-_A) = \sum_{A^\prime\i A} \langle A/A', A'\rangle 
{|\text{Aut}(A')|\cdot |\text{Aut}(A/A')|\over |\text{Aut}(A)|}
 W^-_{A/A'}K^{\bar A'} \otimes W^-_{A'},
\leqno (6.3.3)$$

$$\Delta(K_\alpha) = K_\alpha \otimes K_\alpha, \quad \Delta(K^\chi) = 
K^\chi \otimes K^\chi.
\leqno (6.3.4)$$

the counit given by

$$\epsilon(K_\alpha) = 1, \quad \epsilon(W^\pm_\A)=0, \quad \epsilon(K^\chi) = 1,
 \leqno (6.3.5)$$

and the antipode given by

$$S(K_\alpha) = K_\alpha^{-1} = K_{-\alpha}, \quad S(K^\chi) = K^{\chi^{-1}},
 \leqno (6.3.6)$$

$$ S(W^+_A) = \sum^{\infty}_{n=1} (-1)^n \sum_{A_0 \subset \ldots \subset A_n=A}
 \prod^n_{i=1} \langle A_i/A_{i-1},A_{i-1}\rangle \frac{\prod^n_{j=0}
 |\text {Aut} (A_j/A_{j-1})|} {|\text {Aut} (A)|} \cdot 
\leqno (6.3.7)$$
$$\cdot W^+_{A_0} W^+_{A_1/A_0} \ldots W^+_{A_n/A_{n-1}} \cdot  K^{-1}_{\bar A},$$

$$ S(W^-_A) = \sum^{\infty}_{n=1} (-1)^n \sum_{A_0 \subset \ldots
 \subset A_n=A} \prod^n_{i=1} \langle A_i/A_{i-1},A_{i-1}\rangle
 \frac{\prod^n_{j=0} |\text {Aut} (A_j/A_{j-1})|} {|\text {Aut} (A)|} \cdot 
\leqno (6.3.8)$$
$$ K^{-\bar A} \cdot 
W^-_{A_n/A_{n-1}}\ldots W^-_{A_1/A_0}W^-_{A_0}.$$

(c) The correspondence
$$W_A^+\mapsto \sum_{A^\prime\i A} \langle A/A', A'\rangle 
{|\text{Aut}(A')|\cdot |\text{Aut}(A/A')|\over |\text{Aut}(A)|}
 Z^+_{A'}\otimes \check K_{\bar A'}\check Z^+_{A/A'}, \leqno (6.3.9)$$

$$W^-_A \mapsto \sum_{A^\prime\i A} \langle A/A', A'\rangle 
{|\text{Aut}(A')|\cdot |\text{Aut}(A/A')|\over |\text{Aut}(A)|}
 Z^-_{A/A'}K^{\bar A} \otimes Z^-_{A'},
\leqno (6.3.10)$$

$$K_\alpha \mapsto K_\alpha \otimes\check K_\alpha, \quad K^\chi
 \mapsto K^\chi\otimes \check K^\chi \leqno (6.3.11)$$

defines an algebra homomorphism 
$$\kappa: DD^{alg}(B(\Cal A)) \hookrightarrow HD^{alg}(B(\Cal A))
 \otimes \check HD^{alg}(B(\Cal A)).$$
\endproclaim

\demo {Proof} (a) This follows from (6.1.7) together with the identities
$$W^\pm_A W^\pm_B = \sum_C \langle B,A\rangle g_{AB}^C W^\pm_C, 
\quad W^\pm_A K_\beta = 
(\bar A|\beta)^{\mp 1}K_\beta W^\pm_A, \quad W^\pm_A K^\chi =
 \chi(\bar A)^{\mp 1} K^\chi W^\pm_A,$$
$$K_\alpha K^\chi = K^\chi K_\alpha,\quad K^\chi K^{\chi'} =
 K^{\chi \chi'},$$
which are verified in the same way as (6.2.7-9). The statements (b),
 (c) follow by a straightforward
application of definitions, checking that they indeed make sense on
 the elements given by infinite sums 
we consider and applying the fact that, in general, $DD(\Xi)$ is a
 Hopf algebra and the map $\kappa$ from
proposition 6.2.12 is a homomorphism.

\enddemo

We now define the restricted Drinfeld double $U(\Cal A)$ by imposing 
in $DD^{alg}(B(\Cal A))$
the relations $K^\alpha = K^{-1}_{\pi_J(\alpha)}$. The notation is 
chosen to suggest the analogy
 with the quantum universal enveloping algebra $U_q(\bold g)$ of a
 Lie algebra $\bold g$.
 In fact, when $\Cal A = \Gamma-\text{mod}$
is the category of $\F_q$-representations of a Dynkin quiver $\Gamma$,
 then $U(\Cal A) = U_q(\bold g)$,
where $\bold g$ is the semisimple Lie algebra corresponding to $\Gamma$.

\proclaim{(6.3.12) Proposition} (a) The Hopf algebra structure on
 $DD^{alg}(B(\Cal A))$ descends to $U(\Cal A)$.

(b) The homomorphism $\kappa$ defined in (6.3.9-11) descends to a 
homomorphism
$$\varkappa: U(\Cal A) \hookrightarrow \text{Heis}(\Cal A)\otimes 
{\text{Heis}}^\vee
(\Cal A).$$\endproclaim

The proof is straightforward, by checking the compatibility of the
 relations $K^\alpha = K^{-1}_{\pi_J(\alpha)}$ with the Hopf algebra
 structure and with $\kappa$.

\vskip .3cm

\noindent{\bf  (6.4) The Heisenberg double of the algebra of
 automorphic forms.} We now further
 specialize to the case when $\Cal A$ is the category of 
coherent sheaves on a curve $X/\F_q$ 
and we keep
the notations and assumptions of \S 2-3. Recall that we
 have the exact sequence
$$0\rightarrow \text{Pic}(X)\buildrel i\over\rightarrow 
\Cal K_0(X)\buildrel \text{rk}\over
\rightarrow \Z \rightarrow 0, \leqno (6.4.1)$$
where rk is the generic rank homomorphism, and $i(L) = 
\bar L - 1$. The image of $i$ is precisely $I$, the kernel 
of the form $(\alpha|\beta)$. A natural splitting of the
 sequence (6.4.1) and thus a natural choice of a complement $J$ 
to $I$, is provided by the maps
$$\det: \Cal K_0(X)\rightarrow \text{Pic}(X), \quad \varepsilon:
 \Z \rightarrow \Cal K_0(X)$$
given by $\det(\bar V) = \overline{\bigwedge^{\text{rk}(V)}
 (V)}$, if $V$ is a vector bundle, and by 
$\varepsilon(1) = \bar{\Cal O}$. Thus the Cartan elements of
 the restricted Heisenberg double are:
$$K = K_{\bar\Cal O}, \quad c_L = K_{\bar L - \bar \Cal O}, \,
 L\in \text{Pic}(X), \quad
d_\chi = K^{\chi(\det)}, \, \chi: \text{Pic}(X) \rightarrow \C^*.$$
The elements $c_L$ are central, while the $d_\chi$ are the 
``rotation of the loop" elements. 
They commute with the other generators by the rule:
$$d_\chi W^\pm_{\Cal F} = \chi (\det(\bar\Cal F))^{\mp 1}
 W^\pm_{\Cal F} d_\chi, \quad
\Cal F \in \text{Coh}_X.$$
In particular, for $\lambda\in\C^*$ we will denote $d(\lambda) 
= d_{\chi_\lambda}$, where
$\chi_\lambda: \text{Pic}(X)\rightarrow \C^*$ is the character 
given by $\chi_\lambda(L) = \lambda^
{\deg (L)}$. This element $d(\lambda)$ is a familiar degree operator:
$$ d(\lambda) W^\pm_V = \lambda^{\mp \deg(V)} W^\pm_V d(\lambda), 
\quad V\in \text{Bun}(X),$$
$$d(\lambda) W^\pm_{\Cal F} = \lambda^{\pm h^0(\Cal F)} W^\pm_{\Cal F}
 d(\lambda), \quad \Cal F \in \text{Coh}_{0,X}.$$
We write Heis for $\text{Heis}(\text{Coh}_X)$ and ${\text{Heis}}^\vee$ for 
${\text{Heis}}^\vee(\text{Coh}_X)$. The Cartan elements of
 ${\text{Heis}}^\vee$
will be denoted by $\check K, \check c_L, \check d_\chi$,
according to our  general conventions.

For a cusp form $f\in\text{Cusp}_n$ we introduce the generating functions
$$E^+_f(t) = \sum_{V\in \text{Bun}_n(X)} f(V) t^{\deg(V)} Z^+_V
 \quad \in \quad \text{Heis}[[t, t^{-1}]],
\leqno (6.4.2)$$
$$E^-_f(t) = \sum_{V\in \text{Bun}_n(X)} f(V^*) t^{-\deg(V)} Z^-_V
 \quad \in \quad \text{Heis}[[t, t^{-1}]],$$
$$\Psi^+_f(t) = \sum_{\Cal F \in \text{Coh}_{0,x}} \overline{\chi_f}
([\Cal F]) t^{h^0(\Cal F)} 
|\text{Aut}(\Cal F)| Z^+_{\Cal F} \quad \in \quad \text{Heis}[[t]],$$

$$\Psi^-_f(t) = \sum_{\Cal F \in \text{Coh}_{0,x}} {\chi_f}([\Cal F])
 t^{-h^0(\Cal F)} 
|\text{Aut}(\Cal F)| Z^-_{\Cal F} \quad \in \quad
 \text{Heis}[[t^{-1}]].$$

Similar generating functions in ${\text{Heis}}^\vee$
 will be denoted by $\check E^\pm_f(t), \check \Psi^\pm_f(t)$.
 Notice that, as in (4.5.1), we have Euler product expansions
$$\Psi^\pm_f(t) = \prod_{x\in X} \Psi^\pm_{f,x}(t^{\deg(x)}), \quad
\Psi^+_{f,x}(t) = \sum_{\mu = (\mu_1 \geq ... \geq \mu_n\geq 0)} 
\overline{\chi_f}([\Cal F_{x,\mu}])
|\text{Aut}(\Cal F_{x,\mu})| t^{|\mu|} Z^-_{\Cal F_{x,\mu}}$$
and similarly for $\Psi^-_{f,x}(t)$. We will also use the
 logarithmic generating functions
$$a_f^\pm(t) = \log \, \Psi^\pm_f(t) = \sum_{x\in X}
 a^\pm_{f,x}(t), \quad a^\pm_{f,x}(t) =
\log\, \Psi^\pm_{f,x}(t).$$
The relations between generating functions with the
 same superscript ($+$ or $-$) are the same as
 given 
in Theorem 3.3, so we concentrate on relations involving
 generating functions with different signs 
in the superscript.

\proclaim{(6.5) Theorem} (a) In the algebra Heis we
 have the following identities for any
$f\in\text{Cusp}_n$, $g\in\text{Cusp}_m$:
$$E^+_f(t_1) \Psi_g^-(t_2) = \Psi^-_g(t_2) E^+_f(t_1),
 \leqno (6.5.1)$$
$$\Psi^+_g(t_1) E^-_f(t_2) = { \text{LHom}(f,g, q^{m/2}
t_1c/t_2)\over \text{LHom}(f,g, q^{(m/2)-1}t_1c/t_2)}
 E^-_f(t_2) \Psi^+_g(t_1),\leqno (6.5.2)$$
$$\bigl[ E^+_f(t_1), E^-_g(t_2)\bigr]  = \delta_{f,g}
 \delta(t_1c/t_2) K^n \Psi^+_f(q^{-n/2}t_1), \leqno (6.5.3)$$
$$\Psi^+_f(t_1) \Psi^-_g(t_2) = { \text{LHom}(g,f, q^{(m+n/2)-1 }
t_1c/t_2)\over \text{LHom}(g,f, q^{(m+n/2}t_1c/t_2)}\Psi^-_g(t_2) 
\Psi^+_f(t_1).\leqno (6.5.4)$$

(b) In ${\text H}\text{eis}^\vee$ we have the following identities:
$$\check \Psi^-_g(t_1) \check E^+_f(t_2) = { \text{LHom}(g,f,
 q^{m/2}t_2/t_1)\over \text{LHom}(g,f, q^{(m/2)-1}t_2/t_1)}
\check E^+_f(t_2) \check \Psi^-_g(t_1), \leqno (6.5.5)$$

$$\check E^-_f(t_1) \check \Psi^+_g(t_2) = \check \Psi^+_g(t_2)
\check E^-_f(t_1),\leqno (6.5.6)$$

$$\bigl[ \check E_f^-(t_1), \check E^+_g(t_2)\bigr] = \delta_{f,g} 
\delta(t_2/t_1) K^{-n} \check \Psi^-_f(q^{n/2} t_1),\leqno (6.5.7)$$

$$ \check \Psi_g^-(t_1) \check \Psi^+_f(t_2) = { \text{LHom}
(g,f, q^{(m+n/2)-1}t_2/t_1)\over \text{LHom}(g,f,
 q^{(m+n)/2}t_2/t_1)}\check\Psi^+_f(t_2) \check\Psi^-_g(t_1).\leqno (6.5.8)$$
\endproclaim

\demo{Proof}
We will prove only part (a), part (b) being similar.
 The equality (6.5.1) follows from
Proposition 6.2.12 since the only morphism from a
 torsion sheaf to a vector bundle is the zero
morhism.

Let us prove (6.5.2). By Proposition 6.2.12, 
$$\Psi^+_g(t_1)E^-_f(t_2) = \sum_{\Cal F\in \text{Coh}_{0,X}}
 \sum_{V\in\text{Bun}_n(X)}
\overline{\chi_g}([\Cal F]) f(V^*) \cdot |\text{Aut}(\Cal F)|
\cdot \leqno (6.5.9)$$
$$\cdot t_1^{h^0(\Cal F)} t_2^{-\deg(V)} \sum_{W, 
\Cal F^{\prime\prime} \in\text{Coh}_X}
g_{\Cal F V}^{W\Cal F^{\prime\prime}} K_{\bar\Cal F -
 \bar\Cal F^{\prime\prime}}
\langle \bar\Cal F^{\prime\prime} - \bar W, \, \bar\Cal F -
 \bar\Cal F^{\prime\prime}\rangle
\cdot Z^-_W Z^+_{\Cal F^{\prime\prime}}.$$
In order that $g_{\Cal F V}^{W\Cal F^{\prime\prime}}\neq 0$, 
i.e., that there exist at least one exact sequence
$$0\rightarrow W\rightarrow V\buildrel\varphi\over\rightarrow 
\Cal F \rightarrow
\Cal F^{\prime\prime} \rightarrow 0, \leqno (6.5.10)$$
$\Cal F^{\prime\prime}$ must lie in $\text{Coh}_{0,X}$, and $W$ 
in $\text{Bun}_n(X)$. Therefore
$K_{\bar\Cal F - \bar\Cal F^{\prime\prime}}$ lies in the center
 and is the same as
$c_{\bar\Cal F - \bar\Cal F^{\prime\prime}}$. Also, we have
$$\langle \bar\Cal F^{\prime\prime} - \bar W, \, \bar\Cal F -
 \bar\Cal F^{\prime\prime}\rangle 
= \langle \bar W, \bar\Cal F - \bar\Cal F^{\prime\prime}\rangle^{-1} = 
q^{-n(h^0(\Cal F) -
h^0(\Cal F^{\prime\prime}))/2}.$$
By using (6.2.13), we now write (introducing $\Cal F'$ to be
 $\text{Im}(\varphi)$
in (6.5.10)):

$$\Psi^+_g(t_1) E^-_f(t_2) = \sum_{V,W\in\text{Bun}_n(X)}\,
 \sum_{\Cal F, \Cal F',
 \Cal F^{\prime\prime}
\in \text{Coh}_{0,X}} g_{W\Cal F'}^V g_{\Cal F'
 \Cal F^{\prime\prime}}^{\Cal F}
\cdot { |\text{Aut}(\Cal F')| \cdot |\text{Aut}(\Cal F^{\prime\prime})|
\cdot |\text{Aut}(W)|\over
|\text{Aut}(V)|}\cdot$$

$$\cdot \overline{\chi_g}([\Cal F]) f(V^*) t_1^{h^0(\Cal F)} 
t_2^{-\deg(V)} c_{\bar\Cal F'} q^{-nh^0(\Cal F)/2} Z^-_W 
Z^+_{\Cal F^{\prime\prime}}.$$

We now use the fact that $\overline{\chi_g}$ is a character of
 the Hall algebra, together
with the properties of the (dual) Hecke operators (4.3) to transform
 this sum into

$$\sum_{W, \Cal F', \Cal F^{\prime\prime}}\overline{\chi_g}([\Cal F'])
 \overline{\chi_g}([\Cal F
^{\prime\prime}]) |\text{Aut}(\Cal F')|\cdot |\text{Aut}(\Cal 
F^{\prime\prime})|\cdot (T_{\Cal F'} f)
(W^*) \cdot$$

$$\cdot (t_1c/q^{n/2})^{h^0(\Cal F')} t_1^{h^0(\Cal F^{\prime\prime})}
 t_2^{-h^0(\Cal F')}
t_2^{-\deg(W)} Z^-_W Z^+_{\Cal F^{\prime\prime}}.$$
Since $f$ is an eigenfunction of Hecke operators, this can be factored
 into the product

$$\left( \sum_{W, \Cal F^{\prime\prime}}
 \overline{\chi_g}([\Cal F^{\prime\prime}]) \cdot
|\text{Aut}(\Cal F^{\prime\prime})|\cdot f(W^*) t_1^{h^0(\Cal F^{\prime\prime})} 
t_2^{-\deg(W)} Z^-_W Z^+_{\Cal F^{\prime\prime}}\right)\times $$

$$\times \left(\sum_{\Cal F'}\overline{\chi_g}([\Cal F']) \chi_f([\Cal F']) \cdot 
|\text{Aut}(\Cal F')|\cdot \bigl( t_1c/t_2q^{n/2})\bigr)^{h^0(\Cal F')}\right) =$$

$$= E^-_f(t_2) \Psi^+_g(t_1) \cdot \Lambda(t_1c/t_2q^{n/2}),$$
where $\Lambda(t)$, given by (4.1.10), is equal to the ratio

$${ \text{LHom}\left(f,g, q^{n+m\over 2}t\right) \over \text{LHom}\left(f,g,
 q^{{n+m \over 2} -1} t\right)},$$
as we saw in (4.1.14). This proves the equality (6.5.2). 

\vskip .2cm

We now prove (6.5.3). By definition,
$$E^+_f(t_1) E^-_g(t_2) = \sum_{V\in\text{Bun}_n(X)\atop 
W\in\text{Bun}_m(X)} f(V) g(W^*)
 t_1^{\deg(V)} t_2^{-\deg(W)} Z^+_V Z^-_W =$$
$$\sum_{V,W,M,\Cal F} g_{VW}^{M\Cal F} K_{\bar W-\bar M} \cdot \langle \bar\Cal F - 
\bar M,\,
\bar W - \bar M\rangle \cdot f(V) g(W^*) t_1^{\deg(V)} 
t_2^{-\deg(W)} Z^-_M Z^+_{\Cal F},$$
where $M, \Cal F$ run over all the isomorphism classes of coherent sheaves.
 However,
 since $f$ and $g$ are cusp forms, the reasoning 
similar to that in the proof of (3.3.4)
(see (4.4)) shows that we can restrict the summation
 to $M$ equal either $W$ or 0, 
without changing the sum. The contribution from $M=W$ 
gives $E^-_g(t_2) E^+_f(t_1)$,
while computing the contribution from $M=0$ we find 
(note that $m=n$ because $M=0$):
$$ \bigl[ E^+_f(t_1),\, E^-_g(t_2)\bigr] =
 \sum_{V,W,\Cal F} g_{W\Cal F}^W 
{|\text{Aut}(\Cal F)|\over |\text{Aut}(V)|} K_{\bar W} f(V) g(W^*)
\cdot t_1^{\deg(V)} t_2^{-\deg(W)} \langle \Cal F, W\rangle Z^+_{\Cal F}$$

$$= \sum_{W, \Cal F} \bigl(T^\vee_{\Cal F} f\bigr) (W)
 \cdot {|\text{Aut}(\Cal F)|\over |\text{Aut}(V)|}
\cdot g(W^*) K^m c_{\det(W)} t_1^{\deg(W)+h^0(\Cal F)}
 t_2^{-\deg(W)} q^{-mh^0(\Cal F)/2} Z^+_{\Cal F}$$
$$=\left(\sum_{\Cal F} \chi^\vee_f ([\Cal F])\cdot
 |\text{Aut}(\Cal F)|\cdot (t_1/q^{m/2})^{h^0(\Cal F)}
Z^+_{\Cal F}\right) \cdot K^m \cdot$$
$$\cdot \left(\sum_W { f(W) g(W^*)\over |\text{Aut}(W)|}
 \biggl({t_1\over t_2}\biggr)^{\deg(W)}
c_{\det(W)} \right).$$
If $f\neq g$, then for any $d$ we have, by (2.6.11) and (2.7.4):
$$\sum_{W\in\text{Bun}_{n,d}(X)} { f(W) g(W^*)\over
 |\text{Aut}(W)|} \quad =\quad 0, \leqno (6.5.11)$$
and when $f=g$ we get the formula (6.5.3), because of
 the assumptions (2.6.14) we made on the
normalization of $f\in \text{Cusp}_n$.

\vskip .2cm

Finally, let us prove (6.5.4). Keeping $f$ and $g$ fixed,
 let
$$l(t)\quad  = \quad \log { \text{LHom}\left(g,f, 
q^{ {m+n\over 2} - 1} t\right)\over 
\text{LHom}\left(g,f, q^{ {m+n\over 2}} t\right) }
\quad = \quad \sum_{x\in X} l_x(t^{\deg(x)}),$$
where $l_x$ is the logarithm of the ratio of the 
Euler factors corresponding to $x$. By taking logarithm
of (6.5.4) and using the notation (6.4.3), we find 
that it is enough to prove, for each $x\in X$, the
equality
$$\bigl[ a^+_{f,x}(t_1),\, a^-_{g,x}(t_2)\bigr] = l_x(ct_1/t_2).$$
By introducing the Taylor coefficients of $a^\pm_{f,x}$ and $l_x$:
$$a^\pm_{f,x}(t) = \sum_{d=0}^\infty a^\pm_{f,x,d} t^{\pm d}, \quad
l_x(t) = \sum_{d=0}^\infty l_{x,d} t^d, \leqno (6.5.12)$$
we can rewrite the desired equality as
$$\bigl[ a^+_{f,x,d}, \, a^-_{g,x,d'}\bigr] =
 \delta_{dd'} c^d l_{x,d}.\leqno (6.5.13)$$
This is a statement about the Heisenberg double
 of the Hopf subalgebra $\Cal H_x \i B(\text{Coh}_X)$
generated by $[\Cal F],\, \Cal F\in \text{Coh}_{0,x}(X)$
 and by the element $c_x$. The function $\psi_{f,x}(t)$ 
taking values in this subalgebra, has the coproduct
$$\Delta (\psi_{f,x}(t)) = \psi_{f,x}(t\otimes c_x)
 (1\otimes \psi_{f,x}(t)).$$
This is proved in the same way as (3.3.3). 
So for $a_{f,x}(t) = \log(\psi_{f,x}(t))$ we
have:
$$\Delta (a_{f,x}(t)) = a_{f,x}(t\otimes c_x) + 1\otimes a_{f,x}(t),$$
or, in the coefficient language,
$$\Delta (a_{f,x,d}) = a_{f,x,d} \otimes c_x^d \, + \, 1\otimes a_{f,x,d}.$$
So let us look at the following general situation.
 Let $V$ be a $\Z$-graded vector space. 
Consider the commutative Hopf algebra $\Xi = S^\bullet(V)
 \otimes \C[c, c^{-1}]$
with comultiplication
$$\Delta(v) = v\otimes c^{\deg(v)} \, + \, 1\otimes v,
 \quad \Delta(c) = c\otimes c.$$
For $v\in V$ let $Z^+_v$ be the corresponding element of
 $HD(\Xi)$. For $\phi\in V^*$ let $Z^-_\phi$
be the element of $HD(\Xi)$ corresponding to the
 linear form $\Xi\rightarrow \C$
which on each $S^\bullet(V)\otimes c^m$ is the derivation
 $\partial/\partial \phi$ in the symmetric
algebra. let also $\check Z^+_v, \check Z^-_\phi$ be  the
 similar elements of of $\check HD(\Xi)$.

\proclaim {(6.5.14) Lemma} In the described situation
 we have, for homogeneous $v\in V, \phi\in V^*$,
$$[Z^+_v, Z^-_\phi] = \phi(v) c^{\deg(v)}, \quad
 [\check Z^-_\phi, \check Z^+_v] = \phi(v).$$
\endproclaim

The first formula in this lemma implies the equality 
(6.5.13) (we have added the second formula to give a hint
 of the proof of (6.5.8)). Indeed, take $V$ to be the space
 of primitive elements of the commutative and cocommutative
 Hopf algebra $\Cal H_x/(c-x-1)$. It is graded by $\deg[\Cal F]
 = h^0(\Cal F)/\deg(x)$. 
The algebra $\Xi$ defined above is just $\Cal H_x$ itself.
 Taking $v=a_{f,x,d}$, we find that
$Z^+_v = a^+_{f,x,d}$, while $a^-_{f,x,d} = Z^-_\phi$
 where $\phi$ is the linear functional on $V$ given by
$$\phi(w) = \sum_{\Cal F\in\text{Coh}_{0,x}(X)}
{ a_{f,x,d}(\Cal F) w(\Cal F)\over |\text{Aut}(\Cal F)|}
= (w, \overline{a_{f,x,d}}).$$
To deduce (6.5.13), it remains to recall the formula (4.5.7) 
for the scalar product of the $a_{f,x}(t)$.

\vskip .2cm

As for Lemma 6.5.14, it immediately reduces to the particular
 case when $\dim(V)=1$, i.e., 
$\Xi = \C[x, c, c^{-1}]$ with $\Delta(x) = 
x\otimes c^d \,+ \,1\otimes x$ and
$\Delta(c)=c\otimes c$. The treatment of this case is elementary
 and is left to the reader.

Theorem 6.5 is proved.

\enddemo

\vskip .2cm

\noindent{\bf  (6.6) The Drinfeld double of the algebra
 of automorphic forms.} Let $U=U(\Cal A)$
be the restricted Drinfeld double of the Hopf algebra 
$B(\Cal A)$, $\Cal A =\text{Coh}_X$, see (6.3). To make
 the formulas below more symmetric, we extend $U$ by adding
 square roots of the central elements $c_L, L\in\text{Pic}(X)$.
 This means that we choose an identification $\text{Pic}(X)\simeq
 \Z \oplus\bigoplus Z/m_i$ and embed $\text{Pic}(X)$ into
 the group $\sqrt{\text{Pic}(X))} = 
{1\over 2} \Z \oplus\bigoplus Z/2m_i$ in an obvious way.
 We set
$$U(c^{1/2}) = U \otimes_{\C[\text{Pic}(X)]} \C\left[
 \sqrt{\text{Pic}(X))}\right].$$
Thus, for any $L\in \text{Pic}(X)$ there is a central 
element $c_L^{1/2}\in U(c^{1/2})$. 

For a cusp for $f\in\text{Cusp}_n$ we introduce the following 
generating functions with coefficients
in $U(c^{1/2})$:
$$Y_f^+(t) = \sum_{V\in\text{Bun}_n(X)} f(V) t^{\deg(V)} W^+_V,
 \quad Y^-_f(t) = 
\sum_{V\in\text{Bun}_n(X)} f(V^*) t^{-\deg(V)} W^-_V,
 \leqno (6.6.1)$$
$$\Phi^+(t) = \sum_{\Cal F\in\text{Coh}_{0,X}}
 \overline{\chi_f}([\Cal F]) t^{h^0(\Cal F)} \cdot 
|\text{Aut}(\Cal F)|\cdot c_{\bar\Cal F}^{1/2} W^+_{\Cal F}, $$
$$\Phi^-(t) = \sum_{\Cal F\in\text{Coh}_{0,X}} {\chi_f}([\Cal F])
 t^{-h^0(\Cal F)} \cdot 
|\text{Aut}(\Cal F)|\cdot c_{\bar\Cal F}^{1/2} W^-_{\Cal F}.$$

Now we can formulate the main result of this paper. 

\proclaim{(6.7) Theorem} The above generating functions 
satisfy the following identities,
for any $f\in\text{Cusp}_n$, $g\in\text{Cusp}_m$:
$$Y^\pm_f(t_1) Y^\pm_g(t_2) = q^{mn(1-g_X)} {\text{LHom}(f,g, t_2/t_1)\over
\text{LHom}(f,g, t_2/qt_1)} Y^\pm_g(t_2) Y^\pm_f(t_1), \leqno (6.7.1)$$

$$\bigl[ \Phi^\pm_f(t_1), \Phi^\pm_g(t_2)\bigr] = 0, \leqno (6.7.2)$$

$$Y^\pm_f(t_1) \Phi^+(t_2) = \left( {\text{LHom}\left(f,g,
 q^{m\over 2} c^{\mp 1/2}t_2/t_1\right)\over
\text{LHom}\left(f,g, q^{{m\over 2}-1} c^{\mp 1/2}
 t_2/t_1\right)}\right) ^{\pm 1} 
\Phi^+_g(t_2) Y^\pm_f(t_1), \leqno (6.7.3)$$

$$Y^\pm_f(t_1) \Phi^-(t_2) = \left( {\text{LHom}
\left(g,f , q^{m\over 2} c^{\pm 1/2}t_1/t_2\right)\over
\text{LHom}\left(g,f, q^{{m\over 2}-1} c^{\pm 1/2} t_1/t_2\right)}
\right) ^{\mp 1} 
\Phi^-_g(t_2) Y^\pm_f(t_1), \leqno (6.7.4)$$

$$\bigl[ Y^+_f(t_1), Y^-_g(t_2)\bigr] = \delta_{f,g} 
\left\{ \delta\left({ct_1\over t_2}\right) K^n \Phi^+_f(q^{-n/2}c^{1/2}t_1)
 - \delta\left({t_2\over t_1}\right) K^{-n} \Phi_f^-(q^{n/2} c^{-1/2}t_2)
\right\},\leqno (6.7.5)$$

$$\Phi^+_f(t_1) \Phi^-_g(t_2) = { \text{LHom}\left(g,f, 
q^{{m+n\over 2}-1} c^2 t_1/t_2\right)
\text{LHom}\left(g,f, q^{{m+n\over 2}} c t_1/t_2\right)\over
\text{LHom}\left(g,f, q^{{m+n\over 2}} c^2 t_1/t_2\right)\text{LHom}
\left(g,f, q^{{m+n\over 2}-1} c t_1/t_2\right)} \Phi^-_g(t_2)
 \Phi^+_f(t_1), \leqno (6.7.6)$$

$$KY^\pm_f(t) = q^{\pm n(g_X-1)} Y^\pm_f(t) K, \quad K\Phi^\pm_f(t) 
= \Phi^\pm_f(t) K.\leqno (6.7.7)$$

\endproclaim

\demo{Proof} We use the embedding $\varkappa:
 U(\Cal A)\hookrightarrow \text{Heis}(\Cal A)
\otimes { \text{Heis}}^\vee (\Cal A)$ given by 
  Proposition 6.3.12 (b). From Theorem 3.5 we find that
 on the generating functions the embedding is  as follows:
$$Y^+_f(t) \mapsto E^+_f(t\otimes \check c) \bigl(1\otimes
 \check K^n \check\Psi^+_f(q^{-n/2}t)\bigr)
\,\, + \,\, 1\otimes \check E^+_f(t),$$

$$Y^-_f(t) \mapsto E^-_f(t)\otimes 1 \,\, + \,\,
K^{-n} \Psi^-_f(q^{-n/2}t) \otimes \check E^-_f(t),$$

$$\Phi^+_f(t) \mapsto \Psi^+_f(c^{1/2}t \otimes \check c^{3/2})
 \check \Psi^+_f(c^{1/2}\otimes \check c^{1/2}t),$$

$$\Phi^-_f(t) \mapsto \Psi^-_f(c^{-1/2}t \otimes \check c^{-1/2})
 \check \Psi^-_f(c^{-1/2}\otimes \check c^{-1/2}t),$$

$$c_L \mapsto c_L\otimes \check c_L, \quad K\mapsto K\otimes \check K.$$

Our result follows from these formulas and from Theorem 6.5 giving relations
in the Heisenberg doubles. \enddemo

\vskip .2cm

\noindent{\bf  (6.5) Final remarks.} Comparing (6.7) with the
 relations (5.1) describing quantum affine algebras, we see
 that there is indeed  an almost complete similarity. The only
 difference worth mentioning is related to the zero-modes of
 the Cartan generators. In quantum affine algebras there
are as many such zero modes $\varphi^+_i(0) = \varphi_i^-(0)^{-1}$
 as there are simple roots. In the automorphic case there is only 
one such generator $K$, coming from the K-theory of the curve $X$ 
(the generating functions $\Phi^\pm_f(t)$ have constant term 1).
 This suggests that it may be more natural to consider, as, e.g., 
 in [MW], the space $\Pi$ of all cusp eigenforms on $\text{Bun}(X)$, 
i..e., of functions
$$f\cdot t^{\deg}, \quad V\mapsto f(V) t^{\deg(V)},
 \quad f\in \text{Cusp}, t\in \C^*$$
and make it into a 1-dimensional complex manifold 
(disjoint union of copies of $\C^*$
 paramterized by the set Cusp). Then instead of writing
 $E^\pm_f(t), \Phi^\pm_f(t)$, we can write $E^\pm(\pi),
 \Phi^\pm(\pi)$ where $\pi = f\cdot t^{\deg}\in\Pi$. So
 the whole algebra will look more like an $sl_2$-current
 algebra but with currents being defined on the space $\Pi$. We have 
therefore a kind of conformal field theory on the space of
 cusp forms, which strengthens certain
analogies from [Kap].

 \vskip .2cm

(b) As with the generating functions $\Psi^\pm_f(t)$ 
of Heis, the $\Phi^\pm_f(t)$
can be written as Euler products $\prod_{x\in X} 
\Phi_{f,x}^\pm(t)$. In the same way as we
proved (6.7.6), one obtains that for different $x$ 
the $\Phi^\pm_{f,x}(t)$ commute with each
other, while $$\Phi^+_{f,x}(t_1) \Phi^-_{g,x}(t_2)
\Phi^+_{f,x}(t_1) ^{-1}\Phi^-_{g,x}(t_2)^{-1}$$
is the ratio similar to that in (6.7.6) but formed
 by the Euler factors at $x$ of the
$L$-functions appearing there. This means that for
 any $x,f$ the coefficients $\bold a^\pm_{f,x,d}$
of the expansion of $\bold a_{f,x}^\pm(t) := \log 
(\Phi^\pm_{f,x}(t))$ form, apart from normalization,
free bosons:
$$\bigl[ \bold a^+_{f,x,d}, \, \bold a^-_{g,x,d'} \bigr] =
 \delta_{d+d',0} \cdot\left(\sum_{i,j} 
\biggl({\lambda_{x,j}(g)
\over \lambda_{x,i}(f)}\biggr)^d\right)\cdot 
q^{d(m+n-2)/2} c_x^d (q^d-1) (c_x^d - 1).$$
According to the point of view going back to
 Y.I. Manin and B. Mazur, one should visualize any
1-dimensional arithmetic scheme $X$ as a kind
 of 3-manifold and closed points  $x\in X$
as oriented circles in this 3-manifold. Thus the
 Frobenius element (which is only a conjugacy class
in the fundamental group) is visualized as the
 monodromy around the circle (which, as an element
of the fundamental group, is also defined only
 up to conjugacy since no base point is chosen
on the circle),  Legendre symbols  as linking
 numbers and so on.
From this point of view, it is natural to think
 of the operators (algebra elements)
$\bold a_{f,x,d}$ for fixed $f$ and varying $x,d$
 as forming a free boson field $A_f$
on the ``3-manifold" $X$; more precisely, for $\pm d > 0$,
 the operator
$\bold a^\pm_{f,x,d}$ is the $d$th Fourier component
of $A_f$ along the ``circle" $\text{Spec}(\F_q(x))$. 
The bosons $\bold a^\pm_{f,x,d}$ and their
sums over $x\in X$ (i.e., the Taylor components of
 $\log \,\Phi^\pm_f(t)$) will
be used in a subsequent paper to construct
 representations of $U$ in the spirit of [FJ].

 \vfill\eject

\centerline {\bf  References}

\vskip 1cm

[AF] A. Yu. Alexeev, L.D. Faddeev, Quantum $T^*G$ as 
a toy model for conformal field theory, {\it Comm. math. Phys.}, 
{\bf  141} (1991), 413-443.

\vskip .3cm

[Be] A. Beilinson, Coherent sheaves on $P^n$ and problems of linear
algebra, {\it Funct. Anal. Appl.} {\bf  12} {(1978)}, N. 3, p. 68-69.

\vskip .3cm

[BerZ] I. Bernstein, A. Zelevinsky, Induced representations
of reductive $p$-adic groups I, {\it Ann. ENS}, {\bf  10} (1977), 441-472.

\vskip .3cm

[Bo] A.I. Bondal, Representations of associative algebras
 and coherent sheaves, {\it  Math. USSR Izvestija}, {\bf  34} (1990), 23-42.

\vskip .3cm

[Bu] D. Bump, The Rankin-Selberg method: a survey, in:
``Number theory, trace formulas and discrete groups"
(Symposium in honor of A. Selberg) p. 49-109, Academic Press, 1989.

\vskip .3cm

[CP] V. Chari, A. Pressley, A Guide to Quantum groups,
Cambridge Univ. Press, 1995. 

\vskip .3cm

[De] P. Deligne, Constantes des equations fonctionelles des fonctions L, in: 
Lecture Notes in Math. {\bf  349} , p. 501-597,
Springer-Verlag, 1973. 

\vskip .3cm

[DF] J. Ding, I. Frenkel, Isomorphism of two realizations of quantum
affine algebra $U_q(\widehat{gl(n)})$, {\it Comm. Math. Phys.} {\bf  156}
(1993), 277-300. 

\vskip .3cm

[Dr 1] V.G.  Drinfeld, A new realization of Yangians and quantized enveloping
algebras, {\it Sov. Math. Dokl.} {\bf  36} (1988), 212-216.

\vskip .3cm

[Dr 2] V.G.  Drinfeld, A new realization of Yangians 
and quantized enveloping
algebras (in Russian), Preprint FTINT, Kharkov, 1986, \# 30.

\vskip .3cm

[Dr 3] V.G. Drinfeld, Quantum groups, Proc. ICM-86 
(Berkeley), vol.1, p. 798-820, Amer. math. Soc.,
1987.

\vskip .3cm

[FJ] I. Frenkel, N. Jing, Vertex representations of
 quantum affine algebras, {\it Proc. Nat. Acad. USA} {\bf  85}
(1988), 9373-9377. 

\vskip .3cm

[FK] Y. Flicker, D. Kazhdan, Geometric Ramanujan conjecture
and Drinfeld reciprocity law, in: ``Number theory,
 trace formulas and discrete groups"
(Symposium in honor of A. Selberg) p. 201-218, 
Academic Press, 1989.

\vskip .3cm

[FLM] I. Frenkel, J. Lepowsky, A. Meurman, Vertex Operators and the Monster,
Academic Press, 1986. 

\vskip .3cm

[Fr] P. Freyd, Abelian Categories, Harper \& Row Publ. New York, 1964.

\vskip .3cm

[GL] W. Geigle, H. Lenzing, A class of weighted projective
 curves arising in representation theory of
 finite-dimensional algebras, Lecture Notes in Math.
{\bf  1273}, p.265-297, Springer-Verlag, 1988. 

\vskip .3cm

[GN] I.M. Gelfand, M.A. Najmark, The principal series
 of irreducible representations of the complex unimodular group,
{\it Dokl. AN SSSR}, {\bf  56} (1947), 3-4 (reprinted in
Collected Papers of I.M. Gelfand, vol 2, p. 128-129,
 Springer-Verlag, 1988). 

\vskip .3cm

[GKV] V. Ginzburg, M. Kapranov, E. Vasserot, Langlands 
reciprocity for
algebraic surfaces, {\it Int. Math. Research Notes}, 
{\bf   2} (1995), 147-160.

\vskip .3cm

[Gr] J.A. Green, Hall algebras, hereditary algebras 
and quantum groups,
{\it Invent. Math.} {\bf  120} (1995), 361-377. 

\vskip .3cm

[Gro] I. Grojnowski, Affinizing quantum algebras: from D-modules to 
K-theory, Preprint 1994. 

\vskip .3cm

[Ha 1] G. Harder, Chevalley groups over function fields
and automorphic forms,
 {\it Ann. Math.} {\bf  100}
(1974), 249-300. 

\vskip .3cm

[Ha 2] G. Harder, Minkowskische Reduktiontheorie
 \"uber Funktionenk\"orper, {\it Invent. math.} {\bf  7} (1968), 33-54.

\vskip .3cm

[HN] G. Harder, M.S. Narasimhan, On the cohomology groups of moduli spaces
of vector bundles on curves, {\it Math. Ann.} {\bf  212} (1975), 215-248.

\vskip .3cm

[JPS] H. Jacquet, I.I. Piatetski-Shapiro, J. Shalika, 
Rankin-Selberg convolutions, {\it Amer. J. Math.} {\bf  105}
(1983), 367-464.

\vskip .3cm

[Kac] V. Kac, Infinite-dimensional Lie algebras, Cambridge University
Press, 1992. 

\vskip .3cm

[Kap] M. Kapranov, Analogies between the Langlands
 correspondence and
the topological quantum field theory, in: ``Functional analysis
 on the eve of 21st century" (in honor of I.M. Gelfand), vol 1,
 p. 119-151, Birkh\"auser, Boston, 1995.

\vskip .3cm

[Kas] R.M. Kashaev, Heisenberg double and the pentagon relation, 
preprint q-alg/9503005.

\vskip .3cm

[La 1] R.P. Langlands, On the functional equations satisfied by Eisenstein
series, Lecture Notes in Math. {\bf  544}, Springer-Verlag, 1976.

\vskip .3cm

[La 2] R.P. Langlands, Euler Products, Yale Univ. Press, 1971.

\vskip .3cm

[La 3] R.P. Langlands, Eisenstein series, in: ``Algebraic
groups and discontinuous subgroups" (Proc. Symp. Pure math. 
{\bf  9}), p. 235-252, Amer. Math. Soc., 1966.

\vskip .3cm

[Lu1] G. Lusztig, Introduction to Quantum Groups
(Progress in Math. 
{\bf  110}), Birkhauser, Boston, 1993. 

\vskip .3cm

[Lu2] G. Lusztig, Canonical bases arising from quantized enveloping 
algebras, {\it J. AMS}, {\bf  3} (1990), 447-498.

\vskip .3cm

[Lu3] G. Lusztig, Quivers, perverse sheaves and quantized 
enveloping algebras, {\it J. AMS}, {\bf  4} (1991), 365-421.

\vskip .3cm

[Lu4] G. Lusztig, Character sheaves I, II, {\it Adv. in Math.}
 {\bf  56} (1985), 193-237, {\bf  57} (1985),
226-265.

\vskip .3cm

[Mac] I.G. Macdonald, Symmetric Functions and Hall
 Polynomials, Oxford, Clarendon Press, 1995.

\vskip .3cm

[MW] C. Moeglin, J.-L. Waldspurger, Spectral Decomposition
and Eisenstein Series, Cambridge Univ. Press, 1995.

\vskip .3cm

[Mor] L.E. Morris, Eisenstein series for reductive
 groups over global function fields I, II, {\it Canad. J. Math.} 
{\bf  34} (1982), 91-168, 1112-1282.

\vskip .3cm

[R1] C.M. Ringel, Hall algebras and quantum groups, 
{\it Invent. Math.}
{\bf  101} (1990), 583-592.

\vskip .3cm

[R2] C.M. Ringel, The composition algebra of a cyclic quiver,
{\it Proc. London Math. Soc.} (3) {\bf  66} (1993), 507-537.

\vskip .3cm

[R3] C.M. Ringel, Hall algebras revisited, in:
 ``Israel Math. Conf. Proc." vol. 7 (1993), p. 171-176. 

\vskip .3cm

[Sh] J. Shalika, Multiplicity one theorem for GL(n),
{\it Ann. Math.}  {\bf  100} (1974), 171-193.

\vskip .3cm

[ST] M.A. Semenov-Tian-Shansky, Poisson Lie groups,
quantum duality principle and the quantum double,
{\it Contemporary math.}, {\bf  178}, p. 219-248,
Amer. Math. Soc, 1994.

\vskip .3cm

[Stu] U. Stuhler, On the cohomology of $SL_n$ over rings of algebraic functions,
in: Lecture Notes in Math. {\bf  967}, p. 316-359, Springer-Verlag, 1982.

\vskip .3cm

[Var] A. Varchenko, Multidimensional Hypergeometric Functions and Representation 
Theory of Lie Algebras and
Quantum Groups, World Scientific, Singapore, 1995.

\vskip .3cm

[X1] J. Xiao, Hall algebra in a root category, Preprint 95-070, Univ. of Bielefeld,
1995.

\vskip .3cm

[X2] J. Xiao, Drinfeld double and Green-Ringel theory of Hall algebras, Preprint
95-071, Univ. of Bielefeld, 1995.

\vskip .3cm

[Ze] A. Zelevinsky, Representations of Finite Classical
Groups (a Hopf Algebra Approach), Lecture Notes in Math.
{\bf  869}, Springer-Verlag, 1981.
                                                                                                                             \vskip 2cm

Department  of Mathematics, Northwestern University, Evanston IL 60208, USA
(kapranov\@chow.math.nwu.edu).

\vskip .4cm

Address till December 1996:  Max-Planck Institut
 f\"ur Mathematik, Gottfried-Claren-Strasse 26,
53225 Bonn, Germany (kapranov\@mpim-bonn.mpg.de)






\enddocument